\documentclass[10pt,a4paper
]{amsart}

\usepackage{graphicx}
\usepackage{psfrag}
\newtheorem{theorem}{Theorem}[section]
\newtheorem*{T}{Theorem}
\newtheorem{lemma}[theorem]{Lemma}

\newtheorem{proposition}[theorem]{Proposition}
\newtheorem*{Pro}{Proposition}

\theoremstyle{definition}

\theoremstyle{remark}

\numberwithin{equation}{section}

\newcommand{\Ce}{\mathbb{C}}

\newcommand{\Er}{\mathbb{R}}

\newcommand{\f}{^\flat}
\newcommand{\rf}{r\f}
\newcommand{\sigf}{\sigma\f}
\newcommand{\thf}{\theta\f}
\newcommand{\smf}{s\f}

\newcommand{\nf}{\nu\f}
\newcommand{\nuf}{\nf}
\newcommand{\ef}{\varepsilon\f}
\newcommand{\co}{\mathcal{O}(1)}
\newcommand{\Uf}{U\f}
\newcommand{\Af}{A\f}
\newcommand{\Bf}{B\f}
\newcommand{\lf}{\lambda\f}
\newcommand{\Lf}{L\f}
\newcommand{\pf}{P\f}
\newcommand{\ppf}{P_\perp\f}
\newcommand{\qf}{q\f}
\newcommand{\Qf}{Q\f}
\newcommand{\Qpf}{Q_\perp\f}
\newcommand{\caf}{\cala\f}

\newcommand{\cala}{\mathcal{A}}
\newcommand{\calb}{\mathcal{B}}
\newcommand{\calc}{\mathcal{C}}
\newcommand{\cald}{\mathcal{D}}

\newcommand{\calg}{\mathcal{G}}
\newcommand{\calh}{\mathcal{H}}
\newcommand{\calj}{\mathcal{J}}
\newcommand{\call}{\mathcal{L}}

\newcommand{\calp}{\mathcal{P}}

\newcommand{\calr}{\mathcal{R}}
\newcommand{\calu}{\mathcal{U}}

\newcommand{\cgf}{\calg\f}
\newcommand{\lu}{\mu_2}

\newcommand{\ve}{\varepsilon}

\newcommand{\norm}[2][{}]{\lVert#2\rVert_{{#1}}}
\newcommand{\vnorm}[2][{}]{\left\Vert#2\right\Vert_{#1}}
\newcommand{\abs}[2][{}]{\lvert#2\rvert_{#1}}
\newcommand{\vabs}[2][{}]{\left\vert#2\right\vert_{#1}}
\newcommand{\one}{\boldsymbol{1}}
\newcommand{\im}{\operatorname{Im}}

\newcommand{\sgn}{\operatorname{sgn}}

\begin{document}

\title[Thual-Fauve pulse\ \today]{The Thual-Fauve
Pulse: skew stabilization}

\author{Piero de Mottoni}
\address{Universit\`a di Roma, Tor Vergata, 21/3/1943 --- 25/11/1990}

\author{Michelle Schatzman}
\address{UMR 5585 CNRS Analyse Num\'erique \\
Universit\'e Lyon 1 \\
69622 Vil\-leur\-ban\-ne Cedex \\
France}
\email{schatz@maply.univ-lyon1.fr}
\urladdr{http://numerix.univ-lyon1.fr/$\sim$schatz}


\begin{abstract}Consider the quintic complex Ginzburg-Landau equation
\begin{equation*}
\begin{split}
u_t =&  (m+i\alpha\mu_0) u_{xx} -(m+i\alpha\mu_1) u  \\
&+ (1+i\alpha\mu_2)|u|^2u -(1+i\alpha\mu_3)|u|^4u, x\in \Er. 
\end{split}
\end{equation*}
The parameter $m$ is close to $3/16$ so that for $\alpha=0$, this equation
possesses a unstable pulse-like solution. For $\vabs{\alpha}$ small
the equation possesses pulse-like solutions of the form $e^{i\omega
t} e^{i\phi(x)} 
r(x)$, with $r$ a positive function decreasing exponentially at
infinity and $\phi$ asymptotic to $-C|x|+D$ at infinity.

 These solutions are
linearly unstable for $\abs{\alpha}\le \alpha_c$; if a certain rational
function of $\mu_2$ and $\mu_3$ is strictly positive
and not too large, they become stable for $\vabs{\alpha}\ge \alpha_c$:
when the initial data is a pulse 
plus a small perturbation, its limit for large times is the same
pulse, possibly translated in space and in phase.

This article gives a rigorous proof of a conjecture of Thual and Fauve
\cite{ThuFau}; it relies on a very detailed asymptotic analysis of the
eigenvalues of the linearized operator, depending on the parameters
$3/16 -m$ and $\mu_j$.
\end{abstract}

\maketitle

\section{Introduction}\label{sec:Introduction}

In \cite{ThuFau}, Thual and Fauve proposed a model of localized
structures generated by subcritical instabilities; in their article,
they mentioned several examples of such localized structures in
systems far from equilibrium: local regions of turbulent motion
surrounded by laminar flow as in \cite{Tr} chapter 19, spatially localized
traveling waves at 
convection onset in binary fluid mixtures as in  \cite{MoFiSt} or \cite{HeAhCa}, a 
Faraday experiment in a narrow annular dish as in \cite{KoBeSu}.

The equation proposed by Thual and Fauve is the quintic
Ginzburg-Landau equation
\begin{equation}
\frac{\partial u }{ \partial t} = m_0 \frac{\partial ^2u}{
\partial x^2} + m_1 u + m_2 |u|^2u + m_3 |u|^4 u,\label{1.1}
\end{equation}
where the $m_i$ are complex
coefficients with $\Re m_0>0$.

Thual and Fauve assumed $\Re m_3 <0$, in order to stabilize large
amplitudes, and $\Re m_1<0$, so that the zero solution should be stable.
Choosing $\Re m_2 >0$ and adequate relations on the other
parameters of the equation ensured that there would exist non
zero stable homogeneous solutions of the form $r e^{i\omega t + kx}$.

They integrated numerically~\eqref{1.1} and found thus a
solution of~\eqref{1.1} of the form
\begin{equation}
u(x,t) = e^{i\omega t} r(x) e^{i\phi(x)},\label{1.2}
\end{equation}
where $r$ takes positive values and decays exponentially at infinity,
while the phase $\phi$ is asymptotic to $-C|x| + D$ at $|x|=\infty$. 
They obtained numerically standing wave solutions of
~\eqref{1.3}, which they found to be experimentally stable
for large enough absolute values of
the imaginary parts of $m_1$, $m_2$ and $m_3$.

Let us first simplify~\eqref{1.1} thanks to some scale consideration:
by changing the units of $t$, $u$ and $x$, we can see that it is
possible to choose $\Re m_2=1$, $\Re m_3=-1$ and $m=\Re m_0=-\Re
m_1$. Thanks to the $S^1$ equivariance, we can also take $m_0$ to be
real, with very little loss of generality.
For reasons which will be justified in the course of the
article, our results hold when $m$ is slightly less than $3/16$; let
it be said only at this point that for $m=3/16$, and all the $m_i$'s
real, there exist four real heteroclinic solutions of
\begin{equation*}
-mr''+mr -r^3 + r^5=0,
\end{equation*}
which are distinct up to translation. We let $\tilde r$ be one of
these heteroclinic solutions which takes the value $\sqrt 3/2$ at
$-\infty$ and $0$ at $+\infty$; the three other ones are obtained by
mirroring $\tilde r$, $x$ or both.

For $m$ slightly less than
$3/16$, there are two homoclinic solutions with a very large
``shelf'', i.e. a region where the solution is very close to $\pm
3/4$; it is the presence of this large shelf which makes this choice
of parameters interesting.

For real values of $m_i$, pulses must be linearly unstable. Let us
sketch the argument which proves this statement: if $u$ is a
pulse-like solution according to the above definition, it will be
proved at lemma~\ref{Lemma2.1} that $u$ can be taken real, even and
positive; then $u$ solves the following differential equation:
\begin{equation}
-m u_{xx} +m u -u^3 +  u^5=0,\quad x\in \Er.\label{1.3}
\end{equation}
Moreover, the linearized operator at $u$ is the unbounded
operator in $L^2(\Er)$ given by
\begin{equation}
D(A) = H^2(\Er),\quad As = - m s'' +  m s - 3 u^2 s + 5 u^4s.\label{1.4}
\end{equation}
The operator $A$ can also be seen as a Schr\"odinger operator in
$\Er$, with potential $m  -3 u^2 + 5 u^4$; in particular, its
essential spectrum lies above $m$.
Differentiating~\eqref{1.3} with respect to $x$, we find
\begin{equation*}
Au'=0,
\end{equation*}
which expresses the translation invariance of~\eqref{1.3}.
Thus $u'$ is an eigenfunction relative to the eigenvalue $0$; since
$u'$ changes sign, the maximum principle implies that $0$ cannot be
the lower bound of the 
spectrum of $A$. We denote by $\lambda<0$ this lower bound.

The eigenmode pertaining to $\lambda$ is the shrinking-swelling mode: under most
perturbations, the pulse shrinks to $0$ or swells to infinity. 

In order to prove our main results, we introduce a few notations: 
\begin{align}
\nu &= 1-\frac{16m}{ 3},\label{1.6}\\ \quad L &= \frac{1}{
4}\ln\frac{4}{\nu}.\label{1.7}
\end{align}

In the first part of this introduction, we restrict ourselves to a
simplified case of~\eqref{1.1}: 
\begin{equation}
u_t =-m u + mu_{xx} +(1+i\alpha)|u|^2u
-|u|^4u,\label{1.5}
\end{equation}
where $m$ and $\alpha$ are real parameters;

In this article, we prove the following existence theorem:

\begin{T}[Theorem \ref{Theorem4.5}]
For all $p>0$, for all $\nu$ small enough, 
there exists $C>0$ such that for all $\alpha$ satisfying
\begin{equation*}
|\alpha|\le \alpha_m=\frac{1}{2}\sqrt\nu\,\bigl(1
-CL^{-p}\bigr),
\end{equation*}
there exists a pulse solution of~\eqref{1.5}, i.e. a solution of the form
\begin{equation*}
u(x,t)=e^{i\omega t} r(x) e^{i\phi(x)},
\end{equation*}
where $r$ is a positive function which decays exponentially at
infinity, and $\phi(x)$ is asymptotic at infinity to $-C|x|+D$,
with $C$ a positive constant and $D$ a real constant.
\end{T}

Then, we prove the stabilization property conjectured by Thual and
Fauve:

\begin{Pro}[Proposition~\ref{Proposition5.9}]
Let $u$ be the solution defined at theorem 4.5. Let
$\vabs{\alpha}\le \alpha_m$ be the above defined number; there exists
a number $\alpha_c$ whose asymptotic is
\begin{equation*}
\alpha_c(\nu) = \frac{1}{2}\sqrt\nu\left(1 - \frac{\pi^2}{48 L^2}+
O\bigl(L^{-5/2}\bigr)\right) . 
\end{equation*}
such that the solution $u$ is stable iff $\alpha_c < \abs{\alpha} <\alpha_m$, and
unstable if $\abs{\alpha}\le \alpha_c$.

More precisely, if $\abs{\alpha}<\alpha_c$, the spectrum of the
linearized operator at $u$ contains exactly 
one negative eigenvalue; when 
$\alpha=\alpha_c$, the eigenvalue $0$ is of algebraic multiplicity $3$
and geometric multiplicity $2$, with a non trivial Jordan block of
dimension $2$; when $\alpha_c < \abs{\alpha} <\alpha_m$, the linearized
operator at $u$ has the semisimple double eigenvalue $0$, and the
remainder of the spectrum is bounded away from the imaginary axis.
\end{Pro}

The general case i.e.
\begin{equation}
u_t = m u_{xx} -(m +i\alpha\mu_1) u +(1+i\alpha\mu_2)\vabs{u}^2 u
-(1+i\alpha\mu_3)\vabs{u}^4 u.\label{eq:46}
\end{equation}
is treated in section~\ref{sec:general-situation}
with very few analytical details; under the condition
\begin{equation*}
\chi(\alpha)=\biggl[\mu_2 -\frac{9\mu_3}{8}
\biggr]\biggl[\frac{\pi^2\mu_2}{4} 
-\frac{3\pi^2\mu_3}{16}
+\frac{9\mu_3}{16}\biggr]\frac{1} {\bigl(2\mu_2
-15\mu_3/8\bigr)^2}>0,
\end{equation*}
skew stabilization also takes place, as is shown in
Proposition~\ref{thr:1}. 

Let us give a very rough idea of the reason for existence, and of the skew
stabilization mechanism. In both cases, we will restrict ourselves to
the simplified equation~\eqref{1.5}.

For the existence, we start from so-called kinks: they are solutions
of~\eqref{1.5} of the form 
\begin{equation}
u(x,t) = e^{i\omega t} K(x - ct),\label{1.8}
\end{equation}
with
\begin{equation}
K(x) = r(x) e^{ikx}.\label{1.9}
\end{equation}
to make things precise, 
we demand that $r$ take positive values, increasing from $r(-\infty) = 0$ to
$\bar r=r(+\infty)>0$.                  
If we substitute~\eqref{1.8} and~\eqref{1.9} into~\eqref{1.5}, we
obtain the equation 
\begin{equation}
i\omega r - c(r' + ikr) -m(r'' + 2 i k r' - k^2 r) +m r -(1+i\alpha)r^3 +
r^5=0.\label{1.10}
\end{equation}
Define $\tilde \omega = \omega - kc$; the vanishing of the imaginary
part of~\eqref{1.10} implies 
that 
\begin{equation}
r' = \frac{(\tilde\omega -\alpha r^2)r}{2k}.\label{1.11}
\end{equation}
Since we assumed that $r$ is increasing from $0$ to $\bar r$ and is 
strictly positive, relation~\eqref{1.11} implies the sign conditions:
\begin{equation}
\tilde \omega/k >0, \quad \alpha /k>0.\label{eq:43}
\end{equation}

If we differentiate this relation with respect to $x$, we obtain 
\begin{equation}
r'' = \frac{(\tilde\omega -\alpha r^2)(\tilde\omega - 3 \alpha r^2)r}{
4k^2}.\label{eq:56}
\end{equation}
Substituting the expressions~\eqref{1.11} and~\eqref{eq:56}
in the real part of~\eqref{1.10},
we obtain a polynomial of 
degree 5 in $r$; since $r$ is assumed to be different from $0$, we obtain through
straightforward algebra
\begin{equation*}
k^2 = \frac{3\alpha^2}{4m}, \quad 2\tilde \omega + ck = \frac{3\alpha}{
2},\quad \tilde \omega 
^2 -\alpha \tilde \omega +\alpha^2(m + 3\alpha^2/4)=0.
\end{equation*}
The second degree equation in $\tilde \omega$ has real roots if and
only if $4\alpha^2 - 
\nu<1/3$, which we assume from now on. These roots are given by
\begin{equation}
\bar \omega = \frac{\alpha}{ 4}\left(2 \pm\sqrt{1 + 3(\nu - 4
\alpha^2)}\right).\label{eq:41}
\end{equation}
We see also that 
\begin{equation}
\bar r=\sqrt{\tilde \omega/\alpha}.\label{eq:42}
\end{equation}
If we linearize the equation for $v=u(x-ct)
\exp\bigl(-ik(x-ct)-i\omega t\bigr)$ around $\bar r$, the condition
for stability of $0$-wave number modes is
\begin{equation*}
2\bar r^2 -4 \bar r^4 <0.
\end{equation*}
Comparing this relation with~\eqref{eq:42}, we can see that we have to
choose the $+$ sign in~\eqref{eq:41}.

Hence, the velocity $c$ is given by
\begin{equation*}
c = \frac{\sqrt 3 (4\alpha^2 -  \nu)}{ 1 + \sqrt{1 + 3(\nu - 4 \alpha^2)}}.
\label{1.12}
\end{equation*}
where we have used the sign condition~\eqref{eq:43}.

If $4\alpha^2 > \nu$, the velocity $c$ is positive, and the zero state
gains over the non zero state; we will 
say that this is the inflow situation; on the contrary, if
$4\alpha^2<\nu$, the non zero 
state gains over the zero state; we will say that we have an outflow situation.

Assume from now on that $\alpha$ and $\nu$ are small. If there existed a
pulse-like solution, it could probably be approximated by a
combination of a kink $K(x+L 
-ct) e^{i\omega t}$ centered at $-L$ and an antikink $K(L-x-ct)$
centered at $L$, provided 
that we know how to glue their phases together. In the outflow case,
two competing effects 
take place: on one hand, the choice of parameters tends to produce an
expanding pulse; 
on the other hand, an attraction effect between the walls is expected
as in \cite{CaPe}, \cite{FuHa} 
and \cite{Fu}. It is natural to expect that this attraction effect
should be exponentially small 
with $L$. Thus, it is reasonable to conjecture that the evolution of
the pulse will be given by
\begin{equation*}
\dot L = - c -C_1e^{-C_2 L}.
\end{equation*}
The pulse will be in an equilibrium if the two competing effects
balance, i.e.
\begin{equation}
L \sim \frac{1}{C_2} \ln\frac{1}{\nu - 4\alpha^2}.\label{1.13}
\end{equation}
This analysis assumes an almost scalar pulse, so that it is easy to
glue together the 
phase of the kink and of the antikink. Of course, the pulse obtained
by this argument is 
not stable: any amount of swelling or shrinking of $L$ destabilizes it.

Let us consider the inflow case: $4\alpha^2 >\nu$. Now, nothing scalar can stop the
pulse from collapsing. However, Malomed and Nepomnyashchy argue in \cite{MaNe} that
the pulse does not 
collapse because of phase incompatibility: a minimum distance is
necessary to match the 
phases of the kink and of the antikink. By formal asymptotics, they
obtain a half length 
of the pulse given by
\begin{equation}
L \sim \frac{C_3}{4\alpha^2 - \nu},\label{1.14}
\end{equation}
in a very small range $\alpha^4 \ll 4\alpha^2 - \nu\ll \alpha^2$.

However, the second author of the present article tends to believe
that~\eqref{1.14} holds for 
a much larger range of $\alpha$. The article \cite{MaNe} contains no analysis of
the stability of the pulse 
obtained; however the claim of stability (in a mathematical sense) seems quite
reasonable.

\begin{figure}
\begin{center}
\psfrag{dMS}{de Mottoni - Schatzman}
\psfrag{MN}{Malomed - Nepomnyaschy}
\psfrag{longueur}{$L$}
\psfrag{alpha}{$\alpha$}
\includegraphics[width=0.90\textwidth]{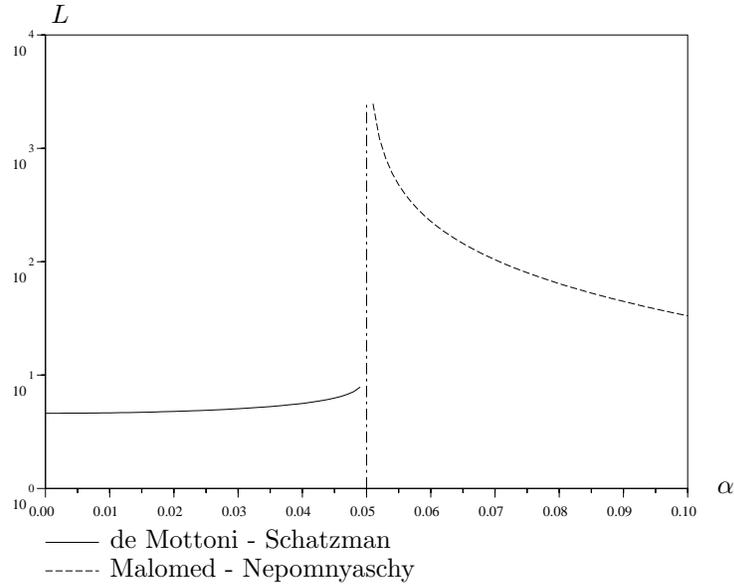}
\end{center}
\caption{The length $L$ of the pulse as a function of $\alpha$..}
\label{fig:length}
\end{figure}

With our results and the results of Malomed and Nepomnyashchy, we can
plot a graph of 
$L$ as a function of $\alpha$, and we obtain Fig.~\ref{fig:length}.

The region around $4\alpha^2=\nu$ seems difficult and interesting;
jumping somewhat too 
fast to conclusions, it would be nice to believe that these two branches
join to form a single branch; proving this might mean considerable effort. 

Let us now give an idea of the mechanism of skew stabilization. We
recall that we assumed that the shelf in the solution is very large,
i.e. $\nu$ is very small. The definition of $L$ is such that for
$\alpha=0$, the width of the shelf of the pulse is approximately $2L$.
 We define a function 
\begin{equation*}
F(u,\omega,\alpha,\nu) = i\omega u -m u'' +m u -(1+i\alpha)|u|^2 u +
|u|^4u.\label{1.15}
\end{equation*}
A function $u e^{i\omega t}$ is a solution of~\eqref{1.5} iff $u$ and
$\omega$ are such that $F(u,\omega,\alpha,\nu)$ vanishes. Moreover, it
is linearly stable if $D_u F(u,\omega, \alpha,\nu)$ has its spectrum
included in the right-hand side complex half-plane; moreover, we
require for stability that $0$ should be a semi-simple eigenvalue of
finite multiplicity while the remainder of the spectrum is bounded
away from the imaginary axis.

The symmetries of the problem imply that there are two eigenfunctions
of $D_u F$, relative to the eigenvalue $0$: $t=iu$ and $u'$; there is
also the eigenfunction $w$ relative to the shrinking-swelling mode,
and the aim of the game is to show that the corresponding eigenvalue $\lambda$
crosses $0$ for some appropriate value of $\alpha$. Let us define
\begin{equation*}
\tilde F(u,\tau,L,\beta)=F\bigl(u,e^{-2L}\beta \tau, e^{-2L}\beta, e^{-4L}\bigr).
\end{equation*}
At this point, we will use a measure of cheating in order to explain
what is going on: if $u$ and $\tau$ can be seen as smooth functions of $L$, we
differentiate with respect to $L$, and we find
\begin{equation}
\left\{\begin{split}
D_u \tilde F(u(\cdot, L, \beta),& \tau(L, \beta),L, \beta) = e^{-2L}
\beta \frac{\partial \tau}{\partial L} u + 2 ie^{-2L}\beta\bigl(\tau u
-\abs{u}^2u\bigr) + \\
&\qquad\frac{4\nu}{1-\nu}\bigl[i e^{-2L}\beta \tau u
-\abs{u}^2u +\abs{u}^4 u -ie^{-2L} \beta \abs{u}^2u\bigr].
\end{split}\right.\label{eq:55}
\end{equation}
Since we expect the interaction between the rotation mode and the
shrinking-swelling mode to be the reason for the skew stabilization,
we let $\bigl\{\hat w, \hat t\bigr\}$ be a dual basis to the basis
$\bigl\{w,t\bigr\}$, i.e. 
\begin{equation*}
\begin{split}
&\int \hat t^{\,T} t\, dx =\int \hat w^{\,T} w\, dx =1, \quad \int \hat
t^{\,T} w\,dx=\int \hat w^{\,T} t\, dx=0,\\
&\text{$y$ and $t$ are even, } \hat t^{\,T} D_u \hat
F(u(L,\beta),\tau(L,\beta),L,\beta)=0,\\& \hat 
w^{\,T} D_u F(u(L,\beta),\tau(L,\beta),L,\beta)=\hat w^{\,T} \lambda(L,\beta).
\end{split}
\end{equation*}
We multiply~\eqref{eq:55} on the left by $\hat w^{\,T}$, we integrate,
and we find that
\begin{equation*}
\begin{split}
\lambda(L,\beta) \int \hat w^{\,T} \frac{\partial u}{\partial L}\, dx
&= -2e^{-2L} \beta \int \hat w^{\,T} i\abs{u}^2 u\, dx \\ &+
\frac{4\nu}{1-\nu} \int\hat w^{\,T} \bigl[-\abs{u}^2 u +\abs{u}^4 u
-ie^{-2L}\beta \abs{u}^2u\bigr]\, dx.
\end{split}
\end{equation*}
Since $\alpha$ is small, the problem is very close to being
self-adjoint, it is not unreasonable to take $w$ as an approximation
of $\hat w$; moreover, in the approximation of the large shelf, $u$ 
is very close to $\tilde r(\vabs{x}-L)$; therefore, $w$ and $\partial 
u/\partial L$ can be taken very
close to $-\tilde r'(\vabs{x}-L)\sgn x$. These considerations imply that
\begin{equation*}
\int \hat w^{\,T} \frac{\partial u}{\partial L}\, dx \sim 2 \int_{-L}^\infty \vabs{\tilde
r'}^2\, dx = \frac 38.
\end{equation*}
We see also that
\begin{equation*}
\int \hat w^{\,T} \bigl[-\abs{u}^2 u +\abs{u}^4 u\bigr]\, dx \sim 2
\int_{-L}^\infty \tilde r' \tilde r^3  -\tilde r' \tilde r^5\, dx
\sim -\frac{9}{64}.
\end{equation*}
Thus, we obtain the following ``equivalent'' for $\lambda(L, \beta)$:
\begin{equation*}
\lambda(L,\beta)\sim -\frac{3e^{-4L}}{2}  -2e^{-2L} \beta \int \hat
s^{\,T} i\abs{u}^2 u\, dx. 
\end{equation*}
If we are able to calculate with enough precision the integral in the
right hand side of the above equation, and if it turns out to be
negative, we will hope for a stabilization effect.

However, we have cheated too much for this argument to go through a
rigorous analysis; the main source of inexactitude comes from the
assumption that for $\alpha$ small, we can find an even solution $u$
of $\tilde F(u, \tau, L, \beta)$ which is close to the real even
solution of $\tilde F(u,0,L,0)=0$; indeed, $u$ can be found close to
the real even solution of $\tilde F(u,0,L+y,0)=0$, where $y$ is
related to $\alpha$ by the relation
\begin{equation*}
y=-\frac14 \ln\bigl(4\beta -2\bigr) =-\frac14 \ln\bigl(4 e^{2L}\alpha -2\bigr).
\end{equation*}
In the course of the proof of skew stabilization, we will discover 
the following spatial scales in the problem: we have already $1$ and
$L$; the third scale is $\ln L$, which is the order of magnitude of
the critical $y$ for which the skew stabilization occurs.

We can now explain the organization of the paper; since we are in an
almost scalar situation, we let $\alpha=\sqrt\ve$ and we define 
\begin{equation*}
G(\xi, \eta,\tau,\varepsilon,\nu) = \Re F(\xi + i\sqrt\varepsilon \eta, \sqrt\varepsilon \tau,
\sqrt\varepsilon, \nu) +\frac{i}{ \sqrt\varepsilon}\Im  F(\xi + i\sqrt\varepsilon
\eta, \sqrt\varepsilon 
\tau, \sqrt\varepsilon, \nu).\label{1.16}
\end{equation*}
In our continuation analysis, we need a starting point $(u,\omega)$ such that
\begin{equation}
F(u,\omega,0,\nu)=0.\label{1.17}
\end{equation}
In section~\ref{sec:scal-equat}, we prove that a solution $u$ of~\eqref{1.17} exists and
decays exponentially at 
infinity only if $m$ belongs to the interval $(0, 3/16)$. This solution is, up to phase
and space translations, the unique positive even solution of
\begin{equation*}
-m r'' +m  r - r^3 + r^5=0,
\end{equation*}
as is proved at Lemma~\ref{Lemma2.1}.

Then, we proceed to study precisely this $r=r(x,\nu)$, which has an
explicit expression; 
however, we use mostly the asymptotic for $r(x,\nu)$, when $\nu$ is
small (Lemma~\ref{Lemma2.2}). Thus, we see that the shelf of $r$ has
indeed a half-length of
$L$ defined
by~\eqref{1.7}. We study the linearized 
operator $A$ around $r$; it is an unbounded self-adjoint operator in $L^2(\Er)$ given by
\begin{equation*}
D(A) = H^2(\Er), \quad Au = -m u'' +(m - 3 r^2 + 5 r^4)u.
\end{equation*}
We prove that the spectrum of $A$ contains a group of eigenvalues $\{\lambda, 0\}$, and
that
\begin{equation*}
\lambda \sim -\frac{3\nu}{ 2}
\end{equation*}
is the lowest eigenvalue of $A$; moreover, this group of eigenvalues
is bounded away from 
the remainder of the spectrum, uniformly in $\nu$. The eigenvalue $0$
corresponds to the 
translation mode, with eigenfunction $r'$, and the eigenvalue
$\lambda$ corresponds to 
the shrinking-swelling mode, with eigenfunction $s$. This analysis
is made possible by 
the following fact: if
\begin{equation*}
\sigma(x,\nu) = - 4\nu \frac{\partial r(x, \nu)}{ \partial \nu} = \frac{\partial}{\partial
L} r(x,\nu),
\end{equation*}
then, with an appropriate normalization of $s$,
\begin{equation*}
|\sigma - s|_{H^2(\Er)} = O(\nu \sqrt L),
\end{equation*}
as is proved at Theorem~\ref{Theorem2.8}.

The existence will be proved using the scaled equation 
\begin{equation*}
G(\xi, \eta,\tau, \varepsilon, \nu) = 0.
\end{equation*}
The scaled equation is more interesting from the point of view of
continuation, because 
for $\alpha=0$, $\omega =0$, and it is easy to check that
\begin{equation*}
\Im F(r,0,0,\nu)=0,
\end{equation*}
which does not provide any information.
However,
\begin{equation*}
\Re G(\xi, \eta,\tau, 0,\nu)=\Re F(\xi, 0,0,\nu)
\end{equation*}
so that $r$ satisfies
\begin{equation*}
\Re G(r,\eta,\tau, 0, \nu)=0,
\end{equation*}
and $\eta$ and $\tau$ are yet undetermined. The second equation is
\begin{equation}
\Im G((r,\eta,\tau, 0, \nu) = \tau r - r^ 3 - m \eta'' +m \eta -r^2 \eta + r^4
\eta=0,\label{1.18}
\end{equation}
and it is studied in details in Section~\ref{sec:second-equation}.
Let $\theta$ and $q$ denote the values of $\eta$ and $\tau$ which
satisfy~\eqref{1.18}. It is natural to study the operator $B$ in
$L^2(\Er)$ defined by
\begin{equation*}
D(B) = H^2(\Er), \quad B u = -m u'' +(m - r^2 + r^4) u.
\end{equation*}
We can see immediately that $0$ is the lowest eigenvalue of $B$ and the corresponding
eigenvector is $r$; this is not surprising, since it is the analytical
translation of the 
$S^1$ equivariance of~\eqref{1.5}. In other words, $ir$ is the phase rotation mode. 
Therefore,~\eqref{1.18} will have a solution if and only if $r$ is
orthogonal to $r^3 - \theta 
r$; this determines $\theta$; if we impose that $q$ be orthogonal to
$r$, it is uniquely 
determined.

For later  purposes, we need an asymptotic on the second eigenvalue $\mu_2$ of
$B$; it is proved at Theorem~\ref{Lemma3.4} that $\mu_2\sim C/L^2$, where $C$ is a
positive constant.

Section~\ref{sec:Existence-pulses} is devoted to the existence proof. Preliminary
computations showed
that continuation is not good enough to obtain a satisfactory range of
existence; what is 
needed is an ansatz for the pulse; it is obtained by taking $\nf = \nu
e^{-4y}$, where 
$y$ is some positive number bounded by $L^p$; the corresponding $r$,
$\theta$ and $q$ are 
denoted by $\rf$, $\thf$ and $\qf$. Now, $\ef$ has to be determined: this
is a version of the Lyapunov-Schmidt method
of bifurcation theory. Our choice is to require that $\Re
G(\rf,\qf,\thf ,\ef,\nu)$ is 
orthogonal to $s\f$, the shrinking-swelling mode corresponding to
$\nf$. An asymptotic for 
$\ef$ is given by
\begin{equation*}
\ef\sim \frac{\kappa}{ 4},
\end{equation*}
where 
\begin{equation*}
\kappa =(\nu -\nf)/(1 - \nf).
\end{equation*}
Let $\Uf$ be the vector  of components $(\rf,\qf,\thf,\ef)$. The idea
is to observe that 
this $\Uf$ is almost a solution of $G(\Uf, \nu)=0$. Existence is
proved using a version  
of the
implicit function theorem with estimates, proved in
Section~\ref{sec:Append-An-impl} . The 
pulse obtained this way 
is denoted by $u$.

In other words, we approximate the pulse at $\alpha=\sqrt{\ef}$ and
$\nu$ by the pulse at $\alpha=0$ and $\nuf$.

Let us denote by $\co$ any quantity bounded by a finite power of $L$.

The main result of this article, i.e. the proof of stabilization
(section~\ref{sec:Stability-pulse}) uses the details of the proof of existence. If 
$(1-\kappa)\cald $ is the differential of $F$ with respect to $u$ at
$(u,\omega, \alpha,\nu)$, then,
\begin{equation*}
\cald =\cala^\flat +\sqrt\kappa\,\calb+\kappa \calc.
\end{equation*}
Here,
\begin{equation*}
\cala^\flat=  \begin{pmatrix}
\Af&0\\ 0&\Bf
\end{pmatrix}
,
\end{equation*}
where $\Af$ (resp. $\Bf$) is $A$ (resp. $B$) at $\nf$ instead of
$\nu$, and $\calb$, $\calc$ are 
$2\times 2$ matrix of multiplication operators $C_{ij}$, $1\le
i,j \le 2$ such that 
\begin{equation*}
\begin{split}
&\calb=\begin{pmatrix}
0&\calb_{12}\\\calb_{21}&0
\end{pmatrix},\quad \calc=\begin{pmatrix}
\calc_{11}&0\\0&\calc_{22}
\end{pmatrix},\\
&\vnorm[L^\infty]{\calb_{12}}+\vnorm[L^\infty]{\calb_{21}}+\vnorm[L^\infty]{\calc_{11}}+\vnorm[L^\infty]{\calc_{22}}
=\co. 
\end{split}
\end{equation*}

The idea 
is to consider the restriction of $\cald$ to the generalized
eigenspace corresponding to 
the eigenvalues of $\calb$ which are close to zero. This eigenspace is
of dimension $3$; a basis of it is $\{s, iu, u'\}$, where $iu$ spans the phase
rotation mode, $u'$ spans 
the space translation mode, and $s$ corresponds to the shrinking-swelling mode. In this
basis, the matrix of $\calb$ is given by
\begin{equation*}
\begin{pmatrix}
M_{11} & 0 &0\\ M_{21}& 0 & 0\\ 0&0&0
\end{pmatrix}
.
\end{equation*}
The sign of $M_{11}$ determines the stability of the pulse: if $M_{11}>0$, the
pulse is stable (up 
to space and phase translation); if $M_{11}\le 0$, the pulse is unstable. Thus
we have to give 
an asymptotic for $M_{11}$. For this purpose, we embed the operator
$\calc$ into an holomorphic family $\cald(c)$ of operators depending
on $c$, and we prove estimates using the strong properties of such
families. In particular, we give a precise description of the
expansion of a basis of eigenfunctions relative to the very small
eigenvalues of $\cald(c)$, and of the dual basis, and we validate
these expansions. With the residue theorem, we are able to describe
$M_{11}$ with sufficient precision, and the symmetries of the problem
lead us to an almost completely explicit value for it
(Lemma~\ref{Lemma5.4}, Theorem~\ref{Theorem5.5} and
Theorem~\ref{Theorem5.7}). 
We conclude this asymptotic analysis at Proposition~\ref{Proposition5.9}.

Thanks to a result of Henry \cite{He}, Chapter 5, Exercise 6, the
linearized  stability
implies the following non linear stability result: take an initial
condition for~\eqref{1.5} which is equal to a pulse plus a small perturbation
; then, if $|\alpha|\ge \alpha_c$, the asymptotic state of the solution
of~\eqref{1.5} is a pulse possibly translated in space and in phase.

In section~\ref{sec:general-situation} we give the analogous
asymptotic for the case when the $\mu_j$'s do not vanish; the appendix
(section~\ref{sec:Append-An-impl}) gives an implicit function theorem with estimates; this
theorem is the key to the existence result; in other words: our
existence result is based on an ansatz: if the ansatz is good enough, then it is
indeed a good approximation of the solution. When small parameters are
involved, a correct argument deserves a proof.

A rather curious fact is that the number $\pi^2/6$ appears in the
calculation of the expansion of $M_{11}$; in partial differential
equations, it is usually related to a trace, but we have been unable
to uncover such an origin; therefore, its presence may be a coincidence.

Some of these results were announced in \cite{MoSc2} which contains a
number of errors. A preprint~\cite{Mosc3} was circulated but never published as an
article; the present article contains for the first time the approximate explanation of
the skew-stabilization of the pulse and also the case of general
coefficients as in~\eqref{eq:46}.

There is considerable interest in the Ginzburg-Landau models; scanning the
literature, one can find for example~\cite{PuRo} which lists
fluctuations in lasers, order-disorder
transitions, population dynamics and ordering in uniaxial ferromagnetic films
as domains where Ginzburg-Landau of the third degree has been used as a model.
W. Eckhaus~\cite{Ec} states that Ginzburg-Landau of the third degree
is ``universal'' for modulation equations,
which is another way of saying that it behaves as a normal form.

The article \cite{SaHo} describes a large number of solutions
of~\eqref{1.1}, perturbation expansions 
for large values of $\Im m_j$ and gives conjectures on the behavior
of the solutions of~\eqref{1.1} in different regions of the parameter space.

Ginzburg-Landau of the fifth order is much less generic than the third
order Ginzburg-Landau. Its main
merit is that it allows 
for subcritical bifurcation of the constant amplitude solutions. 

Thual and Fauve explained in~\cite{ThuFau} the phenomenon they obtained in terms of the
general picture of 
subcritical bifurcation, and also as a perturbation with respect to a
nonlinear Schr\"odinger equation.
                                                   
Shortly later, Malomed and Nepomnyashchy \cite{MaNe} considered the same
equation~\eqref{1.1} and 
explained the existence of a pulse by formal asymptotics. A careful
examination of their 
results shows that they worked in a different range of parameters from ours. 

Hakim, Jakobsen and Pomeau \cite{HaJaPo} have given a general idea of the
bifurcation picture 
in a situation which is
close to the present one; however, it is difficult to compare the situations, since 
their
statements are not described with complete precision. One of their
statements is the 
subcritical character of the bifurcation.
While the bifurcation of space independent solutions is clearly subcritical,
the bifurcation of the pulse is subcritical only by the fact that
initially the solution is unstable, and it is stabilized along the solution curve;
however, the typical picture of subcriticality as in  Fig.1 of
\cite{ThuFau} has not been found in the present article.

Kapitula \cite{Kap} gave a general theory for the
existence of heteroclinic traveling 
wave solutions of the quintic Ginzburg-Landau equation with convective
terms. Kapitula 
basically studies the persistence under perturbations of heteroclinic
orbits close to  
orbits which can be obtained almost explicitly when all coefficients
are real. The proof 
relies on a very precise study of the perturbed invariant manifolds
for a flow associated to the system. Possibly, his methods could be
adapted to give existence of homoclinic solutions of~\eqref{1.3}.

According to a very striking phrase by Yves
Meyer~\cite{Me}, the present work is like the success of the lock breaker:
he/she has to use many pick locks, and the place looks messy;
however, it is expected that once the door is open, and the lock can be
dismantled and studied in details, someone will be able to devise a nice
key which will open it in a single move. 

Kapitula \cite{Kap} used an entirely different set of techniques to
devise a nice key for related existence questions,
but a nice key to stability does not yet exist.

This article makes available a number of pick locks, to be stored in the tool box
of the mathematical lock picker. It is more in the spirit of the SLEP
method of~\cite{MR91d:35024} than in the spirit of the large
literature on the analysis of the stability of traveling
waves: solutions of the quintic Ginzburg-Landau equation have been
analyzed in \cite{MR96a:35196}, \cite{MR99a:35238},
\cite{MR99h:35199}, \cite{MR99g:35118}, where the main difficulty is
the bifurcation from the essential spectrum; the foundational work on
the stability of traveling and standing wave solutions of
semilinear parabolic equations is related to the Evans function; see
in particular \cite{MR86b:35011}, \cite{MR92d:58028} and \cite{MR93g:35115}, which were
followed by a considerable literature, including in particular
\cite{MR93c:35009}, \cite{MR92d:58028}, \cite{MR97e:34105},
\cite{MR98b:35093}, \cite{MR1664760}, \cite{MR99g:35064}. The analysis of
perturbation of periodic states has been taken up by the Evans function
method in \cite{MR94j:35072} \cite{MR99a:35196} and by modulation equation
methods in \cite{MR99i:35156}. Solutions with several fronts or bumps
have been constructed in \cite{MR94k:35042}, \cite{MR94m:35152},
\cite{MR98d:35100}, \cite{MR98h:35227}, \cite{MR98h:35212} and
\cite{MR99b:35012}.

It is possible that the methods of this article are close to those of
\cite{MR99b:35094}, which however does not have the $S^1$
symmetry of the problem considered in the present article.

We would like to thank Stephan Fauve for introducing us to this
problem, which revealed itself as infinitely more complicated than
what we would have expected. The second author is glad to thank R.L.
Pego, B. Malomed and A. Nepomnyashchy for fruitful discussions and
exchanges of ideas.

St\'ephane Descombes read the article in detail in the course of his
first year of graduate studies, and spotted the defects, the typos and the errors.
His criticisms improved considerably the article and he deserves
praise and thanks for his patience.

\section{The scalar equation and the corresponding linearized
operator}\label{sec:scal-equat} 

We study in details the solutions of the equation
\begin{equation*}
u_t - u_{yy} +m u -|u|^2u(1+i \alpha) +|u|^4 u=0
\end{equation*}
under the assumption $\alpha=0$ and
\begin{equation}
m \in \Er\setminus{0}.\label{2.1}
\end{equation}

We look for solutions of the form                            
\begin{equation*}
u(y,t) = v(y) e^{i\omega t}.
\end{equation*}
If we substitute this expression into our equation, we can see that $v$ and
$\omega$ satisfy
\begin{equation}
i\omega v - v_{yy} +m v - |v|^2 v + |v|^4 v=0.\label{2.2}
\end{equation}

Our first and elementary results on this case are summarized in the
following lemma:

\begin{lemma}\label{Lemma2.1} Let 
$v\in L^5_{\text{loc}}(\Er;\Ce)$ solve~\eqref{2.2} in
the sense of distributions; if $v$ does not vanish identically,
then the following assertions hold:
\begin{enumerate}
\renewcommand{\labelenumi}{(\roman{enumi})}
\item $v$ is infinitely
differentiable. 

\item if $v(y)$ tends to
zero as $|y|$ tends to $+\infty$, then $\omega$ vanishes.

\item under the assumption of (ii), the argument of $v$ does not depend
on $y$.

\item under the assumption of (ii), the number $m$ is strictly positive.

\item under the assumption of (ii), let $x= y\sqrt{m}$ and
\begin{equation*}
|v(y)| = r(x);
\end{equation*}
then, up to translation in space,  $r$ is the unique even positive solution of the
ordinary differential equation 
\begin{equation}
-m r'' +m r - r^3 + r^5 =0\label{2.3}
\end{equation}
which vanishes at infinity. In particular, 
the number $m$ belongs to the interval $(0,3/16)$.
\end{enumerate}
\end{lemma}

\begin{proof}
(i) If $v$ belongs to $L^5_{\textrm{loc}}(\Er;\Ce)$,
then $v''$ belongs to 
$L^1_{\textrm{loc}}(\Er;\Ce)$, so that $v'$ is locally absolutely continuous,
and $v$ is a function of class $C^1$ on $\Er$, with values in $\Ce$. By an
obvious induction argument, $v$ is infinitely differentiable.

\medskip

(ii) Equation~\eqref{2.2} implies that $v''$ tends
to $0$ at infinity. 
Therefore, $v'$ tends to zero at $\infty$: by Taylor's formula,
\begin{equation*}
v'(y)=v(y +1)-v(y)-\int_0^1v''(y+s)(1-s)\, ds;
\end{equation*}
the right hand side of this relation tends to $0$ at infinity;
therefore, the left hand side must also tend to $0$.

Assume $\omega \neq 0$. System~\eqref{2.2} is equivalent to a system
of four ordinary 
differential equations of the first order in $\Er^4$,
and $0\in \Er^4$ is a critical point of this system. We apply the
theory of invariant 
manifolds in a neighborhood of $0$; the matrix of the linearized
system at that point has 
two double eigenvalues $\pm \zeta$, where $\zeta^2= m +i\omega$. We
make the convention 
that $\Re\zeta>0$. Therefore, $|v|$, $|v'|$ and $|v''|$ decrease
exponentially  fast at 
infinity. Multiply~\eqref{2.2} by $\bar v$, and integrate over $\Er$;
we obtain after an 
integration by parts
\begin{equation*}
i\omega \int |v|^2\,dy + \int |v'|^2 \,dx +\int\bigl( m |v|^2 - |v|^4 +
|v|^6\bigr)\,dy=0.
\end{equation*}
This shows that $\omega$ vanishes.  

\medskip

(iii) Define 
\begin{equation*}
W(y) = m -|v(y)|^2 +|v(y)|^4.\label{2.5}
\end{equation*}
We let $v = v_1+ iv_2$; then $v_1$ and $v_2$ solve the linear differential equation
\begin{equation}
w'' = W w,\label{2.6}
\end{equation}
and they cannot both vanish identically. 
If $w$ is a solution of~\eqref{2.6} which does not vanish identically,
then $w$ and $w'$ cannot vanish simultaneously; if $w_1$ and $w_2$ are
two solutions of~\eqref{2.6} which vanish at infinity, their Wronskian
$w'_1 w_2 -w'_2 w_1$ is constant and vanishes at infinity; therefore,
any two solutions of~\eqref{2.6} are proportional. Thus $v_1$ and
$v_2$ are linearly dependent, which proves the desired assertion.

\medskip

(iv) 
Thanks to (iii), if $v$ is a solution of~\eqref{2.2}, with
$\omega=0$, we may assume that it is real without loss of
generality. Thus, it solves~\eqref{2.6}; moreover,
$v$ and $v'$ cannot vanish simultaneously by uniqueness of solution of the Cauchy
problem for~\eqref{2.2}. Therefore, we can find a continuous
determination of the angle $\theta$ such that
\begin{equation}
w=r\cos\theta, w'=r\sin\theta.\label{eq:1}
\end{equation}
It is immediate that this determination is also of class $C^\infty$.
Let us multiply~\eqref{2.2} by $2 v'$, and integrate,
remembering that $\omega =0$ and $v$ is real; we obtain
\begin{equation}
-|v'|^2 +m |v|^2 - |v|^4/2 +|v|^6/3=\textrm{constant}.\label{2.4}
\end{equation}
By taking the limit of the left hand side of~\eqref{2.4}, we see that
the constant on the right hand side of~\eqref{2.4} vanishes.
We use~\eqref{eq:1}, and we can see after substituting
into~\eqref{2.4} and dividing by $r^2$ that
\begin{equation*}
-\sin^2\theta +m \cos^2\theta =\frac{r^2}{2} -\frac{r^4}{3}.
\end{equation*}
Assume that $m$ is strictly negative.
As $x$ tends to infinity, the upper limit of the left hand side of the
above equation is at most 
equal to $-\min\bigl(1,\abs{m}\bigr)$, while the right hand side of this
equation tends to $0$; thus, we have a contradiction.

\medskip

(v) The proof of (iv) shows that $v$ can be taken real; define a
function $r$ by
\begin{equation*}
v(y) = z r(y\sqrt{m});
\end{equation*}
then $r$ is a real solution of 
\begin{equation*}
 -m r'' + m r - r^3 + r^5=0
 \end{equation*}
which decays at infinity. Moreover, if we define
\begin{equation*}
\Phi(r) = m r^2 - \frac{r^4}{2} + \frac{r^6}{3},
\end{equation*}
\eqref{2.4} implies
\begin{equation*}
m \abs{r'}^2 = \Phi(r).
\end{equation*}
\begin{figure}
\begin{center}
\includegraphics[width=10cm]{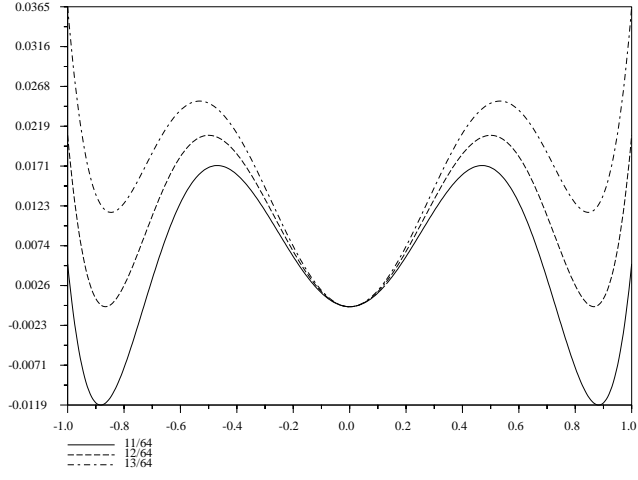}
\end{center}
\caption{The graph of $\Phi$ for $m=11/64, m=12/64, m=13/64$.}
\label{fig:plotphi}
\end{figure}

In order to have a non trivial solution of~\eqref{2.3} which tends to $0$ at
$\pm \infty$, we  
must choose $m$ so that $\Phi$ can vanish for non zero values of $r$; letting 
\begin{equation*}
\Psi(X)= m - \frac{X}{2} + \frac{X^2}{3}=X^{-1} \Phi(\sqrt X),
\end{equation*}
it is immediate that $\psi$ can vanish for positive values of $X$ if and only if 
\begin{equation*}
m <  \frac{3}{16}.
\end{equation*}
Henceforth, we will write
\begin{equation*}
m = \frac{3}{16}(1-\nu) \Longleftrightarrow \nu =  1 -\frac{16 m}
{3}.
\end{equation*}
The level curves of $(r,s) \mapsto s^2 - \Phi(r)$ are represented at
Fig.~\ref{fig:contour}; this figure 
shows that there
exist exactly two homoclinic orbits of~\eqref{2.3} through 0 for $m<3/16$; on one
of them, $r$ takes only 
positive values, and on the other one $r$ takes only negative
values. Therefore, we choose 
the homoclinic orbit which is situated in the half plane $r\ge 0$. 
The solution $r$ is still defined only up to translation. When $r$
reaches a maximum at 
some point $x_0$, $r'(x_0)=0$, and $\psi(r^2(x_0))$ vanishes.Therefore,
\begin{equation}
r^2(x_0) = \frac{3(1-\sqrt\nu)}{ 4}.\label{eq:2}
\end{equation}
Since~\eqref{2.3} is invariant by the reflexion $x\mapsto 2x_0 - x$,
we can see that $r(x) = 
r(2x_0 - x)$; if $u$ has more than one maximum, it is periodic, which
is possible only if 
$u$ vanishes identically. Therefore, $r$ has exactly one maximum. 

\begin{figure}
\begin{center}
\hbox{\includegraphics[width=6.8cm]{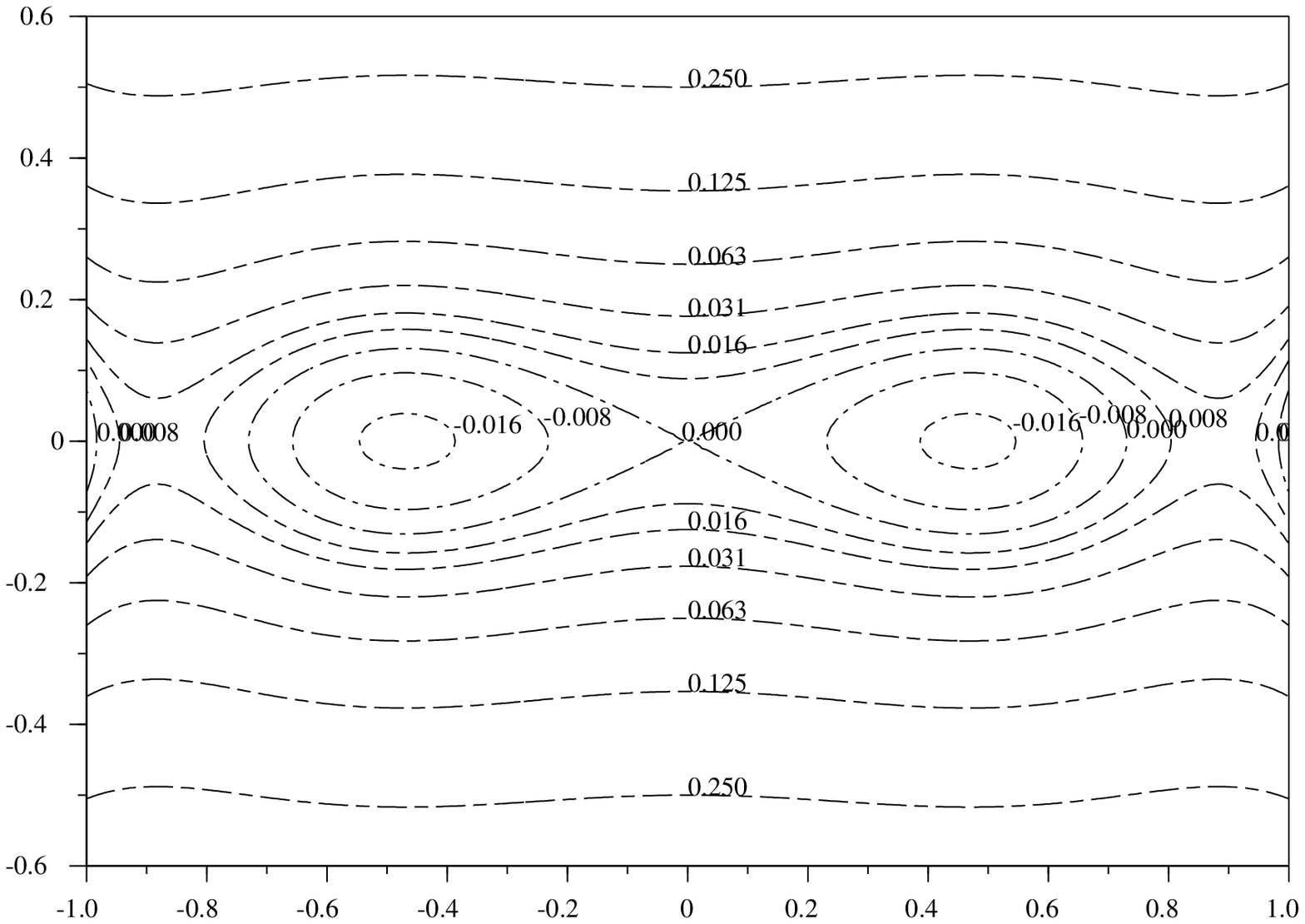}
\hfill\includegraphics[width=6.8cm]{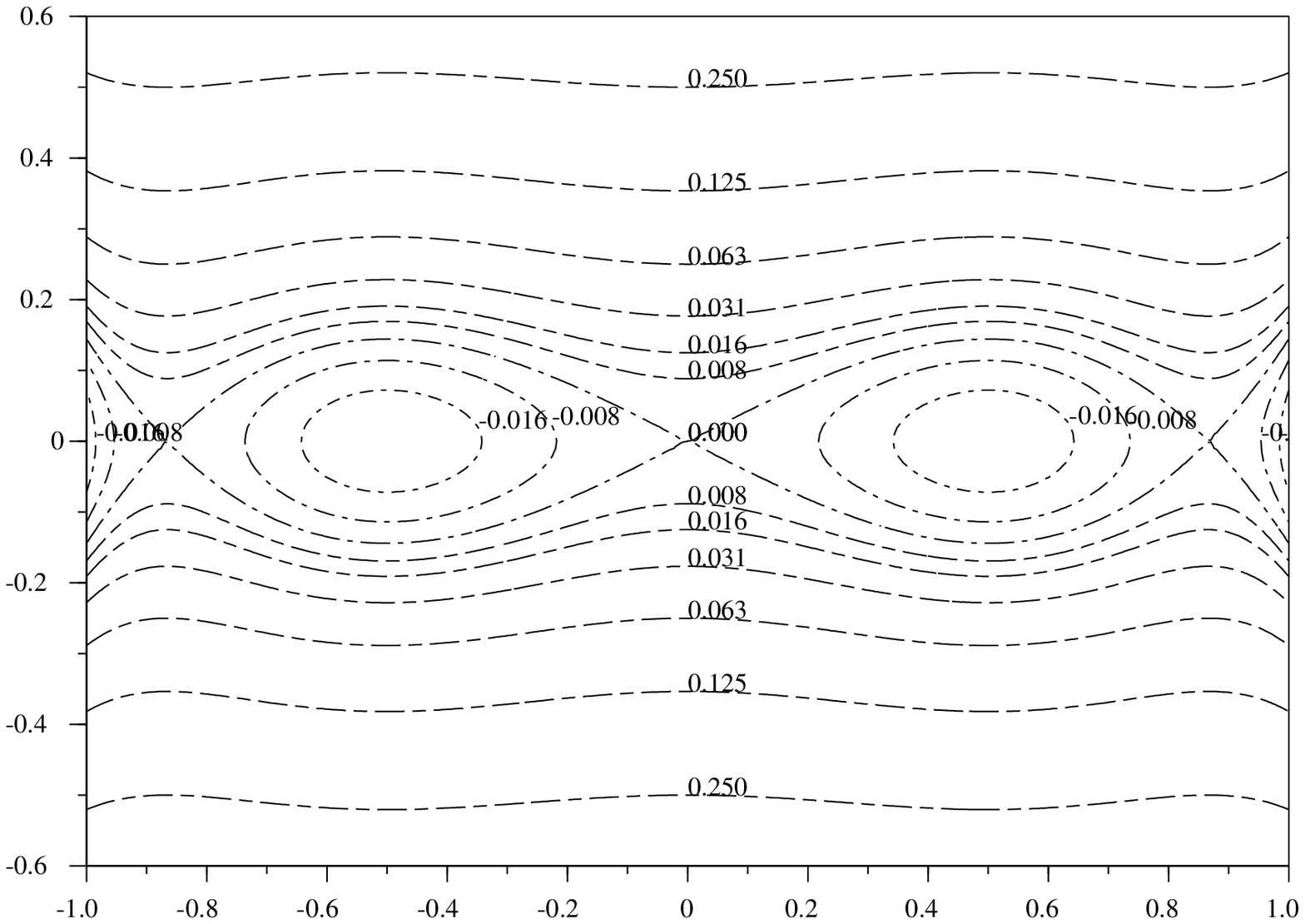}}
\end{center}
\begin{center}
\includegraphics[width=6.8cm]{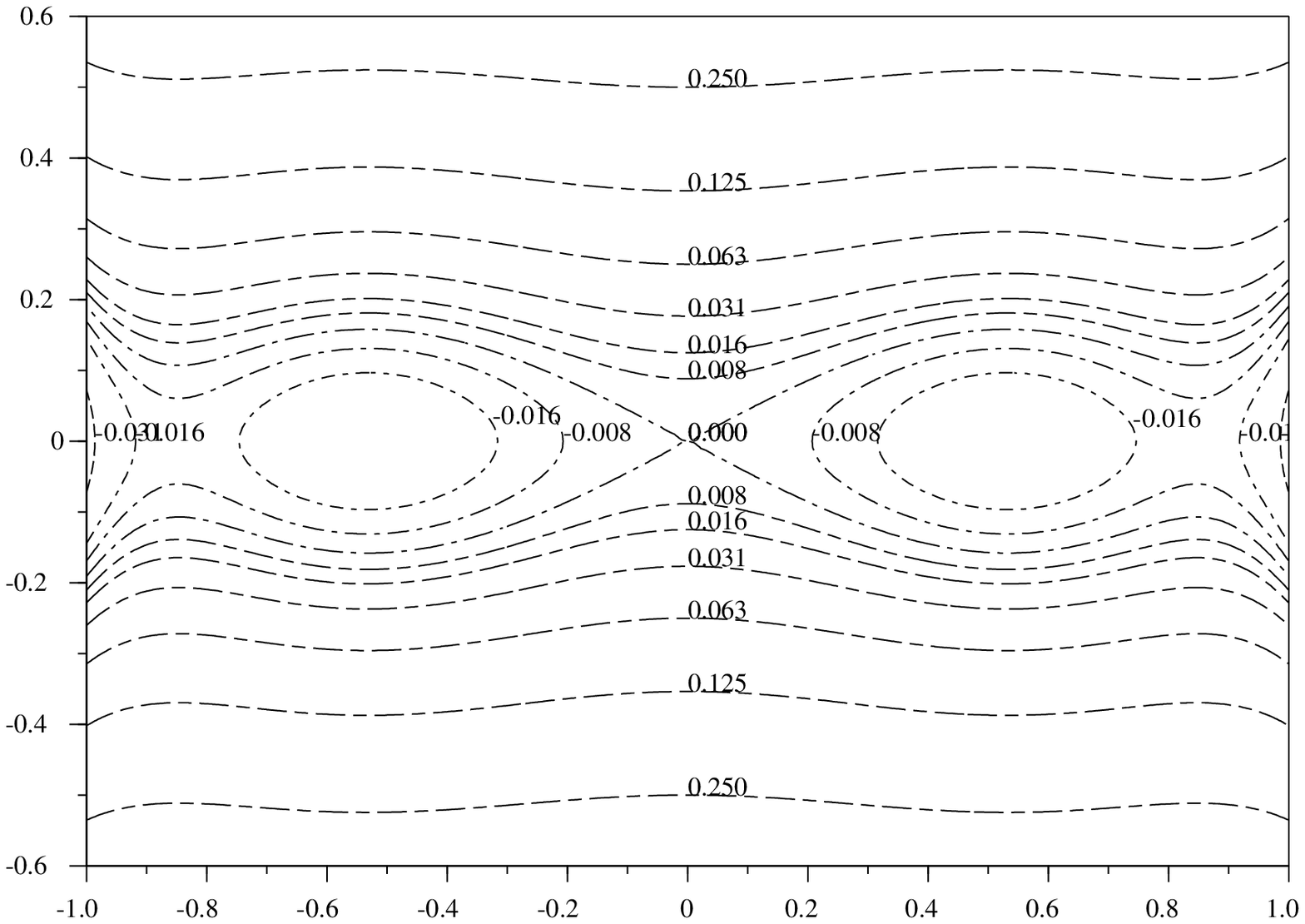}
\end{center}
\caption{The level lines of $(r,s)\mapsto s^2-\Phi(r)$ for $m=11/64$,
$m=12/64$, $m=13/64$.}
\label{fig:contour}
\end{figure}

If we choose the translation parameter so that $r$ attains its maximum
at $0$,  then $r$ 
is even and positive over $\Er$, which concludes the proof of the lemma.
\end{proof}

In what follows, the dependence of $r$ or $R=r^2$ over $\nu$ will be
emphasized from time to 
time, in which case, we will write $r(x, \nu)$ or $R(x,\nu)$ instead
of $r(x)$ or $R(x)$. 

The value of $R$ can be found explicitly; it is equal to 
\begin{equation*}
R(x) =\frac{3}{4}\frac{1-\nu}{1 + \sqrt\nu \cosh 2x} .
\end{equation*}

\begin{figure}
\begin{center}
\includegraphics[width=10cm,height=5cm]{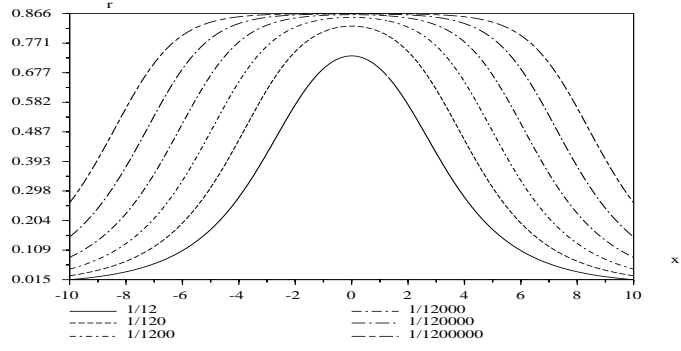}
\end{center}
\caption{ The graph of $R$ over the interval $[-10,10]$, for $\nu = 10^{-p}/12,
0 \le p\le 5$.}
\label{fig:plotr}
\end{figure}

If we draw a graph of the values of $r$  or $R$ for different values of $\nu$,
we can observe  
that the `width' of $r$ or $R$ increases logarithmically with $1/\nu$;
see~\ref{fig:plotr}. 

In fact,  as $\nu$ decreases to $0$ $R$ is
very close to a two front solution, and the distance between these two
fronts can be estimated very precisely. This important observation in 
stated and proved in next Lemma:

\begin{lemma}\label{Lemma2.2}
Let 
\begin{equation}
L=\frac{1}{4} \ln\frac{4}{\nu} \Longleftrightarrow
\nu=4e^{-4L},\label{2.9}
\end{equation}
and
\begin{equation}
\tilde R(x) =\frac{3}{4} \frac{1}{ 1 +  e^{2x}}.\label{2.10}
\end{equation}
Then, we have the estimates for all $x\in \Er$
\begin{gather}
\vabs{R(x+L) - \tilde R(x)} \le \nu e^{-2x}.\label{2.11}\\
\vabs{R'(x+L) - \tilde R'(x)} \le O(1) \nu e^{-2x}.\label{2.12}
\end{gather}
Moreover, $\tilde r=\sqrt {\tilde R}$ is a solution of the differential equations
\begin{equation}
-\frac{3}{16}\tilde r'' +\tilde W \tilde r =0\label{2.13}
\end{equation}
where
\begin{equation}
 \tilde W = \frac{3}{16} -
 \tilde R + \tilde R^2,\label{2.14}
 \end{equation}
  and
\begin{equation}
\tilde r' = - \tilde r \left(1 - \frac{4\tilde r^2}{
3}\right).\label{2.15}
\end{equation}
\end{lemma}

\begin{proof} We infer from~\eqref{2.10} the relation
\begin{equation}
e^{-2x} = \frac{\tilde R(x)}{3/4 -  \tilde R(x)};\label{2.17}
\end{equation}
it is also immediate that $\tilde R$ satisfies the ordinary
differential equation
\begin{equation}
\tilde R' = -2\tilde R(1 - 4\tilde R/3),\label{2.21}
\end{equation}
which can also be found by passing to the limit as $\nu$ tends to $0$
in the differential equation satisfied by $R$:
\begin{equation}
R' = {\pm 2} R\sqrt{\psi(R)/m}.\label{eq:6}
\end{equation}
The definition~\eqref{2.9} of $L$ and a direct computation give
\begin{equation}
R(\cdot+L) = \frac{(1-\nu)\tilde R(1-4\tilde R/3)}{1-4\tilde R/3 + 4\nu
\tilde R^2 /9},\label{2.18}
\end{equation}
from which we infer the identity
\begin{equation}
R(\cdot+L) - \tilde R =  -\frac{\nu\tilde R(1 - 2\tilde R/3)^2}{ 1 -
4\tilde R/3 + 4\nu \tilde 
R^2/9} .\label{2.19}
\end{equation}
Thanks to~\eqref{2.17}, we deduce immediately~\eqref{2.11}
from~\eqref{2.19}. If we differentiate~\eqref{2.19} with respect to
$x$, we find
\begin{equation*}
\begin{split}
&R'(\cdot +L)-\tilde R'\\&=-\nu\tilde R'\left[(-4/3 +8\nu\tilde
R/9)\frac{\tilde R(1-2\tilde R/3)^2}{(1 - 4\tilde R/3 + 4\nu\tilde
R^2/9)^2} -\frac{(1-8\tilde R/3 +4\tilde R^2/3}{1 - 4\tilde R/3 +
4\nu\tilde R^2/9}\right]. 
\end{split}
\end{equation*}
Then, with the help of~\eqref{2.21} and~\eqref{2.17},
relation~\eqref{2.12} is clear. 
By a passage to the limit as $\nu$ tends to $0$, or by a direct
computation, we see
that~\eqref{2.13},~\eqref{2.14} and~\eqref{2.15} hold.
As $x$ tends to $+\infty$, $\tilde R$ tends to 0 and $ V$ tends to $3/16$. 
\end{proof}

The linearized operator around $r$ is an unbounded operator $A$ in
$L^2(\Er)$ defined by 
\begin{equation*}
D(A) = H^2(\Er), \quad A  u = -m u'' + Vu,\quad V = m - 3r^2 +
5r^4.
\end{equation*}

The spectral properties of $A$ are important for what follows:

\begin{lemma}\label{Lemma2.3}
The operator $A$ is self-adjoint; its continuous spectrum is
included in the interval $[m, +\infty)$; the infimum of the spectrum of $A$ is a
strictly negative
number $\lambda$, to which corresponds an even strictly positive
eigenfunction $s$. 
The function $r'$ is an eigenfunction of $A$ corresponding to the
eigenvalue $0$. Both 
$\lambda$ and $0$ are simple eigenvalues.
\end{lemma}

\begin{proof}
Clearly, $A$ is self-adjoint.
We have 
\begin{equation*}
V(\pm \infty)= m>0;
\end{equation*}
The essential spectrum of $A$ is contained in $[m, +\infty)$ (see
Kato, \cite{Kat},Theorem 
5.26),
and  the spectrum of $A$ is bounded from below by $\min_x V(x)$; for
all $\beta>0$, the 
intersection of the spectrum of $A$ with $(-\infty, m-\beta]$ contains only
eigenvalues of finite multiplicity; since $A$ is a Schr\"odinger operator in
one-dimensional space, these eigenvalues are all simple; moreover,
since $r$ is even, $V$ 
is even, and the eigenfunctions of $A$ are even or odd; in particular, an
eigenfunction corresponding to the lowest eigenvalue $\lambda$ of
$A$ does not change 
sign, thanks to the maximum principle; therefore, it is even, since such an
eigenfunction cannot vanish together with its derivative.

Denote by $s$ an eigenfunction of $A$ corresponding to $\lambda$;
this $s$ can be 
chosen so that
\begin{equation*}
s(x) >0, \forall x \in \Er.
\end{equation*}

Differentiating~\eqref{2.3} with respect to $x$, we obtain
\begin{equation*}
A r'=0.
\end{equation*}
It is clear that $r'$ belongs to $L^2(\Er)$; therefore, it is an
eigenvector of $A$, 
corresponding to the eigenvalue $0$. On the other hand, $r'$ is odd;
therefore, we have 
proved that the infimum of the spectrum of $A$ is a strictly negative
number, $\lambda$. 
\end{proof}

We will obtain later much more precise information on the spectrum of
$A$; but this 
information is easier to obtain if we take into account the
nonlinearities of our problem 
; in particular, we will asymptotically describe the eigenfunction
$s$ in terms of the 
derivative of $r$
with respect to $\nu$, a feature which is clearly specific to our
nonlinear situation. 
Thus, we have to study the dependence with respect to $\nu$ of a
number of objects which 
appear in this section.

\begin{lemma}\label{Lemma2.4}
\begin{enumerate}
\renewcommand{\labelenumi}{(\roman{enumi})}
\item The function 
\begin{equation*}
\begin{matrix}
]0,1[&\to & H^2(\Er)\\
\nu& \mapsto & r(x,\nu)
\end{matrix}
\end{equation*}
is infinitely differentiable.

\item The eigenvalue $\lambda(\nu)$ and the eigenprojection  $P(\nu)$ on $\Er s(\nu)$
are infinitely differentiable, respectively with values in $\Er$ and in the space of
continuous operators in $L^2(\Er)$.
\end{enumerate}
\end{lemma}

\begin{proof}
(i) Denote by $H^2_{\textrm{even}}(\Er)$ (resp. $L^2_{\text{even}}(\Er)$) the subspace of even 
functions belonging to $H^2(\Er)$ (resp. $L^2(\Er)$). Define a mapping $f$ from
$H^2_{\textrm{even}}(\Er) \times (0,1)$ to $L^2_{\textrm{even}}(\Er)$ by
\begin{equation*}
f(v,\nu) = -m v'' + m v - v^3 + v^5.
\end{equation*}
This mapping is of class $C^\infty$, we have $f(r, \nu)=0$, and 
$D_1 f(r,\nu) = A$, which is an isomorphism from $H^2_{\text{even}} $
to $L^2_{\text{even}}(\Er)$. Therefore, the implicit function theorem
applies, and $r$ is a $C^\infty$ 
function of $\nu$ with values in $H^2_{\textrm{even}}(\Er)$.

(ii)If we divide $A$ by $1-\nu$ we obtain a
Schr\"odinger operator whose differential 
part is constant, and whose potential part depends in a $C^\infty$ fashion on $\nu$,
according to part (i) of this proof; therefore, we can apply the
results of Kato, \cite{Kat} IV.3, which enable us to 
conclude.
 \end{proof}

Define three new functions $\sigma(x,\nu)$  and $\sigma_2(x,\nu)$ by
\begin{align}
S(x,\nu)&=-4\nu \frac{\partial R}{\partial \nu}
=\frac{\partial}{\partial L} R\bigl(x,4e^{-4L}\bigr)\label{eq:33}\\
\sigma(x,\nu)&=-4\nu \frac{\partial r(x,\nu)}{\partial \nu} =
\frac{\partial}{\partial L}r(x,4e^{-4L})
,\notag
\\
\sigma_2(x,\nu)&=-4\nu \frac{\partial \sigma(x,\nu)}{ \partial \nu} =
\frac{\partial^2}{\partial L^2}r(x,4e^{-4L}).
\notag
\end{align}
The following lemma states that $S+ R_x$, $\sigma +r_x$ and
$\sigma_2 + 2 \sigma_x + r_{xx}$
are very close to $0$; this is really a
consequence of the non linearity;  the benefit will
be even greater when we will show below that
$\sigma$ is an excellent approximation of the eigenfunction of $A$
 corresponding to its lowest eigenvalue.
 
\begin{lemma}\label{Lemma2.5}
We have the following estimates
for all $x\in \Er^+$, and all $\nu\in (0,1)$:
\begin{align}
\vabs{\frac{S(x,\nu)+R_x(x,\nu)}{R(x,\nu)}}&\le \frac{O(1)
\nu}{1-4\tilde R(x-L)/3},\label{eq:34}\\
\vabs{\sigma(x,\nu) +\frac{\partial r}{\partial x}(x,\nu)}
&\le O(1)\sqrt\nu\,e^{-x},\label{2.26}
\\						   
\intertext{and}					   
\vabs{\sigma_2(x,\nu) +2\frac{\partial \sigma}{\partial x}(x,\nu)+
\frac{\partial ^2 r}{\partial x^2}(x,\nu)}	   
&\le O(1)\sqrt\nu\,e^{-x},\label{2.27}
\end{align}
\end{lemma}

\begin{proof}The differentiation of~\eqref{2.18} with respect to $L$
gives
\begin{equation*}
\frac{1}{R\bigl(\cdot +L,4e^{-4L}\bigr)}\frac{1}{\partial L}
R\bigl(\cdot +L,4e^{-4L}\bigr) =-\frac{4\nu}{1-\nu}\frac{(1-2\tilde
R/3)^2}{1-4\tilde R/3 + 4\nu \tilde R^2/9},
\end{equation*}
which implies immediately~\eqref{eq:34}.
An analogous computation implies
\begin{equation}
\frac{\partial}{\partial L}\bigl[r(x+L,4e^{-4L})\bigr] =
r_x(x+L,4e^{-4L}) + \sigma(x+L,4e^{-4L}),\label{2.28}
\end{equation}
and
\begin{equation}
\begin{split}
\frac{\partial^2}{\partial L^2}\bigl[r(x+L,4e^{-4L})\bigr] &=
r_{xx}(x+L,4e^{-4L})\\
&\quad + 2\sigma_x(x+L,4e^{-4L})+\sigma_2(x+L,
4e^{-4L}).
\end{split}
\label{2.29}
\end{equation}
Differentiating~\eqref{2.18}, we find
\begin{equation*}
\frac{\partial}{\partial L} r(\cdot+L,4e^{-4L}) =\frac{2\nu \tilde r
(1 - 4 \tilde R/3)(1 - 2 \tilde R/3)^2}{
(1 - 4 \tilde R/3+ 4\nu \tilde R^2/9)^2},
\end{equation*}
and
\begin{equation*}
\begin{split}
\frac{\partial ^2}{\partial L^2} r(\cdot+L,4 e^{-4L})
&=-\frac{8 \nu \tilde r(1 - 4\tilde R/3)(1 - 2\tilde R/3)^2
(1 - 4 \tilde R/3 - 4\nu\tilde R^2/9)}{
(1 - 4 \tilde R/3+ 4\nu \tilde R^2/9)^3}\\
& -\frac{4\nu^2\tilde r
(1-4\tilde R/3)^2(1-2\tilde R/3)^4}{(1-4\tilde R/3 + 4\nu \tilde R^2/9)^4}
.
\end{split}
\end{equation*}
Observe that
\begin{equation*}
\frac{\nu \tilde r}{1-4\tilde R/3}
=\nu\sqrt{\frac{3}{4}}e^{-x}\bigl(1+e^{-2x}\bigr)^{1/2} ,
\end{equation*}
so that:
\begin{equation}
\forall x\ge -L, \quad \frac{\nu \tilde r}{1-4\tilde R/3}\le
O(1)\sqrt{\nu}e^{-x-L}; \label{eq:4}
\end{equation}
moreover we have also
\begin{equation}
\forall x\ge -L, \quad \frac{\nu }{1-4\tilde R/3} \le O(1); \label{eq:5}
\end{equation}
We use the information given by~\eqref{eq:4} and~\eqref{eq:5}, to
infer that
\begin{equation*}
\begin{split}
\forall x\ge -L,\quad&\vabs{\frac{\partial}{\partial L}
r(\cdot+L,4e^{-4L})}\le
O(1)\sqrt{\nu}e^{-x-L},\\&\vabs{\frac{\partial^2}{\partial L^2} 
r(\cdot+L,4e^{-4L})}\le O(1)\sqrt{\nu}e^{-x-L} ;
\end{split}
\end{equation*}
the conclusion of the lemma follows immediately.
\end{proof}

For later purposes, we define
\begin{equation*}
S(x+L,4e^{-4L})=\frac{\partial}{\partial L} R(x+L, 4e^{-4L}).
\end{equation*}

From here, we will obtain an upper estimate for the lower
bound of the spectrum of $A$:
 
\begin{lemma}\label{Lemma2.6}
We have
\begin{equation*}
\frac{(A\sigma, \sigma)}{(\sigma, \sigma)} \sim -\frac{3\nu}{2}.
\end{equation*}
\end{lemma}
 
\begin{proof}
Let us compute $A\sigma$ by differentiating the equation
\begin{equation*}
-m r'' + m r - r^3 + r^5=0
\end{equation*}
with respect to $L$; we obtain
\begin{equation}
A\sigma =\rho\label{2.32}
\end{equation}
where
\begin{equation}
\rho= \frac{4\nu}{1 - \nu}(r^5 - r^3).\label{2.33}
\end{equation}
 
To prove our assertion, we will have to calculate
$(\sigma, \sigma)$ and $(\sigma, \rho)$.
From Lemma \ref{Lemma2.5}, we can see that
\begin{equation*}
\int\sigma^2(x, \nu)\,dx = 2\int_0^\infty r_x(x,\nu)^2\,dx
+O(\sqrt\nu),
\end{equation*}
and from estimate~\eqref{2.12}, we have
\begin{equation*}
2\int_0^\infty r_x(x,\nu)^2\,dx = 2 \int_{-L}^\infty \tilde
r'(x)^2\,dx + O(\sqrt\nu);
\end{equation*}
but relations~\eqref{2.15} and~\eqref{2.21} imply that 
\begin{equation*}
\vabs{\tilde r'}^2 = -\frac12\tilde R'\left(1-\frac{4\tilde R}{3}\right);
\end{equation*}
we infer from relation~\eqref{eq:2} that
\begin{equation}
2\int_{-L}^\infty \vabs{\tilde r'}^2 \, dx = \frac38 +O\bigl(\sqrt{\nu}\bigr),\label{eq:17}
\end{equation}
so that
\begin{equation}
(\sigma, \sigma) = \frac{3}{8} + O(\sqrt \nu).\label{2.34}
\end{equation}
  
On the other hand,
\begin{equation}
\begin{split}
(r^5 - r^3,\sigma)&=- 2 \int_0^\infty(r^5 - r^3)r_x\,dx
+ O(\sqrt \nu)\\
&=2\left(\frac{R(0, \nu)^3}{6} - \frac{R(0,\nu)^2}{4}\right) + O(\sqrt \nu)\\
&=-\frac{9}{64} + O(\sqrt \nu).
\end{split}\label{eq:3}
\end{equation}
 
Thus, we obtain from~\eqref{2.33} and~\eqref{eq:3} the relation
\begin{equation}
(A\sigma,\sigma)=-\frac{9\nu}{16} + O(\nu^{3/2})  .\label{2.35}
\end{equation}
The conclusion of the lemma is a direct consequence of~\eqref{2.34}
and~\eqref{2.35}.
\end{proof}
 
\begin{theorem}\label{Theorem2.7}
There exists a constant $K$ such that for all
small enough $\nu$ the spectrum of $A$ is partitioned in its intersection
with the interval $[-K,K]$ and its intersection with interval $[2K,
+\infty)$.
\end{theorem}
 
\begin{proof}Define
\begin{equation*}
\tilde V(x)=\frac3{16}-3\tilde R +5\tilde R^2.
\end{equation*}
It is clear that $\tilde r'$ satisfies the Schr\"odinger equation
\begin{equation*}
-\tilde r''' + \tilde V \tilde r'=0.
\end{equation*}
As $\tilde r'$ is strictly negative for all $x$, this means that $0$
is the lower bound of the spectrum of the operator $\tilde A$ defined by
\begin{equation*}
D(\tilde A)=H^2(\Er), \tilde A v =-v''+\tilde V v.
\end{equation*}
As $\tilde V(-\infty)=3/4$ and $\tilde V(+\infty)=3/16$, the essential
spectrum of $\tilde A$ is included in $[3/16, +\infty)$, and thus the
second eigenvalue of $\tilde A$ is some number $\Lambda>0$. 

We use now the twisting trick of E. B. Davies \cite{daviestwist}: define
the operators
\begin{equation*}
\tilde \calj=\begin{pmatrix}
-\partial^2 + \tilde V(\cdot-L) &0\\
0& -\partial^2 +\tilde V(-\cdot -L)
\end{pmatrix}, \quad \calj=\begin{pmatrix}
-\partial^2 + 3/4&0\\
0&A
\end{pmatrix};
\end{equation*}
let $\eta$ be an infinitely differentiable function from $\Er$ to
$[0,\pi/2]$ which vanishes for $x\le -1$ and is equal to $\pi/2$ for
$x\ge 1$ and let $\eta_L=\eta(\cdot/L)$. Define now a unitary
transformation in $L^2(\Er)^2$ by its matrix 
\begin{equation*}
U_L=\begin{pmatrix}
\cos\eta_L&-\sin\eta_L\\\sin\eta_L&\cos\eta_L
\end{pmatrix}.
\end{equation*}
We compute $U_L \calj U_L^*$ and we find that
\begin{equation*}
\begin{split}
&U_L \frac{\partial^2}{\partial x^2} U_L^*=
\frac{\partial^2}{\partial x^2} +
-\frac{\eta'(\cdot/L)}{L}\begin{pmatrix}
0&-1\\1&0
\end{pmatrix}\frac{\partial}{\partial x}
 -\frac{\eta''(\cdot/L)}{L^2}
\begin{pmatrix}
0&-1\\1&0
\end{pmatrix} -\frac{\bigl(\eta'(\cdot/L)\bigr)^2}{L^2}\\
\intertext{and}
&U_L\begin{pmatrix}
V&0\\0&3/4
\end{pmatrix} U_L^*=\begin{pmatrix}
V\cos^2\eta_L  + \bigl(3\sin^2\eta_L\bigr)
/4&
\bigl(\cos\eta_L\sin\eta_L\bigl) (V-3/4)\\
\bigl(\cos\eta_L\sin\eta_L\bigr) (V-3/4)&
V\sin^2\eta_L  + \bigl(3\cos^2\eta_L\bigr)
/4
\end{pmatrix}.
\end{split}
\end{equation*}
An elementary calculation shows now that 
\begin{equation*}
U_L\calj U_L^* = \tilde\calj + B_L\frac{\partial}{\partial x} + C_L
\end{equation*}
where $B_L$ and $C_L$ are matrix valued functions whose norm tends to
$0$ as $L$ tends to infinity. It is plain that the resolvent sets of
$\tilde \calj$ and $\tilde A$ are identical.

If $E$ and $F$ are Banach spaces, $\call(E,F)$ is the space of linear
bounded 
operators from $E$ to $F$ and $\vnorm[\call(E,F)]{\>\>}$ the
corresponding norm; if $E=F$, we abbreviate $\call(E,F)$ to $\call(E)$.
If $E$ is the space $L^2(\Er)^n$ we abbreviate $\call(E)$ to $\call$.

Let $z$ belong to the resolvent set of $\tilde A$; then for every $g$
in $L^2(\Er)^2$, 
$(\tilde \calj -z)^{-1}$ belongs to $H^2(\Er)^2$, and the function 
$z\mapsto \norm[\call(L^2)]{(\partial^2/\partial x^2)(\tilde \calj
-z)^{-1}}$ is continuous on the resolvent set of $\tilde \calj$. 

It is equivalent to solve
\begin{equation*}
\left(\tilde \calj +B_L\frac{\partial}{\partial x} + C_L -z\right)u=f
\end{equation*}
and
\begin{equation*}
g+\left(B_L\frac{\partial}{\partial x} + C_L\right)(\tilde \calj
-z)^{-1} g=f;
\end{equation*}
the previous considerations show that for all $z$ in the resolvent set
of $\tilde \calj$, we can find $L(z)$ such that for all $L\ge L(z)$,
$z$ is also in the resolvent set of $\calj$. In particular, if we take
$K=\Lambda/3$, we can find $\bar L$ such that for all $L\ge \bar L$,
the segment $(K, 2K)$ is included in the resolvent set of $\calj$. 
In particular, we see also that for all $L\ge \bar L$, the generalized
eigenspace of $A$ relative to the interval $[-\min V, K]$ contains
exactly two eigenvalues; the minimum of $V$ is equal to $-21/80 +
O(\nu)$; the above considerations show that these two eigenvalues tend
to $0$ as $\nu$ tends to $0$.
\end{proof}

Let us prove now that $\lambda$ is very close to
$(A\sigma,\sigma)/(\sigma,\sigma)$. 
 
Denote by $\Pi$ the projection onto the sum of the eigenspaces relative to
$0$ and $\lambda$. We have
\begin{equation}
\Pi = \frac{1}{2 \pi i }\int_\gamma (\zeta - A)^{-1}\,d\zeta,\label{2.39}
\end{equation}
where $\gamma $ is a circle of radius $3K/2$, which is traveled once in the positive
direction.
 
\begin{theorem}\label{Theorem2.8}
We have the following estimates:
\begin{equation}
|\sigma - \Pi\sigma|_{H^2(\Er)} = O(\nu \sqrt L),\label{2.40}
\end{equation}
and
\begin{equation}
\lambda = 
-\frac{3\nu}{2}+O\bigl(\nu^{3/2}\bigr).\label{2.41}
\end{equation}
In particular, we can choose $s = \Pi\sigma$ as an eigenvector of
$A$ relative to the eigenvalue $\lambda$, as soon as $\nu$ is small
enough.
Moreover, if $\hat \sigma$ is the unique solution of 
\begin{equation}
\Pi\hat\sigma=0, \quad A\hat \sigma =(\one -\Pi)(\lambda \sigma -\rho)\label{eq:28}
\end{equation}
then
\begin{equation}
\vabs[H^2(\Er)]{\Pi\sigma -\sigma -\hat \sigma} =O\bigl(\nu^2 \sqrt L\,\bigr).\label{eq:27}
\end{equation}
\end{theorem}
 
\begin{proof}
We observe that for $|\zeta| = 3K/2$, $\zeta - A$ is an isomorphism from $H^2(\Er)$
to $L^2(\Er)$ which transforms even functions into even functions,
and odd functions into odd functions. Therefore, since $\sigma$ is even by
construction, $\Pi \sigma$ is even too.
 
We rewrite the identity $A\sigma = \rho$ as
\begin{equation*}
\sigma = \frac{\rho}{\zeta} +\frac{(\zeta - A)\sigma}{\zeta}.
\end{equation*}
Hence, by integration along $\gamma$,
\begin{equation}
\Pi \sigma = \sigma + \frac{1}{2 \pi i}\int_\gamma
(\zeta- A)^{-1}\zeta^{-1}\rho\,d\zeta.\label{eq:25}
\end{equation}
Since $|\rho| = O(\nu \sqrt L)$, we obtain
\begin{equation}
|\sigma -\Pi\sigma| = O(\nu \sqrt L),\label{2.42}
\end{equation}
  
If $\nu$ is small enough,~\eqref{2.42} implies that $\Pi\sigma$ does not
vanish identically. 
Thus, by construction it is an eigenfunction of $A$ relative to the eigenvalue
$0$ or to the eigenvalue $\lambda$. But we know that $\Pi\sigma$ is even;
thus, it has to be an eigenfunction of $A$ relative to $\lambda$, and
\begin{equation*}
\lambda = \frac{(A\Pi \sigma, \Pi\sigma )}{(\Pi\sigma, \Pi\sigma)}.
\end{equation*}
But
\begin{equation*}
(A\Pi\sigma, \Pi\sigma) = (\rho, \Pi\sigma) = (\rho, \sigma) + (\rho,
\Pi\sigma - \sigma) = (\rho, \sigma) +O(L \nu^2).
\end{equation*}
and, on the other hand
\begin{equation*}
(\Pi \sigma, \Pi\sigma) = (\sigma, \Pi\sigma)
=(\sigma, \sigma) + (\Pi\sigma-\sigma, \sigma)
=(\sigma, \sigma) + O(\nu \sqrt L).
\end{equation*}
This proves~\eqref{2.41}.
 
In order to obtain~\eqref{2.40} which holds in $H^2$ norm,
while~\eqref{2.40} holds in $L^2$ norm, we subtract the relation
$A\Pi \sigma = \lambda\Pi\sigma$ from the relation $A\sigma = \rho$
and we obtain
\begin{equation*}
m\bigl|\sigma_{xx} - (\Pi\sigma)_{xx}\bigr| =
\bigl|\rho - \lambda \Pi \sigma
-V(\sigma - \Pi \sigma)\bigr| = O(\nu \sqrt L).
\end{equation*}

For the last assertion of the theorem, we remark thanks to
Cauchy's theorem, the last term in~\eqref{eq:25} can be rewritten as
\begin{equation*}
\frac{1}{2 \pi i}\int_\gamma
(\zeta- A)^{-1}\zeta^{-1}\rho\,d\zeta=\frac{1}{2i\pi} \int_\gamma
(\zeta-A)^{-1} \zeta^{-1}(\one -\Pi)\rho\, d\zeta=\tilde \sigma
\end{equation*}
and $\tilde \sigma$ is the unique solution of
\begin{equation}
\Pi\tilde \sigma = 0, \quad A\tilde \sigma =(\one -\Pi)\rho.\label{eq:26}
\end{equation}
Subtracting~\eqref{eq:26} from~\eqref{eq:28}, we find that
\begin{equation*}
\Pi(\tilde \sigma -\hat \sigma)=0, \quad A(\tilde\sigma -\hat \sigma)=\lambda\bigl(\sigma-\Pi\sigma).
\end{equation*}
and this implies immediately estimate~\eqref{eq:27}, which concludes
the proof.
\end{proof}

We need also the equation satisfied by $\sigma_2$. If
we differentiate~\eqref{2.32} with respect to $L$, we obtain
\begin{equation}
A\sigma_2 + (20r^3 - 6r)\sigma^2=\rho_2,\label{2.43}
\end{equation}
where
\begin{equation*}
\rho_2 =-4\nu\frac{\partial \rho}{\partial \nu}
+\frac{3\nu}{4}(\sigma''-\sigma)
=-\frac{16\nu}{(1-\nu)^2}(r^5 - r^3) +\frac{8\nu}{1-\nu}(5r^4
-3r^2)\sigma.
\end{equation*}
Therefore,
\begin{equation}
|\rho_2|= O\bigl(\sqrt L \,\nu\bigr).\label{2.44}
\end{equation}
Moreover an analogous and straightforward argument shows that
\begin{equation}
\frac{\partial \lambda}{\partial L} =0(\nu), \quad
\vabs{\frac{\partial s}{\partial L}
-\sigma_2}=O\bigl(\sqrt L\, \nu\bigr).\label{eq:24}
\end{equation}

The final information we obtain are asymptotics on  $\partial
\lambda/\partial L$ on $\partial s/\partial L$:

\begin{lemma}\label{Lemma2.9} The following estimates hold:
\begin{align}
&\frac{\partial \lambda}{\partial L}=6\nu +O\bigl(\nu^{3/2}\bigr), \label{eq:30}\\
&\vabs[H^2(\Er)]{\frac{\partial
s}{\partial L}-\sigma_2}=O\bigl(\nu\sqrt L\,\bigr). \label{eq:31}
\end{align}
\end{lemma}

\begin{proof}
For the first estimate, we differentiate the relation
\begin{equation*}
A\Pi\sigma =\lambda \Pi\sigma
\end{equation*}
with respect to $L$; we multiply scalarly the result
\begin{equation*}
\frac{\partial A}{\partial L} \Pi\sigma +A\frac{\partial
\Pi\sigma}{\partial L}= \frac{\partial\lambda}{\partial L}\Pi\sigma +
\lambda \frac{\partial \Pi\sigma}{\partial L}
\end{equation*}
by $\Pi\sigma$ and we obtain
\begin{equation*}
\vabs{\Pi\sigma}^2 \frac{\partial
\Pi\sigma}{\partial L}=\left(\frac{\partial A}{\partial
L}\Pi\sigma,\Pi\sigma\right). 
\end{equation*}
Since
\begin{equation}
\frac{\partial A}{\partial L}=\frac{4\nu}{1-\nu} \bigl(A-3r^2
+5r^4\bigr) -6r\sigma +20 r^2 \sigma,
\end{equation}
we may write
\begin{equation*}
\left(\frac{\partial A}{\partial
L}\Pi\sigma,\Pi\sigma\right)=4\nu \bigl((3r^2-r^4)\sigma,\sigma\bigr)
+\bigl((20r^3 -6r)\sigma^2,\sigma +2\hat \sigma\bigr) +O(\nu^2 L).
\end{equation*}
Thanks to~\eqref{2.43},
\begin{equation*}
\begin{split}
&\bigl((20r^3 -6r)\sigma^2,\sigma +2\hat \sigma\bigr)=\bigl(A\sigma_2
-\rho_2,\sigma\bigr)\\
&\quad=(\sigma_2,\rho)-(\rho_2,\sigma) +(\sigma_2,\lambda\sigma-\rho)
+O\bigl(\nu^{3/2}\bigr);
\end{split}
\end{equation*}
here we have used~\eqref{eq:28} and~\eqref{2.40}.  Therefore,
\begin{equation*}
\left(\frac{\partial A}{\partial
L}\Pi\sigma,\Pi\sigma\right)=4\nu
\bigl((3r^2-r^4)\sigma,\sigma\bigr)+2\lambda \sigma_2, \sigma)-
\frac{\partial}{\partial L}(\rho, \sigma) + O(\nu^{3/2}).
\end{equation*}
At this point it is clear that there exists a number $c$ such that
\begin{equation}
\frac{\partial \lambda}{\partial L}=c\nu + O\bigl(\nu^{3/2}\bigr),\label{eq:29}
\end{equation}
but the explicit calculation of $c$ is tedious. It can be avoided by
integrating~\eqref{eq:29} with respect to $L$: we see that
\begin{equation*}
\lambda = -c\nu/4 + O\bigl(\nu^{3/2}\bigr),
\end{equation*}
and we infer from~\eqref{2.41} that $c=6$.

The other asymptotics are obtained as follows: we differentiate the
relation
\begin{equation*}
\Pi\sigma=\int_\gamma (\zeta-A)^{-1}\sigma\, d\sigma
\end{equation*}
with respect to $L$: 
\begin{equation*}
\frac{\partial\Pi\sigma}{\partial L} = \frac{1}{2i\pi}\int_\gamma
(\zeta -A)^{-1} \Bigl[ \frac{\partial A}{\partial L}(\zeta -A)^{-1}
\sigma +\sigma_2\Bigr]\, d\zeta.
\end{equation*}
We deduce from~\eqref{eq:27} that for $\zeta\in\gamma$
\begin{equation*}
\vabs[H^2(\Er)]{(\zeta-A)^{-1}\sigma -(\zeta -\lambda)^{-1}}
=O\bigl(\nu\sqrt L\, \bigr);
\end{equation*}
therefore, it suffices to find an asymptotics for
\begin{equation*}
\frac{1}{2i\pi}\int\gamma (\zeta
-A)^{-1}(\zeta -\lambda)^{-1}\Bigl[\frac{4\nu}{1-\nu}\bigl(\rho +(3r^2 -5r^4)\sigma +
\rho_2 -A\sigma_2 +(\zeta -\lambda)\sigma_2\Bigr]\, d\zeta.
\end{equation*}
But
\begin{equation*}
(\zeta -A)^{-1}(\zeta -\lambda)^{-1} \bigl(-A\sigma_2 +(\zeta
-A)\sigma_2\bigr)= (\zeta -\lambda)^{-1} -\lambda (\zeta
-A)^{-1}(\zeta -\lambda)^{-1}
\end{equation*}
and it is clear now that~\eqref{eq:29} holds.
\end{proof}

\section{The second equation}\label{sec:second-equation}

Let us first observe that we need only consider the case $\alpha >0$:
the mapping 
$\alpha\mapsto -\alpha$ transforms a solution into its conjugate.

There are two different ways of writing the full problem, each with a
different scaling. 
The unscaled problem is
\begin{equation}
-\frac{3}{16}(1 - \nu) u'' + \frac{3}{16}(1 - \nu) u -
(1 + i\alpha) |u|^2 u + |u|^4 u+i\omega u=0.\label{3.1}
\end{equation}
Therefore, we define  
\begin{equation*}
F(u,\omega, \alpha,\nu) = -\frac{3}{16}(1 - \nu) u'' +
\frac{3}{16}(1 - \nu) u -
(1 + i\alpha) |u|^2 u + |u|^4 u+i\omega u.
\end{equation*}
Whenever necessary, we will use the components of $u$, denoting
 them by $u_1$ and $u_2$, so that $u= u_1 + i u_2$, and we define 
similarly the components of $F$, $F_1$ 
and $F_2$; therefore, we have   
\begin{equation*}
F = F_1 + i F_2.
\end{equation*}
Let $\alpha=\sqrt\varepsilon$, $u_1=\xi$, $u_2=\eta\sqrt \varepsilon$, $\omega=\tau
\sqrt\varepsilon$.  
The scaled problem is
\begin{equation}
G(\xi, \eta,\tau,\varepsilon,\nu)=0,\label{3.4}
\end{equation}
where we define
\begin{equation}
G(\xi,\eta,\tau,\varepsilon,\nu)=
F_1(\xi+i\sqrt\varepsilon\,\eta,\sqrt\varepsilon\tau,\sqrt\varepsilon,\nu) +
\frac{i}{{\sqrt \varepsilon}}F_2(\xi+i\sqrt\varepsilon\,\eta,
\sqrt\varepsilon\tau,\sqrt\varepsilon,\nu).\label{3.5}
\end{equation}
 
More explicitly,
\begin{align*}
G_1(\xi,\eta,\tau,\varepsilon,\nu)&=-\tau\varepsilon\eta-
\frac{3}{16}(1-\nu)\xi'' +\frac{3}{16}(1-\nu)\xi\notag\\
&\quad -(\xi-\varepsilon\eta)(\xi^2 + \varepsilon\eta^2)
+(\xi^2 + \varepsilon \eta^2)^2\xi,
\\
G_2(\xi,\eta,\tau,\varepsilon,\nu)&=\tau \xi
-\frac{3}{16}(1-\nu)\eta'' +\frac{3}{16}(1-\nu)\eta\\
&\quad -(\xi + \eta)(\xi^2+ \varepsilon \eta^2)+(\xi^2 + \varepsilon
\eta^2)^2 \eta.
\end{align*}

Clearly,~\eqref{3.4} is equivalent to~\eqref{3.1}. We will solve
the existence problem 
looking at the $G$ formulation, and the stability problem looking at the $F$
 formulation.
 
We have already determined a solution of
\begin{equation*}
G_1(\xi, \eta,\tau,0,\nu)=0.
\end{equation*}
it is given by $\xi=r$, and $\eta$ and $\tau$ are arbitrary.

In next lemma, we state the sense in which we can find a solution for
the second equation 
at $\varepsilon = 0$:

\begin{lemma}\label{Lemma3.1}
There exists a unique $\theta\in \Er$ and a unique $q$ in $H^2(\Er)$
and orthogonal to $r$
which satisfy
\begin{equation}
G_2(r,q,\theta,0,\nu)=0\label{3.7}
\end{equation}
Moreover, $q $ is even.
\end{lemma}

\begin{proof}
The second equation at $\varepsilon = 0$ can be written
\begin{equation*}
\tau\xi -m\eta'' +m\eta -\xi^3 -\xi^2\eta +\xi^4 \eta=0.
\end{equation*}
Therefore, if we let $\xi=r$, we have to find $q$ and $\theta$ such
that~\eqref{3.7} holds, i.e. 
\begin{equation}
Bq=-m q'' + m q - r^2 q + r^4 q=r^3 - \theta r.\label{3.8}
\end{equation}

Define an unbounded operator $B$ in $L^2(\Er)$ by
\begin{equation*}
D(B)=H^2(\Er), \quad B u = -m u'' + W u,
\quad W = m - r^2 + r^4.\label{3.9}
\end{equation*}
Clearly, $B$ is self adjoint.
By construction
\begin{equation*}
B r =0.
\end{equation*}
 
The lower bound of the essential spectrum of $B$ is $m>0$;
since $r$ is positive, the lower bound of the spectrum of $B$
is zero, which is a simple isolated eigenvalue. Therefore,
solving~\eqref{3.8} will be possible if and only if its right hand side
$r^3 - \theta r$
is orthogonal to $\ker B= \Er r$.
This defines $\theta$ thanks to
\begin{equation*}
\theta = \frac{\int r^4 \,dx}{\int r^2\,dx}.
\end{equation*}
Moreover, we define completely $q$ by requiring
\begin{equation}
(q,r)=0.\label{3.12}
\end{equation}

It is clear that the $q$ we obtained belongs to $H^2(\Er)$.
Since the equation is invariant by the transformation $x\mapsto -x$,
$q$ has to be even.
\end{proof}

We will need precise asymptotic information, which we proceed to give
now:

\begin{lemma}\label{Lemma3.2}
The number $\theta$ is a $C^1$ function of $\nu\in (0,1)$; we
have the following asymptotics for $\theta$ and $\partial \theta/\partial L$:
\begin{align}
\theta&=\frac{3}{4} -   \frac{3}{8L} +
O(\sqrt\nu).\label{3.13}
\\
\theta_1&=\frac{\partial \theta}{\partial L} =\frac{3}{8L^2} +
O(L^{-1}\sqrt\nu).\label{3.14}
\end{align}
\end{lemma}

\begin{proof}
Let us compute the difference between $3/4$ and $\theta$:
\begin{equation*}
\frac{3}{4} - \theta =\frac{\int_0^\infty (3r^2/4 - r^4)\,dx}{
\int_0^\infty r^2\,dx}.
\end{equation*}
But an explicit computation and~\eqref{2.11} give
\begin{equation*}
\int_0^\infty R\,dx= \frac{3L}{4} +\frac{3}{8}\ln \frac{4}{3} +
O(\sqrt\nu),
\end{equation*}
that is
\begin{equation}
\int_0^\infty R\,dx = \frac{3L}{4} +O(\sqrt\nu).\label{3.16}
\end{equation}
In the same fashion,
\begin{align*}
\int_0^\infty (3R/4 - R^2)\,dx &= \frac{3}{4} \int_0^\infty R(1 - 4R/3)\,dx\\
&=\frac{3}{4}\int_0^\infty \tilde R(1 - 4\tilde R/3)\,dx + O(\sqrt\nu)\\
&=\frac{3}{8}\tilde R(-L) + O(\sqrt\nu)
=\frac{9}{32} +O(\sqrt\nu).
\end{align*}
Therefore,
\begin{equation*}
\frac{3}{4} - \theta = \frac{3}{8L} + O(\sqrt\nu),
\end{equation*}
which implies immediately~\eqref{3.13}.

Since $r$ is a $C^\infty$ function of $\nu\in(0,1)$ with values in
$H^2(\Er)$, and it is 
integrable as well as $\sigma$, $\theta$ is a $C^1$ function of $\nu$. Let us
compute $\partial \theta/\partial L$; we have
\begin{equation}
\frac{\partial \theta}{\partial L} 
=\left(\int_0^\infty r^2\,dx \int_0^\infty 4r^3\sigma\,dx -
\int_0^\infty r^4 \,dx \int_0^\infty 2r\sigma\,dx\right)\left(\int_0^\infty
r^2\,dx\right)^{-2}.\label{3.17}
\end{equation}
 But
\begin{equation*}
\int_0^\infty 4r^3\sigma\,dx = -\int_0^\infty 4r^3 r'\,dx + O(\sqrt\nu)
=R^2(0) + O(\sqrt\nu)
=\frac{9}{16} + O(\sqrt \nu).
\end{equation*}
Similarly,
\begin{equation*}
\int_0^\infty 2 r \sigma\,dx = - \int_0^\infty 2 r r'\,dx + O(\sqrt\nu)
=R(0) + O(\sqrt \nu)
=\frac{3}{4} +O(\sqrt \nu).
\end{equation*}
Therefore, the numerator of~\eqref{3.17} is equal to  
\begin{equation*}
\begin{split}
&\frac{9}{16} \int_0^\infty R\,dx - \frac{3}{4} \int_0^\infty R^2\,dx +
O(L\sqrt\nu) =\frac{9}{16} \int_{-L}^\infty \tilde R(1 - 4\tilde R/3)\,dx +
O(L\sqrt\nu)\\ &=\frac{27}{128} +O(L\sqrt\nu).
\end{split}
\end{equation*}
With the help of~\eqref{3.16}
, we obtain ~\eqref{3.14}
.
\end{proof}

Let us find now asymptotics for $q$ and related quantities.
Let 
\begin{equation}
\phi = \frac{q}{r}.\label{3.18}
\end{equation}
Then $\phi$ satisfies the ordinary differential equation:
\begin{equation*}
-m(\phi'' r + 2\phi' r') = r^3 -\theta r.
\end{equation*}
Multiplying by $r$ we get
\begin{equation}
m(\phi' r^2)' = \theta r^2 - r^4.\label{3.19}
\end{equation}
Therefore,
\begin{equation*}
m (\phi'r^2)(x) = a+ \int_{-\infty}^x(\theta R - R^2)(y)\,dy.
\end{equation*}
But $a$ must vanish: $q'$ belongs to $L^2(\Er)$ and is
equal to $\phi r' + 
\phi' r$; but
\begin{equation*}
r'= - r  \sgn(x) \sqrt{\Psi(r)/m},
\end{equation*}
hence
\begin{equation*}
\phi r' = - q\sgn(x) \sqrt{\Psi(r)/m},
\end{equation*}
so that $\phi ' r = q' - \phi r'$ belongs to $L^2(\Er)$. This implies
that $\phi' r^2$ 
belongs to $L^2(\Er)$ and hence $a$ vanishes.
Therefore
\begin{equation}
m \phi'(x) = \frac{1}{R(x)}\int_{-\infty}^x (\theta - R)
R\,dy= -\frac{1}{R(x)} \int_x^\infty  (\theta - R)
R\,dy. \label{3.20}
\end{equation}
                                                        
\begin{lemma}\label{Lemma3.3}
The following asymptotics hold for $x\ge 0$:
\begin{align}
&\phi'(x)=-2 -\frac{3\ln(1-4\tilde R(x-L)/3)}{4L\tilde R(x-L)}
+O(\sqrt\nu/L).\label{3.21}
\\
&\phi(0)=\frac{L}{3} + O(1),\label{3.22}\\
&\phi(L)=-\frac{2L}{3} +
O(1)\label{3.22a}
\\
&|q|_{L^1(\Er)} = O(L^2),\quad |q|_{L^\infty(\Er)} =O(L), \quad
|q| = O(L^{3/2}).\label{3.23}
\end{align}
\end{lemma}

\begin{proof}
From~\eqref{3.20}, we obtain
\begin{equation}
m \phi'(x) = -\frac{1}{R(x)}(\theta - 3/4)\int_x^\infty R(y)\,dy -
\frac{1}{R(x)}\int_x^\infty (3/4 - R)R\,dy.\label{3.24}
\end{equation}
For $x>0$, we have thanks to~\eqref{2.11}
\begin{equation*}
\begin{split}
&\int_x^\infty(3/4 - R)R\,dy = \int_{x - L}^\infty(3/4 - \tilde
R)\tilde R\,dy 
+O(\sqrt \nu e^{-2x})\\
&\quad=\frac{3}{8}\tilde R(x-L) + O(\sqrt\nu e^{-2x}).
\end{split}
\end{equation*}
But relation~\eqref{2.11} implies also that for $x>0$
\begin{equation}
\vabs{\frac{R(x)}{\tilde R(x-L)} - 1} \le O(1)\sqrt\nu\,.\label{eq:35}
\end{equation}
Therefore,
\begin{equation}
\frac{1}{R(x)}\int_{-\infty}^x (3/4 - R)R\,dy
=\frac{3}{8} + O(\sqrt\nu).\label{eq:36}
\end{equation}

On the other hand,
\begin{equation*}
\begin{split}
\int_x^\infty R\, dy &=\int_{x-L}^\infty \tilde R \, dy +\sqrt\nu\,
e^{-2x}O(1) \\
&=\frac{3}{8} \ln \bigl(1+4e^{2L-2x}/3\bigr) + \nu
O(1)\frac{4\tilde R(x-L)}{1-4\tilde R(x-L)/3}.
\end{split}
\end{equation*}
Therefore, we have also
\begin{equation}
\frac{1}{R(x)}\int_x^\infty R\, dy =-\frac{3}{8} \frac{\ln
\bigl(1-4\tilde R(x-L)/3\bigr)}{\tilde R(x-L)} +O\bigl(\sqrt \nu\,
\bigr).\label{eq:37} 
\end{equation}
Thus we have obtained the asymptotic for $x\ge 0$
\begin{equation}
\int_x^\infty R(\theta -R) \, dy =R\left(\frac{3}{8}
-\frac{3}{8}\left[\theta -\frac{3}{4}\right] \frac{\ln\bigl(1-4\tilde
R(x-L)/3\bigr)}{\tilde R(x-L)} +O\bigl(\sqrt\nu\, \bigr)\right).\label{eq:38}
\end{equation}
Relation~\eqref{3.21} is an immediate consequence of~\eqref{eq:38}.

In order to have an asymptotic for $\phi(0)$, we will give another
asymptotic for $\phi'$, 
which will be much less precise, but more usable at this level of the discussion;
nevertheless, the precise asymptotic~\eqref{3.21} will be useful later.

For $x\ge L$, 
\begin{equation*}
\ln\bigl(1-4\tilde R(x-L)/3\bigr) = -\frac{4\tilde R(x-L)}{3} + \tilde
R(x-L)^2 O(1).
\end{equation*}
Therefore, for $x\ge L$,
\begin{equation*}
\phi'(x) = -2 +\frac{1}{L} + O(1)\frac{\tilde R(x-L)}{L}+O(\sqrt{\nu}).
\end{equation*}
In the same fashion, for $0\le x\le L$,
\begin{equation*}
-\ln\bigl(1- 4\tilde R(x-L)/3\bigr)
=\ln\bigl(1+e^{2(L-x)}\bigr)=2(L-x)+O(1)e^{2(x-L)}
\end{equation*}
\begin{equation}
\phi'(x) =
\begin{cases}
-2x/L +O(1/L)e^{2(x-L)} & \text{if $0\le x\le L$,}\\
-2 + O(1/L) & \text{if $x\ge L$.}
\end{cases}
\label{3.26}
\end{equation}
By integration in $x$, we get now                                    
\begin{equation}
\phi(x) - \phi(0)=
\begin{cases}
-x^2/L + O(1), & \text{if $0\le x\le L$,}\\
L- 2x + O(1)+(x-L)O(1/L), & \text{if $x\ge L$.}
\end{cases}
\label{3.27}
\end{equation}
The value of $\phi(0)$ is determined by the condition that $q$ should
be orthogonal to 
$r$ (see Lemma~\ref{Lemma3.1}). 
Therefore,
\begin{align*}
\phi(0)\int_0^\infty R(x)\,dx& = \int_0^L\left(\frac{x^2}{L}
+O(1)\right)R(x)\,dx \\
&\quad + \int_L^\infty(-L+2x)R(x)\,dx + \int_L^\infty \biggl(O(1) +
\frac{x-L}{L}O(1)\biggr)R(x)\,dx. 
\end{align*}
At this point we observe that
\begin{align*}
&\int_0^\infty R \,dx = \frac{3L}{4} + O(1), \quad \int_0^\infty \bigl(x^2/L +
O(1)\bigr) R(x)\,dx =\frac{L^2}{4} + O(L),\\
&\int_L^\infty (2 x- L) R(x) \,dx = O(L), \quad \int_L^\infty \bigl(
O(1) +O(1)(x-L)/L \bigr) R(x)\,dx = O(1).
\end{align*}
This yields the first relation in~\eqref{3.22}. The second relation in~\eqref{3.22a}
is a consequence of the first together with~\eqref{3.27}.

Estimates~\eqref{3.23}
are a direct consequence of the definition ~\eqref{3.18}
of $\phi$, of the equivalents~\eqref{3.26} and of~\eqref{3.22} and~\eqref{3.22a}. \end{proof}

Since $0$ is an isolated eigenvalue of $B$ for all $\nu$ in $(0,1)$ and since
$(1-\nu)^{-1} W$ is an infinitely differentiable function of $\nu$ with values in
the space $C^2_b$ of functions bounded as well as their first two
derivatives in space, 
we can see that $q$ is a $C^\infty$ mapping from $\nu\in(0,1)$ to
$H^2(\Er)$, under the normalizing condition~\eqref{3.12}.

Differentiation of $q$ with respect to $L$ will be useful too; define
\begin{equation*}
q_1 = \frac{\partial}{\partial L} q(x, 4e^{-4L}).
\end{equation*}
If we differentiate~\eqref{3.8} with respect to $L$ we obtain
\begin{equation}
B q_1 +(4r^3 - 2r)q\sigma =\hat \rho +3r^2 \sigma -\theta\sigma -\theta_1
r,\label{3.29}
\end{equation}
where
\begin{equation}
\hat \rho = \frac{4\nu}{1-\nu}(r^4q -r^2 q -r^3 +\theta
r).\label{3.30}
\end{equation}
We recall that $\theta_1$ has been defined
at~\eqref{3.14}

We will need an estimate over $\phi''$. From~\eqref{3.19}, we have
\begin{equation*}
m\phi'' = - R + \theta - 2m \phi'r'/r.
\end{equation*}
Assume that $x\le 0$ and Write
\begin{equation*}
a=\frac{\phi'(x)-2}{2};
\end{equation*}
then, thanks to~\eqref{eq:6}, we have
\begin{equation*}
-2m\phi'r'/r=-2(1+a)mR'=-(1-\nu)(1+a)\left(\frac34 -R +O(\sqrt{\nu}\right),
\end{equation*}
so that
\begin{equation*}
\theta-R - -2m\phi'r'/r=\theta-\frac34 -\left(\frac34 -R\right)a
+O\bigl(\sqrt{\nu}\bigr). 
\end{equation*}
But the following expressions are bounded independently of $\nu$ for
$x\le 0$:
\begin{equation*}
\left(\frac34 -\tilde R(\vabs{x}-L)\right)\frac{\ln(1-4\tilde
R(\vabs{x}-L)/3)}{\tilde R(\vabs{x}-L)}\text{ and } \sqrt{\nu}e^x 
\frac{\ln(1-4\tilde
R(\vabs{x}-L)/3)}{\tilde R(\vabs{x}-L)};
\end{equation*}
therefore, the term $a(3/4 -R)$ is an $O(1/L)$, and we conclude with
the help of~\eqref{3.13} that 
\begin{equation}
|\phi''|_\infty = O(L^{-1}).\label{3.31}
\end{equation}

We will give now two estimates on the second lowest eigenvalue $\lu$ of $B$; this
eigenvalue is the lowest eigenvalue of $B$ restricted to odd functions.

\begin{lemma}\label{Lemma3.4}
The second lowest eigenvalue $\mu_2$ of $B$ satisfies
\begin{equation}
\mu_2 = \frac{m\pi^2}{4 L^2} + O (L^{-5/2}).\label{3.32}
\end{equation}
\end{lemma}

\begin{proof}
If $t_2$ is an eigenvalue of $B$, it is odd, thanks to the general
Sturm-Liouville theory. Define an operator $B_0$ in $L^2(\Er^+)$ by
\begin{equation*}
D(B_0)=H^1_0(\Er^+)\cap H^2(\Er^+), \quad B_0u=Bu\mid_{\Er^+};
\end{equation*}
then $\mu_2$ is the first eigenvalue of $B_0$. In particular, for all
$v\in H^1_0(\Er^+)\setminus\{0\}$, the following inequality holds:
\begin{equation*}
\mu_2\le \int_0^\infty \bigl(m\vabs{v'}^2 +
Wv^2\bigr)\left(\int_0^\infty v^2\, dx\right)^{-1}.
\end{equation*}
Define a test function $v$ by
\begin{equation*}
v(x)=\begin{cases}
\sin\bigl(\pi x/2\bigl(L-\sqrt{L}\bigr)\bigr)&\text{if $0\le x\le
L-\sqrt{L}$,}\\
r(x)/r\bigl(L-\sqrt{L}\bigr).&\text{if $x\ge L-\sqrt{L}$.}
\end{cases}
\end{equation*}
For all large enough $L$, $W\le 0$ on $\bigl[0, L-\sqrt{L}\bigr]$;
therefore an integration by parts on $\bigl[L-\sqrt{L},+\infty\bigr)$
shows that
\begin{equation*}
\int_0^\infty \bigl(m\vabs{v'}^2 + Wv^2\bigr)\, dx \le
\frac{\pi^2m}{8\bigl(L-\sqrt{L}\bigr)^3}
-\frac{R'\bigl(L-\sqrt{L}\bigr)}{R(L-\sqrt{L}}.
\end{equation*}
On the other hand,
\begin{equation*}
\int_0^{L-\sqrt{L}} v^2\, dx\ge \frac{L-\sqrt{L}}{2};
\end{equation*}
therefore, we see now that
\begin{equation}
\mu_2 \le \pi^2m /8L^2+O(L^{-5/2}.\label{3.33}
\end{equation}

Conversely, define $L_1>0$ by
\begin{equation*}
R(L_1)=\frac{1}{2}\left(1+\frac{\sqrt{1+3\nu +16\mu_2}}{2}\right);
\end{equation*}
then, thanks to~\eqref{2.11} and~\eqref{3.33}
\begin{equation*}
L_1=L +\frac 12 \ln 2 + O(1/L^2);
\end{equation*}
the function $x\mapsto \mu_2 -W(x)$ changes sign only at $x=L_1$;
therefore, there exists $L_2\le L_1$ such that $t_2$ increases on $[0,
L_2]$ and decreases on $[L_2,\infty)$. 

\textbf{Case 1: $L_2\le L-\sqrt L$.} 
Write $M=\sqrt{\mu_2/m}$; since $t_2$ and $t'_2$ cannot
vanish simultaneously and $t_2(0)$ vanishes, there exists a continuous
determination of the angle $\gamma$ such that
\begin{equation}
t_2=b\sin M\gamma, \quad t'_2 =bM\cos M\gamma, \quad \gamma(0)=0.\label{eq:7}
\end{equation}
A classical calculation gives
\begin{equation*}
1-\gamma'=\frac{W}{mM^2}\sin^2 M\gamma.
\end{equation*}
There exists a constant $C$ such that for all $L$ large enough,
\begin{equation*}
\forall x\in\bigl[0,L-\sqrt{L}\bigr],\quad  0\le -W(x)\le Ce^{-2\sqrt{L}};
\end{equation*}
if we let $\varepsilon^2= C me^{-2\sqrt L}$, $\gamma$ satisfies the differential inequality
by integration of the inequality
\begin{equation*}
\gamma'\le 1+\varepsilon \gamma(x)^2,
\end{equation*}
which we integrate as:
\begin{equation}
\gamma(x)\le \varepsilon^{-1}\tan\bigl(\varepsilon x).\label{eq:8}
\end{equation}

Thanks to~\eqref{eq:7}, $M\gamma(L_2)=\pi/2$; therefore thanks
to~\eqref{eq:8}, 
\begin{equation*}
M\ge \varepsilon\pi/2\tan(\varepsilon L_2);
\end{equation*}
but the definition of $\varepsilon$ and the condition on $L_2$ imply
immediately that
\begin{equation*}
M\ge \frac{\pi}{2(1+C\varepsilon^2)\bigl(L-\sqrt L\bigr)}
\end{equation*}
which implies immediately the desired estimate~\eqref{3.32}.

\textbf{Case 2: $L_2\ge L-\sqrt L$.} We have to give a different
argument: multiply the equation $-mt_2''+W t_2 =\mu_2 t_2$ by $v$ and
perform two integration by parts on the term involving a second
derivative; we eventually obtain
\begin{equation*}
\begin{split}
&\left(\frac{m\pi^2}{4\bigl(L-\sqrt L\bigr)^2}
-\mu_2\right)\int_0^{L-\sqrt L} t_2 v\, dx = -\int_0^{L-\sqrt L}
t_2Wv\, dx +\mu_2 \int_{L-\sqrt L}^{\infty} t_2 v\, dx \\
&\qquad\qquad + mt_2(L-\sqrt
L)\bigl[v'\bigl(L-\sqrt L -0\bigr) - v'\bigl(L-\sqrt L +0\bigr)\bigr].
\end{split}
\end{equation*}
Normalize $t_2$ so that $t_2(L_2)=1$; then, by concavity on $[0,L_2]$,
\begin{equation*}
t_2(x)\ge x/L_2\ge x/L_1;
\end{equation*}
therefore, there exists $C>0$ such that
\begin{equation*}
\int_0^{L-\sqrt L} t_2 v\, dx\ge CL;
\end{equation*}
on the other hand,
\begin{equation*}
\begin{split}
&\vabs{\int_0^{L-\sqrt L} t_2 Wv\, dx} =O\bigl(Le^{-2\sqrt
L}\bigr),\quad  \int_{L-\sqrt L}^{\infty} t_2 v\,
dx=O\bigl(\sqrt L\bigr),\\
&\vabs{v'\bigl(L-\sqrt L -0\bigr) - v'\bigl(L-\sqrt L +0\bigr)}
=O\bigl(e^{-2\sqrt L}\bigr) 
\end{split}
\end{equation*}
Therefore
\begin{equation*}
CL\vabs{\frac{m\pi^2}{4\bigl(L-\sqrt L\bigr)^2}
-\mu_2} =O\bigl(L^{-3/2}\bigr),
\end{equation*}
which implies immediately~\eqref{3.33}.

This proves Lemma~\ref{Lemma3.4}.
\end{proof}

\section{Existence of pulses}\label{sec:Existence-pulses}

We will prove that there are couples $(\xi,\eta)$ which satisfy~\eqref{3.4};
the main idea is 
to use an ansatz suggested by previous long computations; this ansatz
relies on the 
observation that for small values of $\varepsilon$, the pulse is simply
stretched; i.e., the 
new pulse obtained should be very close to a pulse with $\varepsilon=0$
but a smaller value 
of $\nu$. In fact, we will work in reverse: given $\nu$, we define $\nf= \nu
e^{-4y}$, with $y\ge 0$; we let $\Lf = L+y$, $\rf = r(\cdot,\nf)$, $\qf=q(\cdot,\nf)$,
$\thf=\theta(\cdot,\nf)$, $\smf=s(\cdot,\nf)$ and we have to determine an
ansatz $\ef$ for $\varepsilon$. We will 
start with this ansatz, with very little motivation.  The
justification will be provided 
by the proofs which come later.

Once the ansatz is determined, we get an
asymptotic
description of the behavior of $\ef$. 
Then, we go on to apply a version of the
implicit function theorem {\sl cum} estimates which is given in the appendix; for this
purpose, we have to check that $G(\rf,\qf,\thf,\ef,\nu)$ is very small, that 
$D_{(\xi,\eta,\tau,\varepsilon)}G(\rf,\qf,\thf,\ef,\nu)$ is invertible,
with an estimate on the
norm of its inverse; this necessitates an adequate choice of
functional spaces; finally, 
we check that the second derivatives of $G$ are not too large in a neighborhood of
$(\rf,\qf,\thf,\ef)$, and we are able to apply this implicit function theorem.  

Moreover, we obtain precise asymptotics on the second term of the expansion. These
estimates are crucial for the proof of the skew stabilization.

Define 
\begin{equation}
g_1(\varepsilon)  = G_1(\rf,\qf,\thf,\varepsilon,\nu).\label{4.1}
\end{equation}
Another expression of $g_1$  obtained by substituting the value of
$(-r'' + r)^\flat$ from~\eqref{2.3} is 
\begin{equation}
g_1(\varepsilon) = \frac{\nu -\nf}{1-\nf}(r^5- r^3)^\flat +
\varepsilon(r^2 q + 2r^3 q^2 - rq^2 -\theta
q)^\flat + \varepsilon^2(rq^4 + q^3)^\flat.\label{4.2}
\end{equation}
The ansatz for $\varepsilon$ is given by
\begin{equation}
(g_1(\ef),\smf)=0.\label{4.3}
\end{equation}

Our purpose now is to show that this ansatz determines uniquely $\ef$
for small enough 
values of $\nu$, and for all $y\ge 0$. Let us introduce new notations

\begin{align}
\kappa &= \frac{\nu -\nf}{1-\nf},\label{4.4}
\\
a_0&=(r^5 -r^3,s),\notag
\\
a_1&=( r^2 q + 2r^3 q^2 - rq^2 -\theta  q,s),\label{4.6}
\\
a_2&= (rq^4 + q^3,s).\notag
\end{align}
Of course, $a_j^\flat$ will be equal to $a_j$ with $r$, $q$, $\theta$ 
and $s$ replaced by their
analogues $\rf$, $\qf$, $\thf$ and $\smf$.         

Therefore,~\eqref{4.3} is equivalent to the following equation of the second degree:
\begin{equation}
\kappa a_0^\flat + \ef a_1^\flat +(\ef)^2a_2^\flat =0.\label{4.8}
\end{equation}
In order to have an asymptotic for $\ef$ we will need an asymptotic
for $a_j$, $0\le j\le 2$. 

We introduce the notation 
\begin{equation*}
\co
\end{equation*}
 for any quantity whose absolute value or relevant norm
is bounded by some power of $L$.

\begin{lemma}\label{Lemma4.1}
The following asymptotics hold
\begin{align}
a_0 &=-\frac{9}{64} + O(\sqrt\nu),\label{4.9}
\\
a_1&=\frac{9}{16}+ O\left(\frac{1}{L}\right),\label{4.10}
\\
a_2&=\frac{L^4}{36} +O(L^3).\label{4.11}
\end{align}
\end{lemma}

\begin{proof}
We have
\begin{equation*}
\begin{split}
(r^5 - r^3, s) &= (r^5 - r^3,\sigma) + O(\nu L) =- 2\int_0^\infty
(r^5 - r^3)r'\,dx + 
O(\sqrt\nu) \\
&=-\frac{9}{64} +O(\sqrt\nu).
\end{split}
\end{equation*}
This gives~\eqref{4.9}.

The second object is more difficult.  If we define
\begin{equation}
a=(r^2q+2r^3q^2-rq^2-\theta q,\sigma)\label{eq:32}
\end{equation}
estimate~\eqref{2.42} implies that $\vabs{a_1-a}=\nu\co$.
But, thanks to~\eqref{3.8} and~\eqref{3.29}, $\tilde a$ can be
rewritten as
\begin{equation*}
a=\bigl((r^2-\theta)q,\sigma) +\frac{1}{2}(\hat \rho,q)+
\frac{1}{2}(3r^2\sigma -\theta\sigma-\theta_1 r,q) -\frac12(q_1, r^3-Br).
\end{equation*}
Thus, we have
\begin{equation*}
a=-(r^4-\theta r^2, \phi_1) +\frac12 \frac{\partial}{\partial L}(r^4
-\theta r^2, \phi) +\frac{1}2(\hat \rho, q).
\end{equation*}
We transform this expression using~\eqref{3.19}:
\begin{equation}
a = \bigl(m (r^2 \phi')', \phi_1\bigr) -\frac{1}{2} \frac{\partial}{\partial L}
\bigl(m (r^2 \phi')', \phi\bigr) + +\frac{1}2(\hat \rho, q),
\end{equation}
and by integration by parts
\begin{equation}
a =m (r\sigma \phi', \phi') + \frac{1}{2}\frac{\partial m}{\partial
L}\bigl(r\phi', r\phi'\bigr) + \frac{1}2\bigl(\hat \rho,
q\bigr).\label{4.12}
\end{equation}
But 
\begin{align*}
\int r\sigma \vabs{\phi'}^2\,dx
&=-2\int_0^\infty r r' \vabs{\phi'}^2\,dx +O(\sqrt\nu)\\
&=-R\vabs{\phi'}^2\bigm|_0^\infty +2\int_0^\infty
R\phi'\phi''\,dx+O(\sqrt\nu)\\
& =\frac{3}{4}\vabs{\phi'(L)}^2
+2\int_0^L\left(R-\frac{3}{4}\right)\phi'\phi''\,dx +2\int_L^\infty
R\phi'\phi''\,dx.
\end{align*}
From estimate~\eqref{3.31} over $\phi''$, and the following facts
\begin{equation*}
\int_0^L|R-3/4|\,dx =O(1), \quad\int_L^\infty R\,dx = O(1),
\end{equation*}
we can see that
\begin{equation}
\int r\sigma \vabs{\phi'}^2\,dx= 3\vabs{\phi'(L)}^2/4+O(1/L),\label{eq:44}
\end{equation}
and from asymptotic~\eqref{3.26}, we obtain
\begin{equation*}
\int r\sigma \vabs{\phi'}^2\,dx =3 +O(1/L).
\end{equation*}
This implies~\eqref{4.10}.

The last of the three estimates concerns
\begin{equation*}
\int (rq^4 +q^3)s\,dx.
\end{equation*}
Thanks to~\eqref{3.23} and~\eqref{2.40}, the second term will
contribute an $O(L^3)$. Let us  
estimate  $\int r q^4s\,dx$, which will be equivalent to
$CL^4$. More precisely,
\begin{align*}
\int r q^4s\,dx &=\int rq^4\sigma\,dx +\nu \co\\
&=-2\int_0^\infty r r' r^4\phi^4\,dx +\nu \co\\
&=-\frac{R^3}{3}\phi^4\bigm|_0^\infty +\int_0^\infty \frac{R^3}{
3}(\phi^4)'\,dx +\nu \co\\
&=\frac{R^3(0)}{3}\phi^4(L) +\int_0^L\frac{R^3-R^3(0)}{3}(\phi^4)'\,dx
+\int_ L^\infty \frac{R^3}{3}(\phi^4)'\,dx + \nu \co.
\end{align*}
From the asymptotics~\eqref{3.21} on $\phi'$ and~\eqref{3.22} and~\eqref{3.22a}
on $\phi$, we can see that
\begin{align*}
&\int_0^L \vabs{\bigl(R^3 -R^3(0)\bigr)(\phi^4)'\bigr)}\,
dx=O(L^3)\int_0^L\vabs{R^3-R^3(0)} \, dx=O(L^3)\\
\intertext{and}
&\int_L^\infty R^3\vabs{(\phi^4)'}\, dx \le C\int_L^\infty x^3
e^{-2(x-L)}\, dx =O(L^3).
\end{align*}
We obtain now
\begin{equation*}
\int r q^4s\,dx =
\frac{1}{3}\left(\frac{3}{4}\right)^3\left(-\frac{2L}{3}\right)^4 +
O(L^3)  
=\frac{L^4}{36} + O(L^3),
\end{equation*}
which is precisely the last of the three estimate of our list.
\end{proof}

Therefore, for all $\nu$ small enough, and for all $y\ge 0$, the
discriminant of equation~\eqref{4.8} is strictly positive, and there exists a unique
positive solution of this equation, which is given by
\begin{equation}
\ef =\frac{-2a_0^\flat \kappa}{a_1^\flat +\sqrt{(a_1^\flat)^2-4\kappa
a_0^\flat a_2^\flat}}.\label{4.13}
\end{equation}
From this relation we get
\begin{equation}
\ef=-\frac{a_0^\flat \kappa}{a_1^\flat} +\kappa^2 \co\text{ and }
\frac{\ef}{\kappa}=\frac{1}{4}+O\left(\frac1L\right).\label{4.14}
\end{equation}

We will apply now lemma~\ref{Lemma7.1}.
Define
\begin{equation}
U=(\xi,\eta,\tau,\varepsilon),\quad \Uf=(\rf,\qf,\thf,\ef).\label{4.15}
\end{equation}
 Thus, we have to check
that $G(\Uf,\nu)$ is small, and to obtain bounds for the
inverse of $D_UG(\Uf,\nu)$ and for
the second derivative of $G$ in a neighborhood of $\Uf$.

The first result is a bound on $G(\Uf,\nu)$.

\begin{lemma}\label{Lemma4.2}
Assume that $y\in [0, L^p]$, where $p\ge 1$ is an
integer. The following bound holds:
\begin{equation*}
|G(\Uf,\nu)|=\kappa\co.
\end{equation*}
\end{lemma}

\begin{proof}
We have already obtained an expression for 
$G_1(\Uf,\nu)$: it is given by
\begin{equation*}
G_1(\Uf,\nu)= \kappa(r^5-r^3)^\flat +\ef(r^2q+ 2r^3q^2
-rq^2 -\theta q)^\flat +(\ef)^2(rq^4 +q^3)^\flat.
\end{equation*}
We can see that
\begin{equation*}
|(r^5 - r^3)^\flat| + |(r^2q+ 2r^3q^2 -rq^2 -\theta q)^\flat| +
|(rq^4 +q^3)^\flat| = \co.
\end{equation*}
With the help of~\eqref{4.14}, we obtain 
\begin{equation}
|G_1(\Uf,\nu)| = \kappa
\co.\label{4.16}
\end{equation}

An expression for $G_2(\Uf,\nu)$ is given by
\begin{equation*}
G_2(\rf,\qf,\thf,\ef,\nu) =\kappa(r^4q -r^3 -r^2q +r\theta)^\flat
+\ef(2r^2q^3 - rq^2 - q^3)^\flat +(\varepsilon^2 q^5)^\flat.
\end{equation*}
Arguing as above,
\begin{equation}
|G_2(\Uf,\nu)| =\kappa
\co.\label{4.17}
\end{equation}

Relations~\eqref{4.16} and~\eqref{4.17} imply the assertion of the
 lemma.
 \end{proof}

Let us compute now the differential of $G$ with respect to $U$, at $(\Uf,\nu)$:
\begin{equation*}
\begin{split}
D_\xi G_1(\Uf,\nu) &=-m\frac{d^2}{dx^2} +m
+(5r^4 - 3r^2)^\flat +\bigl[\varepsilon(2rq - q^2 +6r^2q^2) +\varepsilon^2
q^4\bigr]^\flat\\
&=\frac{1-\nu}{1-\nu^\flat}A^\flat +\kappa(5r^4-3r^2)^\flat +
\bigl[\varepsilon(2rq - q^2 +6r^2q^2) +\varepsilon^2 q^4\bigr]^\flat,\\
D_\eta G_1(\Uf,\nu) &=\bigl[\varepsilon(r^2 - 2rq -\theta
+4r^3q)+\varepsilon^2(3q^2 + 4rq^3)\bigr]^\flat,\\ 
D_\tau G_1(\Uf,\nu) &= -\bigl[\varepsilon q\bigr]\f,\\
D_\varepsilon G_1(\Uf,\nu) &=\bigl[-\theta q +r^2q - rq^2 +
2r^3q^2 +\varepsilon(2q^3 + 2rq^4)\bigr]^\flat,
\end{split}
\end{equation*}
and, similarly,
\begin{equation*}
\begin{split}
D_\xi G_2(\Uf,\nu) &=\bigl[\theta - 3r^2 - 2rq + 4r^3q
+\varepsilon(-q^2 + 4rq^3)\bigr]^\flat,\\
D_\eta G_2(\Uf,\nu) &=\frac{1-\nu}{1-\nu^\flat}B\f
+\kappa(r^4 - r^2)^\flat +\bigl[\varepsilon(6r^2q^2 - 3q^2 - 2rq)
+5\varepsilon^2q^4\bigr]^\flat,\\
D_\tau G_2(\Uf,\nu) & = \rf,\\
D_\varepsilon G_2(\Uf,\nu) &=\bigl[-rq^2
-q^3 +2r^2q^3+2\varepsilon q^5\bigr]^\flat.
\end{split}
\end{equation*}

Therefore, we may write             
\begin{equation*}
D_U G(\Uf,\nu) = (1-\kappa) \cgf + \kappa \calh,
\end{equation*}
where
\begin{equation*}
\calg= \begin{pmatrix}
A &0&0&z\\
w&B & r &\hat z
\end{pmatrix}
\end{equation*}
and the letters $w$, $z$ and $\hat z$ denote multiplication operators by 
the functions
\begin{equation*}
\begin{split}
w&=\theta - 3r^2 - 2rq +4r^3q,\\
z&=-\theta q +r^2q - rq^2 + 2r^3 q^2\\
\hat z&= -rq^2 - q^3 + 2r^2q^3.
\end{split}
\end{equation*}
The matrix $\calh$ has components $H_{ij}$, $1\le i \le 2, 1\le j, \le 4$ which satisfy
\begin{equation*}
\max_{i,j} |H_{ij}|_\infty = \co.
\end{equation*}

\begin{lemma}\label{Lemma4.3}
Let $E_2$ be the space
\begin{equation*}
E_2 = \{\xi\in H^2_{\textrm{even}}(\Er), \xi\perp s\}\times\{\eta\in
H^2_{\textrm{even}}(\Er), \eta\perp r\}\times \Er^2
\end{equation*}
and let $E_0$ be the space
\begin{equation*}
E_0 =(L^2_{\textrm{even}}(\Er))^2.
\end{equation*}
The space $E_2^\flat$ is defined in an obvious way.
Then, for all $p<\infty$, for all small enough $\nu$ and for $y\le
L^p$, $D_U G(\Uf, \nu)$ is an 
isomorphism from $E^\flat_2$ to $E_0$, and the norm of its inverse is
bounded by $\co$.
\end{lemma}

\begin{proof}
We first study $\calg$ as an operator from $E_2$ to $E_0$.
Given $h=(h_1,h_2)^T\in E_0$, we first solve
\begin{equation}
\calg U_1 =h,\label{4.18}
\end{equation}
where 
$U_1= (\xi_1, \eta_1, \tau_1, \varepsilon_1)^T$ is to be found in $E_2$. 
Equation~\eqref{4.17} is equivalent to the system
\begin{equation}
\left\{\begin{split}
&A \xi_1 + z\varepsilon_1 = h_1, \\
&w\xi_1 + B \eta_1 + r\tau_1 + \hat z \varepsilon_1 =
h_2.
\end{split}\right.\label{4.19}
\end{equation}
Since $\xi_1$ is perpendicular to $s$, the first equation gives
\begin{equation}
(z, s) \varepsilon_1 = (h_1, s).\label{4.20}
\end{equation}
But $(z, s)$ is precisely equal to $a_1$ defined at~\eqref{4.6}. Therefore,
there exists a unique $\varepsilon_1$ satisfying~\eqref{4.20} and
\begin{equation*}
|\varepsilon_1| \le \co |h|.
\end{equation*}
Restricted to even functions perpendicular
to $s$, $A$ is invertible, and its inverse is bounded thanks to
Theorem~\ref{Theorem2.7}. Therefore
\begin{equation*}
|\xi_1| =\co,
\end{equation*}
and by classical elliptic estimates
\begin{equation*}
|\xi_1|_{H^2(\Er)^2} =\co.
\end{equation*}

For the second equation of~\eqref{4.18} to have a solution, the following orthogonality
condition must be satisfied:
\begin{equation*}
h_2 - \hat z \varepsilon_1 - w\xi_1 - r\tau_1 \perp r.
\end{equation*}
Hence
\begin{equation*}
\tau_1 = \frac{(h_2 - \hat z \varepsilon_1 - w \xi_1, r)}{
|r|^2},
\end{equation*}
which implies immediately that
\begin{equation*}
|\tau_1| \le \|h\|_{E_0} \co.
\end{equation*}
The operator $B$ restricted to $r^\perp$ is invertible, and the norm of its inverse
is estimated by $\mu_2^{-1} = O(L^2)$, thanks to
Lemma~\ref{Lemma3.4}. Therefore, the unique 
$\eta_1$ orthogonal to $r$ which solves the second equation of~\eqref{4.18} satisfies
\begin{equation*}
|\eta_1|\le \|h\|_{E_0} \co.
\end{equation*}
By classical elliptic estimates
\begin{equation*}
|\eta_1|_{H^2(\Er)}\le \|h\|_{E_0} \co,
\end{equation*}
and therefore
\begin{equation*}
\|U_1\|_{E_2} \le \|h\|_{E_0} \co.
\end{equation*}
The full system
\begin{equation*}
(1-\kappa) \cgf U_1 + \kappa \calh U_1= h
\end{equation*}
can be rewritten
\begin{equation*}
U_1 + \frac{\kappa}{1- \kappa} (\cgf)^{-1} \calh U_1 = (\cgf)^{-1}
\frac{h}{1- \kappa}.
\end{equation*}
Clearly,
\begin{equation*}
\|(\cgf)^{-1} \calh\|_{\call(E_2^\flat)} =\co.
\end{equation*}
Therefore, for $\nu$ small enough and $y \le L^p$, we can see that
\begin{equation*}
\left\|\frac{\kappa}{1- \kappa}(\cgf)^{-1} \calh\right\|_{\call}(E_2^\flat)\le \frac{1}{2},
\end{equation*}
and our assertion is proved.
\end{proof}

We have to estimate $D^2G(U,\varepsilon,\nu)$ in a neighborhood of
$(\Uf,\nu)^T$. Let us choose for this neighborhood the ball of radius $1$
around  $\Uf$, in the space $E_2^\flat$. Without any calculation,
we can see that $D^2G$ contains only multiplications by functions which are polynomials in
$\xi$, $\eta$, $\tau$, $\varepsilon$ and $\nu$. Therefore, the norm of this operator is $\co$
for $y\le L^p$. 

Now, these results enable us to apply Lemma~\ref{Lemma7.1}. There
exists for all small enough $\nu$ a 
unique $U\in E_2$   such that
\begin{equation*}
G(U,\nu)=0.
\end{equation*}
Moreover,
\begin{equation}
\left\{
\begin{split}
&\|U-\Uf\|_{E^\flat_2} = \kappa\co,\\
&\|U - \Uf +D_U G(\Uf, \nu)^{-1}G(\Uf,\nu)\|_{E^\flat_2} =\kappa^2\co.
\end{split}
\right.
\label{4.21}
\end{equation}

 It will be useful in what follows to have an
asymptotic for 
\begin{equation*}
U_1 = D_U G(\Uf, \nu)^{-1}G(\Uf,\nu).
\end{equation*}
 The information needed is
summarized in Lemma~\ref{Lemma4.4}. But 
we need to introduce new notations. The orthogonal projections $P$ on $s$, $P_\perp$ on
$s^\perp$, $Q$ on $r$ and $Q_\perp$ on $r^\perp$ are defined by
\begin{equation*}
\left\{
\begin{split}
&P v = \frac{(s, v)s}{|s|^2}, \quad Q v = \frac{(r,v)r}{
|r|^2},\\
&P_\perp v = v - Pv, \quad Q_\perp v = v - Qv
\end{split}
\right.
\label{4.22}
\end{equation*}
We will also abuse notations: when $f$ is orthogonal to $r'$ we denote
by $A^{-1}f$ the unique solution $u$ of $Au=f$ which is orthogonal to
$r'$; the analogous convention will be used for $B$.

\begin{lemma}\label{Lemma4.4}
Let
\begin{equation*}
U_1= D_U G(\Uf,
\nu)^{-1}G(\Uf,\nu)=(\xi_1,\eta_1,\tau_1,\varepsilon_1)^T.
\end{equation*}
The following asymptotics hold  
\begin{align}
\xi_1 &=(A^\flat)^{-1} G_1(\Uf,\nu)+\kappa^2\co\label{4.23}\\
\eta_1&=(B^\flat )^{-1}Q_\perp^\flat \bigl(G_2(\Uf,\nu) -w^\flat \xi_1\bigr)
+\kappa^2\co,\label{4.24} \\
\tau_1 &=\frac{\bigr(G_2(\Uf,\nu) - \xi_1 w^\flat,
r^\flat\bigr)}{|r^\flat|^2}+\kappa^2\co,\label{4.25}\\
 \varepsilon_1 &=\kappa^2 \co.\label{4.26}
 \end{align}
 \end{lemma}

\begin{proof}
The vector $U_1$ satisfies the equation 
\begin{equation}
\cgf U_1 + \frac{\kappa}{1-\kappa}{\calh} U_1 = \frac{G(\Uf,\nu)}{
1-\kappa}.\label{4.27}
\end{equation}
The first equation of~\eqref{4.27} yields
\begin{equation*}
(z^\flat,\smf)\varepsilon_1 +\frac{\kappa}{1-\kappa}\bigl((\calh U_1)_1,\smf\bigr)=0.
\end{equation*}
Since $\|U_1\|_{E^\flat_2} = \kappa\co,$
we can see that~\eqref{4.26} holds. Therefore, 
\begin{equation*}
A^\flat\xi_1 = G_1(U^\flat, \nu) + \kappa^2 \co,
\end{equation*}
which together with the orthogonality $\xi_1 \perp s^\flat$ implies~\eqref{4.23}.
In a similar fashion,
\begin{equation*}
\xi_1 w^\flat +B^\flat \eta_1 +\rf \tau_1 +\hat z^\flat \varepsilon_1 =\frac{G_2(\Uf,\nu)}{
1-\kappa}-\frac{\kappa}{1-\kappa}(\calh U_1)_2,
\end{equation*}
so that
\begin{equation*}
|r^b|^2\tau_1 = \left(\frac{G_2(U^\flat,\nu)}{1-\kappa} - w^\flat \xi_1 ,
r^\flat\right) + \kappa^2 \co.
\end{equation*}
Hence~\eqref{4.25} holds, and~\eqref{4.24} is inferred immediately
from~\eqref{4.25}.
\end{proof}

The only remaining task is to check that the $(\xi,\eta)$ obtained in this fashion yields
an $u=\xi+i\sqrt\varepsilon\,\eta$ which has the adequate asymptotic behavior at $x=\infty$.
The linearized equation at infinity for $u$ is
\begin{equation*}
-\frac{3}{16}(1-\nu)u'' +\frac{3}{16}(1-\nu) u + i\sqrt\varepsilon \tau
u=0.
\end{equation*}
Therefore, if $\beta$ is the square root of
$\bigl(1+i\sqrt\varepsilon\,\tau\bigr)/m$ which has 
positive real part, there is a constant $\gamma$ such that
\begin{equation*}
u(x)\sim \gamma e^{-i\beta|x|}.
\end{equation*}

The phase of $u$ is given by
\begin{equation*}
\arg u(x) =\arctan \frac{\sqrt\varepsilon \eta(x)}{\xi(x)}.
\end{equation*}
According to our ansatz, this is very close to
\begin{equation*}
\arctan\bigl(\sqrt\varepsilon\,\phi^\flat(x)\bigr),
\end{equation*}
whose behavior is exactly the behavior we postulated at the beginning of this article.

Finally, we obtain:

\begin{theorem}\label{Theorem4.5}
For all $p>0$, for all $\nu$ small enough, 
there exists $C>0$ such that for all $\alpha$ satisfying
\begin{equation}
|\alpha|\le \alpha_m=\frac{1}{2}\sqrt\nu\,\bigl(1
-CL^{-p}\bigr),\label{eq:39}
\end{equation}
there exists a pulse solution of~\eqref{1.5}, i.e. a solution of the form
\begin{equation*}
u(x,t)=e^{i\omega t} r(x) e^{i\phi(x)},
\end{equation*}
where $r$ is a positive function which decays exponentially at
infinity, and $\phi(x)$ is asymptotic at infinity to $-C|x|+D$,
with $C$ a positive constant and $D$ a real constant.
\end{theorem}

\begin{proof}
We just have to translate the condition on $\alpha$ from the condition
on $\varepsilon\f$. We know from~\eqref{4.14} that there exist $C'$,
$k$ and $K$
such that
\begin{equation*}
\varepsilon\f \ge\left(\frac14 -\frac{C'}{L+y}\right)\kappa -\nu^2 KL^k
\end{equation*}
Since the function
\begin{equation*}
y\mapsto \left(\frac14 -\frac{C'}{L+y}\right)\frac{1-e^{-4y}}{1-4\nu e^{-4y}}
\end{equation*}
is increasing on $\Er^+$ for $\nu$ small enough, we can see that if
\begin{equation}
\alpha^2\le\left(\frac14 -\frac{C'}{L+L^p}\right)\nu
\frac{1-e^{-4L^p}}{1-4\nu e^{-4L^p}} -\nu^2 KL^k,\label{eq:9}
\end{equation}
the existence theory works; but it is clear that for $\nu$ small
enough, the right hand side of~\eqref{eq:9} is at least equal to
$\nu(1-C/L^p)^2/4$, and the theorem is proved.
\end{proof}

\section{Stability of the pulse}\label{sec:Stability-pulse}

Our purpose is to study now the spectrum of
\begin{equation*}
D_uF(\xi+i\sqrt\varepsilon\,\eta,\sqrt\varepsilon\,\tau,\sqrt\varepsilon,\nu),
\end{equation*}
under the assumption
\begin{equation*}
\text{the parameter $\nu$ is small enough and there exists $p\ge 1$
such that $0\le y \le L^p$ .}\label{5.2} 
\end{equation*}
We define
\begin{equation}\label{5.1}
\cald =D_uF(\xi+i\sqrt\varepsilon\,\eta,\sqrt\varepsilon\,\tau,
\sqrt\varepsilon,\nu)/ (1-\kappa),
\end{equation}
and we observe that
\begin{equation}
\cald=\caf + \sqrt\kappa \calb + \kappa \calc,\label{5.7}
\end{equation}
where
\begin{equation}
\caf = 
\begin{pmatrix}
\Af&0\\ 0&\Bf
\end{pmatrix}, \quad \calb=\begin{pmatrix}
0&\calb_{12}\\\calb_{21}&0
\end{pmatrix}\text{ and }\calc=\begin{pmatrix}
\calc_{11}&0\\0&\calc_{22}
\end{pmatrix}.\label{eq:49}
\end{equation}
More precisely, we have
\begin{align}
\calb_{12}&=\frac{\sqrt\varepsilon(-\tau +\xi^2 - 2\xi\eta +
4\xi^3\eta) + \varepsilon^{3/2}(3\eta^2 + 4\xi
\eta^3)}{\sqrt\kappa\,(1-\kappa)}\label{eq:10}\\
\calb_{21}&=\frac{\sqrt\varepsilon(\tau - 3\xi^2 - 2\xi\eta +
4\xi^3\eta) +
\varepsilon^{3/2}(-\eta^2+4\xi\eta^3)}{\sqrt\kappa\,(1-\kappa)}\label{eq:11}\\ 
\calc_{11}&=\left[\frac{5\xi^4 -3\xi^2 +\varepsilon (2\xi\eta -\eta^2
+6\xi^2\eta^2)+\varepsilon^2\eta^4}{1-\kappa}-(5r^4
-3r^2)\f\right]\frac1\kappa\label{eq:12}\\ 
\calc_{22}&=\left[\frac{\xi^4-\xi^2 +\varepsilon(-2\xi\eta -3\eta^2
+6\xi^2\eta^2)+5\varepsilon^2 \eta^4}{1-\kappa}-(r^4
-r^2)\f\right]\frac1\kappa.\label{eq:13} 
\end{align}

It is immediate
from~\eqref{4.21},~\eqref{4.23},~\eqref{4.24},~\eqref{4.25}
and~\eqref{4.26} that the following estimate holds:
\begin{equation*}
\norm[L^\infty(\Er)]{\calb_{12}} +\norm[L^\infty(\Er)]{\calb_{21}}
+\norm[L^\infty(\Er)]{\calc_{11}} + \norm[L^\infty(\Er)]{\calc_{22}}=\kappa\co.
\end{equation*}

Identifying the multiplication by $\calb$ or $\calc$ to an
operator from $L^2(\Er)^2$ to itself, we see that
\begin{equation}
\norm[\call]{\calb} +\norm[\call]{\calc}=\co.\label{eq:14}
\end{equation}

The first result concerns the dimension of the generalized eigenspace
associated to the eigenvalues of $\cald$ contained inside a conveniently small
circle; for this purpose, we let $K_1>0$ be a constant such that for all
$L$ large enough
\begin{equation*}
\mu_2+\lambda\ge 4/\bigl(K_1 L^2)
\end{equation*}
Theorem~\ref{Theorem2.8} and Lemma~\ref{Lemma3.4} guarantee the
existence of such a $K_1$.

\begin{lemma}\label{Lemma5.1}
Assume~\eqref{5.2}; let
$\gamma$ be the circle of radius $2/\bigl(K_1L^{2p}\bigr)$ and of
center $0$, traveled 
once in the positive direction. Then, the spectrum of
$\cald$ does not intersect $\gamma$; the eigenprojection associated to the
part of the spectrum contained inside $\gamma$ is of rank 3.
\end{lemma}

\begin{proof}
From Lemma \ref{Lemma2.6} and Theorem \ref{Theorem2.8}, we know that
the spectrum of $\caf$ inside $\gamma$ contains only the semisimple
double eigenvalue $0$ and the simple eigenvalue $\lf$. We know also
that $\mu_2\f$ is the closest element of the remainder of the spectrum
of $\caf$ to $\gamma$.
Define
\begin{equation*}
d(\zeta)=\min\bigl(\bigl|\zeta\bigr|,\bigl|\zeta-\lf\bigr|,
\bigl|\zeta-\mu_2\f\bigr|\bigr).
\end{equation*}
We define a twice punctured disk $D$ by
\begin{equation*}
D\f=\{\zeta\in\Ce: \abs{\zeta}<\mu_2\f, \zeta\neq 0, \zeta\neq \lf\}.
\end{equation*}
Then we have the following estimate on the resolvent
$\calr^0(\zeta)=(\zeta-\caf)^{-1}$ of $\caf$:
\begin{equation*}
\forall \zeta\in D\f,\quad \vnorm[\call]{\calr^0(\zeta)}\le 1/d(\zeta).
\end{equation*}

It will be convenient to embed $\cald$ into a holomorphic family of
operators $\cald(c)$ defined by
\begin{equation*}
\cald(c)=\caf +c \calb+c^2\calc,
\end{equation*}
and to do the same type of calculations as in Chapter II of \cite{Kat};
however, the results quoted therein cannot be directly applied, since
our eigenvalues are not uniformly isolated with respect to the
remainder of the spectrum.

The Neumann perturbation series for $\cald(c)$ is the expression
\begin{equation}
\calr(\zeta,c)=\sum_{j=0}^\infty\bigl((\zeta-\caf)^{-1}(c\calb +
c^2\calc)\bigr)^j (\zeta- 
\caf)^{-1};\label{5.21}
\end{equation}
this series converges provided that
\begin{equation}
\vabs{c}\le
\min\bigl(d(\zeta)/\bigl(\vnorm[\call]{\calb}+\vnorm[\call]{\calc}\bigr),1\bigr). 
\end{equation}
We infer from~\eqref{eq:14} that there exist $K_2>0$ and $k_2$ such
that
\begin{equation*}
\vnorm[\call]{\calb}+\vnorm[\call]{\calc}\le K_2 L^{k_2};
\end{equation*}
without loss of generality, we may assume that for all the $L's$ that
we consider, $K_2 L^{k_2}$ is larger than $1$;
thus the function $\calr$ is holomorphic with values in $\call$ on the set
\begin{equation}
\Delta=\bigl\{(\zeta,c): \zeta\in D\f, L^{k_2}K_2\vabs{c}<
d(\zeta)\bigr\},\label{eq:15} 
\end{equation}
and the following estimate holds over $\Delta$:
\begin{equation}
\vnorm[\call]{\calr(\zeta,c)}\le \frac{1}{d(\zeta)-\vabs{c}K_2L^{k_2}}.\label{eq:20}
\end{equation}
If $\gamma$ is the circle defined in the statement of the lemma, any
element of $\gamma$ satisfies for $L$ large enough
$d(\zeta)> 1/(K_1 L^{2p})$.
Therefore, the circle $\gamma$ is
included in the resolvent set of $\cald(c)$ provided that
\begin{equation*}
K_1 K_2 L^{2p+k_2}\vabs{c}<1.
\end{equation*}
This will hold for $c=\kappa$ and $L$ large enough, and the first
assertion of the lemma is proved.

Since $\calr$ is holomorphic in a neighborhood of $c=0$, it admits a
Taylor series of the form 
\begin{equation*}
\calr(\zeta,c)=\sum_{j=0}^\infty \calr^j(\zeta) c^j,
\end{equation*}
with
\begin{equation*}
\calr^j(\zeta)=\frac1{2\pi i}\int_{\gamma'}\calr(\zeta,c) c^{-j-1}\, dc,
\end{equation*}
where $\gamma'$ is a circle about $0$ of radius strictly less than
$d(\zeta)/\bigl(K_2 L^{k_2}\bigr)$. The first term of the expansion is equal to
\begin{equation}
\calr^0(\zeta)=\bigl(\zeta-\caf\bigr)^{-1},\label{eq:16}
\end{equation}
hence our notations are coherent. If we take $\gamma'$ of radius
$d(\zeta)/\bigl(2K_2 L^{k_2}\bigr)$, we have the estimates
\begin{equation}
\vnorm[\call]{\calr^j(\zeta)}\le \frac{2}{d(\zeta)}\left(\frac{2K_2
L^{k_2}}{d(\zeta)}\right)^j, \label{eq:21}
\end{equation}
and also, for
all $c$ such that $\abs{c}< 
d(\zeta)/\bigl(2K_2L^{k_2}\bigr)$, and for all $k\ge 0$,
\begin{equation}
\vnorm[\call]{\calr(\zeta,c)-\sum_{j=0}^{k-1}\calr^j(\zeta)c^j} \le 
\left[\frac{2K_2 L^{k_2}\abs{c}}{d(\zeta)}\right]^k
\frac{1}{d(\zeta)-2K_2 L^{k_2}\abs{c}}.\label{eq:18}
\end{equation}

The projection on the total eigenspace of $\cald(c)$ relative to the eigenvalues
inside $\gamma$ is given by
\begin{equation*}
\calp(c)=\frac{1}{2i\pi}\int_\gamma \calr(\zeta,c)\, d\zeta.
\end{equation*}
We infer from~\eqref{eq:21} that the operator $\calp^j$ defined by
\begin{equation*}
\calp^j=\frac{1}{2i\pi}\int_\gamma \calr^j(\zeta)\, d\zeta
\end{equation*}
satisfies the estimate
\begin{equation*}
\vnorm[\call]{\calp^j}\le 4\bigl(4K_1 K_2L^{2p+k_2}\bigr)^j;
\end{equation*}
for all $c$ such that
$2K_1K_2L^{k_2+2p}\abs{c}<1$ and for all $k\ge 0$, we have
\begin{equation}
\vnorm[\call]{\calp(c)-\sum_{j=0}^{k-1}\calp^j
c^j}\le
\frac
{
4\bigl(2K_1K_2L^{k_2+2p}\abs{c}\bigr)^k 
}
{1-2\sqrt\kappa K_1 K_2 L^{k_2+2p}}.\label{eq:19}
\end{equation}
For $2K_1K_2L^{k_2+2p}\abs{c}<1$, $\calp(c)$ is a
holomorphic function of $c$.
For $c=0$, relation~\eqref{eq:16}
implies that $\calp_0$ is the projection on the space spanned by $s$,
$ir$ and $r'$; therefore, it is of dimension $3$. A classical argument
shows that a continuous family of projections which is of finite rank
for some value of the parameter is of finite rank for all its
values. This proves the second assertion of the lemma.
\end{proof}

Let us define now a basis of $\im \calp(c)$ which depends analytically
on $c$ in a neighborhood of $0$; since we expect to find a
nontrivial Jordan block, we will not 
look for a basis of eigenvectors of $\cald(c)$.

\begin{lemma}\label{Lemma5.2}Define vectors
\begin{equation}
\begin{split}
&w(c)=\calp(c)\begin{pmatrix}
\smf \\0
\end{pmatrix} , \quad t(c)=\calp(c) 
\begin{pmatrix}-c\eta \sqrt{\varepsilon/\kappa}\\ \rf + c^2(\xi-\rf)/\kappa
\end{pmatrix}
,\\
&v(c)=\calp(c) \begin{pmatrix} (\rf)' +c^2(\xi' -(\rf)')/\kappa\\
c\eta'\sqrt{\varepsilon/\kappa} \end{pmatrix} .
\end{split}
\end{equation}
There exist positive numbers $K_3>2K_1K_2$, $k_3>k_2+2p$, $K_4$ and $k_4$ such that
for $K_3 L^{k_3}\abs{c}<1$ the vectors $w(c)$, $t(c)$
and $v(c)$ are holomorphic in $c$ and form a basis of $\im \calp(c)$; in this basis, the
restriction of the operator $\cald(c)$ to $\im \calp(c)$ has matrix
$M(c)$; $M(c)$ is analytic in $c$ and for $K_3 L^{k_3}\abs{c}<1$, it
admits the expansion
\begin{equation*}
M(c)=\sum_{j=0}^\infty M^j c^j,
\end{equation*}
with the estimate
\begin{equation}
\vnorm[\call(\Ce^3)]{M^j}\le K_4 L^{k_4} \bigl(2K_3 L^{k_3}\bigr)^j.\label{eq:23}
\end{equation}
Similarly, $w$, $t$ and $v$ admit expansions
\begin{equation*}
w(c)=\sum_{j=0}^\infty w^j c^j, \quad t(c)=\sum_{j=0}^\infty t^j c^j,
\quad v(c)=\sum_{j=0}^\infty v^j c^j
\end{equation*}
with the estimates
\begin{equation*}
\max\bigl(\vabs{w^j},\vabs{t^j},\vabs{v^j}\bigr)\le K_4 L^{k_4}
\bigl(2K_3 L^{k_3}\bigr)^j. 
\end{equation*}
Finally, when $c=\sqrt\kappa$, $M(c)$ has the form
\begin{equation*}
M\bigl(\sqrt\kappa\,\bigr)=\begin{pmatrix}
M_{11}(\sqrt\kappa\,\bigr)&0&0\\M_{21}(\sqrt\kappa\,\bigr)&0&0\\0&0&0
\end{pmatrix}.
\end{equation*}
\end{lemma}

\begin{proof} It is clear that $w$, $t$ and $v$ are holomorphic
functions of $c$ if $2K_1K_2L^{k_2+2p}\abs{c}<1$. By construction,
they belong to $\im\calp(c)$. For $c=0$, $w(0)$,
$t(0)$ and $v(0)$ are orthogonal vectors spanning $\calp(0)$;
therefore, there is a neighborhood of $c=0$ for which $w(c)$, $t(c)$
and $v(c)$ constitute a basis of $\im\calp(c)$. Let
$\langle\>,\>\rangle$ denote the $\Ce$-bilinear product on
$L^2(\Er;\Ce)$ defined as follows: if 
\begin{equation*}
X=\begin{pmatrix}
X_1\\X_2
\end{pmatrix}\text{ and }Y=\begin{pmatrix}
Y_1\\Y_2
\end{pmatrix}
\end{equation*}
belong to $L^2(\Er;\Ce)^2$, we let
\begin{equation*}
\langle X,Y\rangle=\int_\Er \bigl(X_1(x)Y_1(x)+X_2(x)Y_2(x)\bigr)\, dx.
\end{equation*}
In order to estimate the size of the neighborhood of $0$ where $w(c)$,
$t(c)$ and $v(c)$ constitute a basis of $\im\calp(c)$, we introduce a
quasi Gram matrix $\Gamma(c)$ defined by
\begin{equation*}
\Gamma(c)=\begin{pmatrix}
\langle w(c), w(c)\rangle&\langle w(c), t(c)\rangle&0\\\langle
t(c),w(c)\rangle&\langle t(c), t(c)\rangle &0\\0&0&\langle v(c), v(c)\rangle
\end{pmatrix}.
\end{equation*}
The terms $\langle w(c),v(c)\rangle$ and $\langle t(c),v(c)\rangle$
vanish because
$w$ and $t$ are even
functions, while $v$ is odd. We also remark that $\Gamma(c)$ is
symmetric, not hermitian, and that
\begin{equation*}
\Gamma(0)=\begin{pmatrix}
3/8 +O\bigl(\sqrt{\nu\f}\,\bigr)&0&0\\
0&3(L+y)/2+O(1)&0\\
0&0&3/8 +O\bigl(\sqrt{\nu\f}\,\bigr)
\end{pmatrix}.
\end{equation*}
We infer from the hypothesis $\vabs{c}\le
1$ and from~\eqref{4.21} that there exist $K_3>2K_1K_2$ and $k_3\ge 2k_2+p$ such that
\begin{equation}
\vnorm[\call(\Er^3)]{\Gamma(0)^{-1}\bigl(\Gamma(c)-\Gamma(0)}\le
\vabs{c} K_3 L^{k_3}.\label{eq:22}
\end{equation}
Therefore, if $\vabs{c} K_3 L^{k_3}<1$, $\Gamma(c)$ is invertible, its
inverse is analytical in $c$ and it satisfies the estimate
\begin{equation*}
\vnorm[\call(\Ce^3)]{\Gamma(c)^{-1}} \le
\frac{\vnorm[\call(\Ce^3)]{\Gamma(0)^{-1}}}{1-\abs{c}K_3 L^{k_3}}.
\end{equation*}

Denote by
\begin{equation*}
\begin{bmatrix}
w(c)&t(c)&v(c)
\end{bmatrix}
\end{equation*}
the function from $\Er$ to $3\times 3$ complex matrices whose columns
vectors are respectively equal to $w(c)$, $t(c)$ and $v(c)$.
We define three vectors $\hat w(c)$, $\hat t(c)$ and $\hat v(c)$ by
\begin{equation*}
\begin{bmatrix}
\hat w(c)&\hat t(c)&\hat v(c)
\end{bmatrix}=\begin{bmatrix}
\calp(c) w(c)&\calp(c)t(c)&\calp(c)v(c)
\end{bmatrix}\Gamma(c)^{-1}
\end{equation*}
Then, for $K_3 L^{k_3}\abs{c}<1$, $\hat w(c)$, $\hat t(c)$ and $\hat
v(c)$ are holomorphic in $c$ and they constitute a basis of
$\im\calp(c)^T$ which is dual to the basis $w(c)$, $t(c)$, $v(c)$ of
$\im\calp(c)$, i.e.
\begin{equation*}
\int_\Er \begin{bmatrix}
w(c)&t(c)&v(c)
\end{bmatrix}^T \begin{bmatrix}
\hat w(c)&\hat t(c)&\hat v(c)
\end{bmatrix}\, dx =\one_{\Ce^3}.
\end{equation*}
Thus, $M(c)$ is given by
\begin{equation*}
M(c)=\int_\Er \begin{bmatrix}
\hat w(c)&\hat t(c)&\hat v(c)
\end{bmatrix}^T \begin{bmatrix}
\cald(c) w(c)&\cald(c) t(c)&\cald(c)v(c)
\end{bmatrix}\, dx.
\end{equation*}
It is immediate that $\cald(c)w(c)$, $\cald(c)t(c)$ and $\cald(c)v(c)$
are analytical in $c$ with values in $L^2(\Er;\Ce)^2$: in fact, they
are polynomial in $c$. The study of section~\ref{sec:Existence-pulses}
and in particular estimates~\eqref{4.21}, completed by
estimate~\eqref{eq:19} for $k=0$ and estimate~\eqref{eq:22}  show that
there exist numbers $K_4$ 
and $k_4$ such that the following estimates hold:
\begin{gather*}
\text{if $\vabs{c}=1/\bigl(2K_3 L^{k_3}\bigr)$, then}\\
\max\bigl(\vabs{w(c)},\vabs{t(c)},\vabs{v(c)},\vabs{\hat
w(c)},\vabs{\hat t(c)},\vabs{\hat v(c)}\bigr)\le K_4 L^{k_4}\\
\text{ and
} \max_{1\le i,j \le 3} \vabs{M_{jk}(c)}\le K_4 L^{k_4}.
\end{gather*}
Standard methods of analytic functions give~\eqref{eq:23} and
analogous estimates for $y$, $t$ and $v$.
Finally, differentiating with respect to $\beta$ the relations
\begin{equation*}
F\bigl(e^{i\beta}\bigl(\xi+i\sqrt\varepsilon\,\eta\bigr),\sqrt\varepsilon\,\tau,
\sqrt\varepsilon,\nu\bigr)=0 
F\bigl(\bigl(\xi+i\sqrt\varepsilon\,\eta\bigr(\cdot +\beta),
\sqrt\varepsilon\,\tau,\sqrt\varepsilon,\nu\bigr)=0
\end{equation*}
we find that $t\bigl(\sqrt\kappa\,\bigr)$ and
$v\bigl(\sqrt\kappa\,\bigr)$ are eigenvectors of $\cald(\sqrt\kappa\,)$ relatively to 
the eigenvalue $0$. Thus $w$ and $t$ belong to
$\im\calp$. 
Therefore the form of the
matrix $M\bigl(\sqrt\kappa\,\bigr)$ is clear, and this concludes the proof of the lemma.
\end{proof}

All our effort will be concentrated now on getting the first three terms of the
expansion of $M_{11}(c)$ since we infer
from~\eqref{eq:23} that
\begin{equation*}
\vabs{M_{11}(\sqrt\kappa\,)-M_{11}^0 -\sqrt\kappa\,M_{11}^1 -\kappa
M_{11}^2}\le \kappa^{3/2}\co.
\end{equation*}

\begin{lemma}\label{Lemma5.3} The first terms of the expansion of
$\calp$ have the form
\begin{equation}
\calp^1 = 
\begin{pmatrix}
0& \calp^1_{12}\\ \calp^1_{21}& 0
\end{pmatrix}
, \quad \calp^2= 
\begin{pmatrix}
\calp_{11}^2&
0\\ 0& \calp^2_{22} 
\end{pmatrix}
,\label{5.28}
\end{equation}
and the first terms of the expansion of $w$ and $t$ are given by
\begin{equation}
w^1 = 
\begin{pmatrix}
0\\ s_1
\end{pmatrix}
,\quad w^2 = 
\begin{pmatrix}
s_2\\ 0
\end{pmatrix}
, \quad t^1 =
\begin{pmatrix}
r_1\\0
\end{pmatrix}
, \quad t^2=
\begin{pmatrix}
0\\ r_2
\end{pmatrix}
,\label{5.29}
\end{equation}
with
\begin{equation}
s_1 = (\lf -\Bf)^{-1} \Qpf\calb_{21}\smf \text{ and } r_1 =
-(\Af)^{-1} \ppf \calb_{12} \rf 
- \sqrt{\frac{\varepsilon}{\kappa}}\pf \eta.\label{5.30}
\end{equation}
\end{lemma}

\begin{proof}
The first statement of the lemma is an immediate consequence of the
definition of the $\calp^j$'s and of the special form of $\cala\f$,
$\calb$ and $\calc$. Let us calculate $\calp^1_{21}\smf$:
$\calr^1(\zeta)$ is equal to $\calr^0(\zeta)\calb\calr^0(\zeta)$ and thus
\begin{equation*}
\calr^1(\zeta)
\begin{pmatrix}
\smf\\0
\end{pmatrix}
=\begin{pmatrix}
0\\ (\lf-\zeta)^{-1}(\zeta-\Bf)^{-1} \calb_{21}
\smf
\end{pmatrix}
.
\end{equation*}
Hence, thanks to the theorem of residues,
\begin{equation*}
s_1= \calp_{21}^1\smf =\frac{1}{2 i \pi}\int_\gamma (\zeta - \Bf)^{-1} \calb_{21} \smf
\frac{d\zeta}{\lf - \zeta} = (\lf - \Bf)^{-1}\Qpf \calb_{21}\smf.
\end{equation*}
Similarly, according to the definition of $t(c)$,
\begin{equation*}
r_1 = \calp_{12} \rf -\sqrt{\frac{\varepsilon}{\kappa}} \pf \eta,
\end{equation*}
and
\begin{equation*}
\calp^1_{12}\rf = \frac{1}{i \pi}\int_\gamma (\zeta - \Af)^{-1}
\calb_{12} \rf\frac{d\zeta}{\zeta} = (- \Af)^{-1} \ppf \calb_{12}
\rf.
\end{equation*}
\end{proof}

Now we are able to give expressions for the first coefficients of the power expansion of
$M_{11}(c)$:

\begin{lemma}\label{Lemma5.4}
The coefficients $M^j_{11}$, for $0 \le j \le 2$ and $M^j_{21}$ for
$j=0,1$ are given by
\begin{align}
&M_{11}^0 = \lf, \quad M_{21}^0 = 0,\label{5.31}\\
&M_{11}^1 =0 , \quad M^1_{21} = \frac{(\calb_{21}\smf - \lf s_1,
\rf)}{|\rf|^2}, \label{5.32}\\ 
&M_{11}^2=\frac{(\calb_{12}s_1 + \calc_{11}\smf - M_{21}^1r_1,
\smf)}{|\smf|^2}.\label{5.33}
\end{align}
\end{lemma}

\begin{proof}
From the analyticity properties we infer that
\begin{align}
\caf y^0 &= M_{11}^0 y^0 + M_{21}^0 t^0,\label{5.34}\\
\caf w^1 + \calb y^0 &= \sum_{j=0}^1 \bigl(M_{11}^j y^{1-j} + M_{21}^j
t^{1-j}\bigr)\label{5.35}\\
\caf w^2 + \calb w^1 + \calc y^0 &= \sum_{j=0}^2 \bigl(M_{11}^j y^{2-j} + M_{21}^j
t^{2-j}\bigr)\label{5.36}
\end{align}
From~\eqref{5.34} and the value of $t(0)$ we infer
immediately~\eqref{5.31}. Relation~\eqref{5.32} can be rewritten with 
the help of~\eqref{5.29} as
\begin{equation*}
\begin{pmatrix}
0\\ \Bf r_1
\end{pmatrix}
+\begin{pmatrix}
0\\ \calb_{21}\smf
\end{pmatrix}
= \lf \begin{pmatrix}
0\\ s_1
\end{pmatrix}
+
M_{11}^1 \begin{pmatrix}
\smf\\ 0
\end{pmatrix}
+ M_{21}^1 \begin{pmatrix}
0\\\rf
\end{pmatrix}
.
\end{equation*}
It is immediate that $M_{11}^1$ vanishes. Moreover, we perform the scalar product of
the second component of the above identity with $\rf$;  since $(\Bf s_1, \rf) = (s_1, \Bf
\rf)=0$, we obtain the second relation of~\eqref{5.32}.

Finally,~\eqref{5.36} implies that
\begin{equation*}
\Af s_2 + \calb_{12} s_1 + \calc_{11} \smf = \lf s_2 + M^2_{11} \smf +
M^1 _{21} r_1.
\end{equation*}
If we perform the scalar product of this relation with $\smf$ and if we observe that $(\Af
s_2, \smf) = (s_2, \lf \smf)$, we obtain~\eqref{5.33}. 
\end{proof}

We give now an asymptotic of $M^2_{11}$ in terms of $\rf$, $\qf$ and $\sigf$:

\begin{theorem}\label{Theorem5.5}
The coefficient $M^2_{11}$ has the following asymptotics
\begin{equation}
M_{11}^2 = \frac{\varepsilon}{\kappa}
|\smf|^{-2}\left[\frac{\partial a}{\partial L}\right]\f +
\bigl(\sqrt{\nuf} +\sqrt{ \kappa}\bigr) 
\co,\label{5.37}
\end{equation}
where $a$ has been defined at~\eqref{eq:32}.
\end{theorem}

\begin{proof}
According to~\eqref{2.40} and Lemma~\ref{Lemma4.4}, we have
\begin{equation*}
\calb_{21} \smf = \sqrt{\frac{\varepsilon}{\kappa}}\bigl[(\theta - 3 r^2 - 2rq + 4 r^3
q)\sigma\bigr]^\flat + (\kappa + \nf )\co,
\end{equation*}
and in virtue of~\eqref{3.29} and~\eqref{3.30}, we can see that
\begin{equation}
\calb_{21}\smf = -\sqrt{\frac{\varepsilon}{\kappa}}(B q_1 + \theta_1
r)^\flat + (\kappa + \nf)\co.\label{eq:40}
\end{equation}
Therefore
\begin{equation}
s_1 = \sqrt{\frac{\varepsilon}{\kappa}}\Qpf \qf_1 + (\kappa +\nf)\co\label{eq:45}
\end{equation}
and
\begin{equation*}
M^1_{21} = - \sqrt{\frac{\varepsilon}{\kappa}} \thf_1 +(\kappa + \nf)\co.
\end{equation*}

Similarly, according to Lemma~\ref{Lemma4.4},
\begin{equation*}
\calb_{12}\rf = \sqrt{\frac{\varepsilon}{\kappa}}\bigl[(-\theta + r^2 - 2 rq + 4 r^3
q)r\bigr]^\flat + \kappa\co,
\end{equation*}
and thanks to~\eqref{3.8} we have
\begin{equation*}
\calb_{12} \rf= \sqrt{\frac{\varepsilon}{\kappa}}(Aq)^\flat+(\kappa+\nf)\co.
\end{equation*}
Therefore
\begin{equation*}
r_1 = -\sqrt{\frac{\varepsilon}{\kappa}} \qf + (\kappa+\nf)\co,
\end{equation*}
and in particular,
\begin{equation}
M_{21}^1 (r_1, \smf) =
\frac{\varepsilon}{\kappa}\thf_1(q,\sigma)^\flat + 
(\kappa + \nuf)\co.\label{5.38}
\end{equation}
There remains to calculate
\begin{equation*}
(\calb_{12} s_1 + \calc_{11} \smf, \smf).
\end{equation*}
We infer from~\eqref{eq:45} that
\begin{equation*}
(\calb_{12} s_1, \smf) = \left(\calb_{12}\sqrt{\frac{\varepsilon}{\kappa}} \qf_1,\smf\right)
-\left(\calb_{12}\sqrt{\frac{\varepsilon}{\kappa}}\Qf\qf_1,\smf\right)
+(\kappa + \nuf)\co 
\end{equation*}
and we have the following relation
\begin{equation*}
(\calb_{12}\rf,\smf)=\sqrt{\frac{\varepsilon}{\kappa}}\bigl((Aq)^\flat + \kappa \co,\smf\bigr)
= (\kappa +\nuf)\co.
\end{equation*}
Therefore
\begin{equation*}
\begin{split}
(\calb_{12} s_1, \smf)&=\sqrt{\frac{\varepsilon}{\kappa}}(\calb_{12}\qf_1 ,\smf) +(\kappa +
\nuf)\co\\
&=\frac{\varepsilon}{\kappa}\bigl((-\theta + r^2 - 2 r q + 4 r^3q)q_1,\sigma\bigr)^\flat
+(\kappa + \nuf)\co.
\end{split}
\end{equation*}
On the other hand,
\begin{align*}
(\calc_{11}\smf, \smf)&= \bigl((5 r^4 - 3r^2)\sigma, \sigma\bigr)^\flat
+\frac{1}{\kappa}\bigl([5 \xi^4 - 3 \xi^2 -(5 r^4 - 3 r^2)^\flat]\sigf,\sigf\bigr)\\
&\quad +\frac{\varepsilon}
{\kappa}\bigl((2 rq - q^2 + 6 r^2 q^2)\sigma, \sigma\bigr)^\flat +
\kappa \co.
\end{align*}
But according to~\eqref{4.21}, a Taylor expansion gives
\begin{equation*}
5 \xi^4 - 3 \xi^2 -(5 r^4 - 3 r^2)^\flat=(20 r^3 - 6 r)^\flat (\xi -\rf) + \kappa^2
\co
\end{equation*}
and thanks to Lemma~\ref{Lemma4.4}, this expression is equal to  
$-(20 r^3 - 6r) \xi_1 + \kappa^2
\co$. 
Hence, from~\eqref{2.43}
\begin{equation}
\begin{split}
\bigl([5 \xi^4 - 3 \xi^2 -(5 r^4 - 3 r^2)^\flat]\sigf,
\sigf\bigr)
&=(\Af\sigf_2,\xi_1) + \kappa\nuf\co \\
&= \bigl(\sigf_2, G_1(\Uf,\nu)\bigr)+\kappa\nuf\co\\
&=\bigl(\sigf_2, \kappa(r^5 - r^3)^\flat+\varepsilon(r^2 q + 2 r^3 q - rq^2 -\theta
q)^\flat\bigr)\\ &\quad\quad +\kappa (\nuf+\kappa)\co.
\end{split}\label{eq:51}
\end{equation}
Now we obtain the expansion
\begin{align*}
&(\calc_{11}\smf, \smf) + (\calb_{12}s_1, \smf)\\
&=
\bigl((5r^4 - 3r^2)\sigma, \sigma\bigr)^\flat +\frac{\varepsilon}{\kappa}\bigl((-\theta + r^2
- 2 rq + 4r^3q)q_1, \sigma\bigr)^\flat 
 +(r^5 - r^3,\sigma_2)^\flat\\
& +\frac{\varepsilon}{\kappa}(r^2q + 2 r^3 q^2 - rq^2 -\theta q,
\sigma_2)^\flat +\frac{\varepsilon} 
{\kappa}\bigl((2rq - q^2 + 6r^2q^2)\sigma, \sigma)^\flat   +(\kappa
+\nuf)\co.
\end{align*}
We observe that
\begin{align*}
\bigl(&(5r^4-3r^2)\sigma,\sigma\bigr)+(\sigma_2,r^5
- r^3) \\
&\quad=2\int_0^\infty \bigl(\sigma_2(r^5 - r^3)+\sigma^2(5r^4 -
3r^2)\bigr)\,dx\\
&=-2\int_0^\infty (r^5 - r^3)(r_{xx}+2\sigma_x)\,dx +
2\int_0^\infty r_x^2(5r^4 - 3r^2)\,dx +\sqrt\nu\co\\
&=\sqrt\nu \,\co;
\end{align*}
therefore, with the help of~\eqref{5.38}, we are left with
\begin{equation*}
\vabs{\smf}^2 M_{11}^2 =
\frac{\varepsilon}{\kappa}\left[\frac{\partial}{\partial L}(r^2 q +
2r^3 q^2 - r q^2 - \theta  
q, \sigma)\right]^\flat +
 (\sqrt {\kappa}
+ \sqrt{\nf})\co.
\end{equation*}
This concludes the proof of the theorem.
\end{proof}

Therefore, the only remaining question is to obtain an asymptotic for
$\partial a_1/\partial L$. 

\begin{theorem}\label{Theorem5.7} The following asymptotics holds:
\begin{equation*}
\frac{\partial a}{\partial L} = \frac{3\pi^2}{
32 L^2}+O\bigl(L^{-5/2}\bigr).
\end{equation*}
\end{theorem}

\begin{proof}We recall~\eqref{4.12}; we expect that the dominant term
in $a$ will be 
\begin{equation*}
b=m\frac{\partial}{\partial L} (r\sigma \phi', \phi'),
\end{equation*}
the remaining terms being small with respect to $b$, as we will check
later. The number $b$ is also equal to
\begin{equation*}
b= m \int\bigl(\sigma^2 + r\sigma_2) \vabs{ \phi'}^2\, dx + 2m \int
r\sigma \phi' \phi'_1\, dx.
\end{equation*}
Thanks to estimate~\eqref{2.27},
\begin{equation*}
\int\bigl(\sigma^2 + r\sigma_2) \vabs{ \phi'}^2\, dx = 2\int_0^\infty
\bigl[\sigma^2 + r\bigl(-2\sigma ' -r''\bigr)\bigr]\vabs{\phi'}^2 \,
dx + O\bigl(\sqrt\nu\,\bigr),
\end{equation*}
and by integration by parts:
\begin{equation*}
\int\bigl(\sigma^2 + r\sigma_2) \vabs{ \phi'}^2\, dx=
2\int_0^\infty\bigl[\bigl(\sigma + r'\bigr)^2\vabs{\phi'}^2 +
2r(2\sigma +r') \phi'\phi''_bigr]\,dx + O\bigl(\sqrt\nu\bigr).
\end{equation*}
Now, thanks to~\eqref{2.26}, we obtain
\begin{equation*}
b= 4m \int_0^\infty r\sigma \phi'\bigl(\phi'' + \phi'_1\bigr)\, dx +
\sqrt{\nu} \co.
\end{equation*}

We differentiate~\eqref{3.20} with respect to $L$, and we recall the
definition~\eqref{eq:33} of $S$; then
\begin{equation*}
\begin{split}
&m\bigl(\phi'' +\phi'_1\bigr) = -\frac{\theta_1}{R }
\int_x^\infty R\, dy \\
&-\phi'\frac{\partial m}{\partial L} +
\frac{1}{R^2}(S + R') \int_x^\infty \bigl(\theta -R)R\, dy
-\frac{1}{R(x)}\int_x^\infty (\theta -2R)(S +R')\, dx
\end{split}
\end{equation*}
In the above expression, the principal term is the first term on the
right hand side;
the other terms are estimated as follows: it is immediate that
\begin{equation*}
\phi' \partial m/\partial L = \nu O(1);
\end{equation*}
thanks to~\eqref{eq:34} and~\eqref{eq:36}, 
\begin{equation*}
\frac{S+R'}{R^2} \int_x^\infty (\theta-R)R\, dx = O\bigl(\sqrt \nu\,\bigr);
\end{equation*}
thanks to~\eqref{eq:34},~\eqref{eq:35} and~\eqref{eq:37}, we also have
\begin{equation*}
\vabs{\frac{1}{R(x)}\int_x^\infty (\theta -2R)(S +R')\, dx} \le
\frac{O(1)\nu}{\tilde R(x-L)}\int_{x-L}^\infty \frac{\tilde
R(y)}{1-4\tilde R(y)/3}\, dy = O(\nu);
\end{equation*}
therefore, 
\begin{equation*}
m\bigl(\phi''+\phi'_1\bigr)= -\frac{\theta_1}{R }
\int_x^\infty R\, dy + O\bigl(\sqrt\nu\, \bigr).
\end{equation*}

Thanks to~\eqref{2.12} and~\eqref{2.26},
\begin{equation*}
\sigma r = -\frac12 \tilde R'(\cdot-L) +O(1)\nu^{3/4} e^{-2(x-L)};
\end{equation*}
it is straightforward that
\begin{align*}
&\int_0^\infty \nu^{3/4} e^{-(2x-2L)} \frac{\ln (1-4\tilde
R(x-L)/3)}{\tilde R(x-L)}\, dx = O\bigl(\nu^{1/4}L\bigr)\\
\intertext{and}
&\int_0^\infty O\bigl(\sqrt\nu\, \bigr) \phi' r\sigma \, dx =
O\bigl(\sqrt\nu\, \bigr). 
\end{align*}

Thus, we have proved that
\begin{equation*}
b= -\frac{3\theta_1}{4}\int_0^\infty \tilde R'(\cdot
-L)\frac{\ln(1-4\tilde R(\cdot -L)/3)}{\tilde R(\cdot -L)}\, \phi'\,
dx + O\bigl(L\nu^{1/4}\bigr).
\end{equation*}

In order to give an asymptotic of the integral in the above
expression, we cut it into three pieces: one piece from $0$ to
$L-\sqrt L$, which is
\begin{equation*}
\int_0^{L-\sqrt L}\frac{\ln(1-4\tilde R(\cdot -L)/3)}{\tilde R(\cdot
-L)}\tilde R'(\cdot -L) O(1)\, dx=L^2 O(1)e^{-2\sqrt L};
\end{equation*}
the second piece is
\begin{equation*}
2\int_{L-\sqrt L}^\infty \frac{\ln(1-4\tilde R(\cdot -L)/3)}{\tilde R(\cdot
-L)}\tilde R'(\cdot -L) \, dx
\end{equation*}
which we integrate thanks to the change of variable $y=4\tilde R(x-L)/3$:
it is thus equal to 
\begin{equation*}
2\int_0^{\tilde R(-\sqrt L)} \frac{\ln(1-y)}{y}\, dy= -\frac{\pi^2}{3}
+ \sqrt LO\bigl(e^{-\sqrt L}\bigr).
\end{equation*}
The last piece is
\begin{equation*}
\int_{L-\sqrt L} (\phi'(x)+2) \frac{\ln(1-4\tilde R(\cdot -L)/3)}{\tilde R(\cdot
-L)}\tilde R'(\cdot -L) \, dx;
\end{equation*}
since for $x\ge L-\sqrt L$, $\phi'+2=O\bigl(1/\sqrt L\,\bigr)$, this
last piece is an $O(1/\sqrt L)$.
Finally, we have obtained
\begin{equation*}
b=\frac{3\pi^2}{32L^2} + O\bigl(L^{-5/2}\bigr).
\end{equation*}
There remains to estimate the other terms in $\partial a/\partial
L$. It is easy to see that all of them are of order $\nu \co$, and
therefore negligible before the error given in the above formula. This
completes the proof of the Theorem.
\end{proof}

We are now in grade to state the stabilization property:

\begin{theorem}\label{Theorem5.8}
For $0\le y \le L^p$, the following asymptotic
holds:
\begin{equation}
M_{11}\bigl(\sqrt\kappa\bigr) = -\frac{3\nuf}{2} +\frac{\varepsilon\pi^2}{4(L +
y )^2} + \kappa O(L^{-5/2}) 
+ (\nf)^{3/2}\co.\label{5.41}
\end{equation}
There is a critical $y_c$ such that
$a(y_c,\nu)$ vanishes; an equivalent for $y_c$ is given by
\begin{equation*}
y_c \sim \frac{1}{2}\ln L.
\end{equation*}
The corresponding critical parameter is
\begin{equation}
\alpha_c(\nu) = \frac{1}{2}\sqrt\nu\left(1 - \frac{\pi^2}{48 L^2}+
O\bigl(L^{-5/2}\bigr)\right) .\label{5.42} 
\end{equation}
\end{theorem}

\begin{proof}
The asymptotic~\eqref{5.41} is a consequence of all the previous results.
In order to find the critical value of $y$, we have to solve the
equation
\begin{equation*}
-\frac{3}{2}\nf +\frac{\varepsilon\pi^2}{4(L+y)^2} +\bigl(\nf\bigr)^{3/2}\co +
\kappa O\bigl(L^{-5/2}\bigr)=0. 
\end{equation*}
After replacing $\kappa$ and $\nf$ by their respective values, and
defining the new unknown
\begin{equation*}
Y=e^{-4Y},
\end{equation*}
we find the following equation in $Y$:
\begin{equation*}
- Y +\frac{\pi^2(1-Y)}{24(L-\ln Y /4)^2} +O\bigl(L^{-5/2}\bigr)=0
\end{equation*}
It is clear that there is a solution $Y_c$ of this equation which
satisfies
\begin{equation*}
Y_c=\frac{\pi^2}{24 L^2} + O\bigl(L^{-5/2}\bigr).
\end{equation*}

Using relations~\eqref{4.4}, \eqref{4.14},~\eqref{4.21} and ~\eqref{5.41}, we can
see now that at the critical value of $y$,
\begin{equation*}
\varepsilon =\frac{\nu}{4} \left(1-\frac{\pi^2}{48 L^2}
+O\bigl(L^{-5/2}\bigr)\right). 
\end{equation*}
Relation~\eqref{5.42} is an immediate consequence of the above,
together with $\alpha=\sqrt\varepsilon$.
\end{proof}

The final result of this section is

\begin{proposition}\label{Proposition5.9}
Let $u$ be the solution defined at theorem 4.5. Let
$\vabs{\alpha}\le \alpha_m$ be the number defined at~\eqref{eq:39} and
$\alpha_c$ be the number defined at~\eqref{5.42}. Then, the
solution $u$ is stable iff $\alpha_c < \abs{\alpha} <\alpha_m$, and
unstable if $\abs{\alpha}\le \alpha_c$. More precisely, in the first
case, the spectrum of the linearized operator at $u$ contains exactly
one negative eigenvalue, if $\abs{\alpha}<\alpha_c$; when
$_alpha=\alpha_c$, the eigenvalue $0$ is of algebraic multiplicity $3$
and geometric multiplicity $2$, with a non trivial Jordan block of
dimension $2$; in the second case, the linearized
operator at $u$ has the semisimple double eigenvalue $0$, and the
remainder of the spectrum is included in $\Re z \ge \mu >0$.
\end{proposition}

\begin{proof}
This Proposition summarizes our previous analysis: the part of the
spectrum of $\cald$ 
which is outside of the disk of radius $L^{-2p-1}$ around $0$ is
contained in the half plane 
\begin{equation*}
\Re \zeta \ge C L^{-2p},
\end{equation*}
for $\nu$ small enough, according to Theorem~\ref{Theorem2.8},
Lemma~\ref{Lemma3.4} and the definition~\eqref{5.7} of $\cald$.
Therefore, the change of 
stability is equivalent to the change of sign of $M_{11}(\sqrt\kappa)$ given
by~\eqref{5.41}. Theorem~\ref{Theorem5.8}
implies that $M_{11}$ vanishes for some value $\alpha_c$ of $|\alpha|$,
whose asymptotic is 
given by~\eqref{5.42}. The results on the multiplicity of $0$ will be
clear provided that we 
show that $M_{21}$ does not vanish; but formula~\eqref{eq:40} together
with~\eqref{3.14} proves that this is the case.
This concludes the proof of our proposition.
\end{proof}

The last step is to understand the evolution of small perturbations of
the pulse; the 
linearized evolution is easy to understand: if $\delta_i(t)$ are the coordinates 
in the basis $\{z, iu, u'\}$ of the
projection of the linearized perturbation on the space $\im Q$, then
\begin{equation*}
\begin{pmatrix}
\delta_1(t)\\\delta_2(t)\\\delta_3(t)
\end{pmatrix}
 =
e^{i\omega t}\begin{pmatrix}
e^{-at} & 0 & 0\\ b(e^{-at}-1)/a & 1 & 0\\ 0 & 0 &1
\end{pmatrix}
\begin{pmatrix}
\delta_1(0)\\\delta_2(0)\\\delta_3(0)
\end{pmatrix}
.\label{5.43}.
\end{equation*}
For $a>0$, we can see that the asymptotic behavior of the perturbation is simply
described by a change of phase. The smaller $a$, the larger this change of phase, and the
smaller the convergence to the phase shift.

This linear theory can be justified, using \cite{He},  Chapter 5, Exercise 6. 
Let us denote by $u(x,t)$ the stable pulse that we have found, and by
$\tilde u(x,t)$ the 
perturbed solution of~\eqref{1.5}.
It is clear that $\calb$ is a sectorial operator, because it is a bounded perturbation
of a positive self-adjoint operator.  Then there exists
in the stable case a $\delta>0$ such that if $|\tilde u_0 - u(\cdot,0)|_{H^2(\Er)}\le
\delta$, then, there exist $X(\tilde u_0)$ and $\Psi(\tilde u_0)$ such that
\begin{equation*}
\lim_{t\to \infty} |\tilde u(x,t) - u(x - X(\tilde u_0), t)e^{i\Psi(\tilde
u_0)}|_{H^2(\Er)} =0.
\end{equation*}
Moreover, a straightforward calculation shows that~\eqref{5.43} gives
the correct behavior for 
small $\delta$.

\section{The general situation}\label{sec:general-situation}

It is indeed possible to perform the same kind of computations for the
general case, i.e. equation~\eqref{eq:46}, under the assumption
\begin{equation}
\Im m_j =\sqrt\ve\, \mu_j, \quad j=1, 2, 3.\label{eq:47}
\end{equation}

let us give the main steps of the calculation; the functional
validation is entirely identical, and we will shorten the exposition
by giving only the relevant part of the expansion, and without
bothering to give estimates on the remainders.

Under assumption~\eqref{eq:47}, the functional $F$ becomes
\begin{equation*}
F(u)=i\omega u -mu''+\bigl(m+i\sqrt\ve\, \mu_1\bigr)u
-\bigl(1+i\sqrt\ve\, \mu_2\bigr)\abs{u}^2u
+\bigl(1+i\sqrt\ve\,\mu_3\bigr)\vabs{u}^4u. 
\end{equation*}

The functions $G_1$ and $G_2$, which are analogous to the ones defined
at the beginning of section~\ref{sec:second-equation} are given by
\begin{align*}
G_1(\xi,\eta,\tau,\ve, \nu)&=-\ve\tau\eta -m\xi''+m\xi -\mu_1\ve\eta
-\bigl(\xi^2+\ve\eta^2\bigr)\bigl(\xi-\mu_2\ve\eta\bigl)\\
&\qquad+\bigl(\xi^2+\ve\eta^2\bigr)^2\bigl(\xi-\mu_3\ve\eta\bigr),\\
G_2(\xi,\eta,\tau,\ve, \nu)&=\tau\xi -m\eta''+m\eta +\mu_1\xi
-\bigl(\mu_2\xi +\eta\bigr)\bigl(\xi^2
+\ve\eta^2\bigr)\\
&\qquad+\bigl(\mu_3\xi +\eta\bigr)\bigl(\xi^2+\ve\eta^2\bigr)^2.
\end{align*}

The second equation~\eqref{3.8} is changed into
\begin{equation*}
Bq=-\theta r -\mu_1r +\mu_2r^3-\mu_3r^5, \quad (q,r)=0,
\end{equation*}
and $\theta$ is determined by the orthogonality condition
\begin{equation*}
\int\bigl((\theta +\mu_1)R -\mu_2R^2 +\mu_3 R^3\bigr).
\end{equation*}
For $\phi$ given by~\eqref{3.18}, we have the relations
\begin{equation*}
-m\bigl(r^2\phi'\bigr)'=-\theta r^2 -\mu_1r^2 +\mu_2r^4 -\mu_3r^6
\end{equation*}
and
\begin{equation*}
m\phi'(x)=\frac{1}{R(x)} \int_x^\infty \bigl[-\theta R -\mu_1 R
+\mu_2R^2 -\mu_3R^3\bigr]\, dy.
\end{equation*}
The derivative $q_1$ of $q$ with respect to $L$ satisfies the relation
\begin{equation*}
Bq_1 =-\theta_1 r -\theta \sigma -\mu_1\sigma + 3\mu_2 r^2\sigma
-5\mu_3 r^4\sigma +2rq\sigma -4r^3q\sigma + \nu\co.
\end{equation*}

The asymptotics for $\theta$, $\theta_1$ and $\phi'$ are given by
\begin{align*}
\theta&=-\mu_1 +\frac34 \mu_2 -\frac{9}{16} \mu_3
-\frac{3\mu_2}{8L} +\frac{27\mu_3}{64L} +\sqrt\nu\, O(1),\\
\theta_1&=\frac{3\mu_2}{8L^2} -\frac{27\mu_3}{64 L^2} +\sqrt\nu\, O(1)
\end{align*}
and
\begin{equation}
\phi'(x)\sim 
\frac{c_1}{L}\bigl(1+e^{2(x-L)}\bigr)
\ln\bigl(1+e^{2(L-x)}\bigr)
+c_2 +\frac{3\mu_3}{4} \bigl(1+e^{2(x-L)}\bigr)^{-1}
\label{eq:48}
\end{equation}
where the numbers $c_1$, $c_2$ and $c_3$ are given by:
\begin{equation}
c_1=\mu_2-\frac{9\mu_3}{8},\quad c_2=-2\mu_2
+\frac{3\mu_3}{2}, \quad c_3=\frac{3\mu_3}{4}.\label{eq:52}
\end{equation}

We obtain these asymptotics with the help of the following
calculations:
\begin{align*}
\int_x^\infty \tilde R(y)\, dy &= \frac38 \ln\bigl(1+e^{-2x}\bigr),\\
\int_x^\infty \tilde R(y)^2\, dy
&=\frac{9}{32}\bigl[\ln\bigl(1+e^{-2x}\bigr)-\bigl(1+e^{2x}\bigr)^{-1}\bigr],\\
\int_x^\infty \tilde R(y)^3\,
dy&=\frac{27}{128}\bigl[\ln\bigl(1+e^{-2x}\bigr)
-\bigl(1+e^{2x}\bigr)^{-1}
-\bigl(1+e^{2x}\bigr)^{-2}/2\bigr].
\end{align*}

The existence ansatz is completely analogous to the one defined at the beginning
of section~\ref{sec:Existence-pulses} by~\eqref{4.1},~\eqref{4.2}
and~\eqref{4.3}; the definition of $a_0$ is as
in~\eqref{4.9}, and we let
\begin{align*}
a_1&=\bigl((-\theta-\mu_1+\mu_2r^2-\mu_3r^4)q-q^2r +2r^3
q^2,s\bigr)\\
a_2&=\bigl(\mu_2 q^3 +rq^4 -2\mu_3 r^2q^3,s\bigr),\\
a_3&=-\bigl(\mu_3q^5,s).
\end{align*}

It is straightforward to prove that the new $a_1$ is still equivalent
to $m\bigl(r\sigma\phi',\phi'\bigr)$ as in~\eqref{4.12}, and we just have to
calculate $\phi'(L)$ in order to apply formula~\eqref{eq:44}; we infer
from estimate~\eqref{eq:48} that
\begin{equation*}
a_1\sim \frac{9}{64}\bigl(-2\mu_2 + 15 \mu_3/8\bigr)^2 + O(1/L).
\end{equation*}
Therefore,
\begin{equation}
\ve\f \sim\kappa\bigl(-2\mu_2 + 15 \mu_3/8\bigr)^{-2}.\label{eq:50}
\end{equation}

The remainder of the existence proof is identical to the proof given
in section~\ref{sec:Existence-pulses}, its details are left to the reader.

Next step is to calculate $D_u F$ at the solution obtained by this
existence proof; we will have~\eqref{5.7}, with $\calb$ and $\calc$ of
the same form as in~\eqref{eq:49}, and the new values for the
coefficients of $_calb$ and $\calc$ are given by
\begin{align*}
\calb_{12}&=\frac{\sqrt\ve\bigl[-\tau -\mu_1 -2\xi\eta
+\vabs{u}^2\mu_2 +\eta^2\ve +4\vabs{u}^2\xi\eta-\vabs{u}^4 \mu_3
-4\vabs{u}^2\eta^2 \ve\mu_3\bigr]}{\sqrt\kappa(1-\kappa)}\\
\calb_{21}&=\frac{\sqrt\ve\bigl[\tau+\mu_1 -\vabs{u}^2
-2\xi^2\mu_2-2\xi\eta
+4\xi\eta \vabs{u}^2+\mu_3\abs{u}^4 + 4\mu_3
\abs{u}^2\xi^2\bigr]}{\sqrt\kappa(1-\kappa)},\\
\calc_{11}&=\frac1{\kappa} \Bigl[\frac{-\abs{u}^2 -2\xi^2 +2\ve
\mu_2\xi\eta +\abs{u}^4 +4\abs{u}^2 \xi^2 -4\mu_3
\xi\eta\abs{u}^2\ve}{1-\kappa}
+\bigl(3r^2 -5r^4)\f\Bigr],\\
\calc_{22}&=\frac1{\kappa}\Bigl[\frac{-\abs{u}^2 -2\eta^2\ve
-2\eta\xi\mu_2\ve +\abs{u}^4 +4\abs{u}^2\eta^2\ve
+4\mu_3\abs{u}^2\eta\xi}{1-\kappa} +\bigl(r^2 -r^4\bigr)\f\Bigr].
\end{align*}
The principal part of $\calb_{12}$ is
\begin{equation*}
\sqrt{\frac{\ve}{\kappa}}\bigl(-\theta -\mu_1 -2rq +\mu_2 r^2 +
4r^3q -\mu_3r^4\bigr)\f,
\end{equation*}
and the principal part of $\calb_{21}$ is
\begin{equation*}
\sqrt{\frac{\ve}{\kappa}}\bigl(\theta +\mu_1 -3\mu_2 r^2 -2rq +
4r^3 q+ 5\mu_3 r^4\bigr)\f.
\end{equation*}
The principal part of $\calc_{11}$ is
\begin{equation*}
\bigl(5r^4 -3r^2\bigr)\f +\frac{\ve}{\kappa}\bigl(-r^2 +2\mu_2 rq +
6r^2 q^2 -4r^3q\mu_3\bigr)\f +\frac{5\xi^4 -3\xi^2 -\bigl(5r^4
-3r^2\bigr)\f}{\kappa}. 
\end{equation*}

Lemmas~\ref{Lemma5.1},~\ref{Lemma5.2} and~\ref{Lemma5.3} have an
identical statement and an identical proof. The calculations of
Lemma~\ref{Lemma5.4} are modified only to get $M_{11}^2$: we have
first
\begin{equation*}
\begin{split}
&\bigl(\calc_{11}s\f,s\f\bigr)\sim \bigl(\calc_{11}\sigf,\sigf\bigr)\\
&\sim
\bigl((5r^4-3r^2)\sigma,\sigma\bigr)\f+\frac{\ve}{\kappa}\bigl((-r^2
+2\mu_2 rq+6r^2q^2-4r^3q\mu_3)\sigma,\sigma\bigr)\f\\
&\quad+\frac{\bigl(5\xi^4 -3\xi^2-(5r^4 -3r^2)\f\sigf,\sigf\bigr)}{\kappa}.
\end{split}
\end{equation*}
With the same computation as in~\eqref{eq:51}, the last term of the
principal part is $\bigl(\sigma_2\f,G_1(\Uf,\nu)\bigr)/\kappa$ which is
equal to
\begin{equation*}
\bigl(\sigma_2,r^5
-r^3\bigr)\f+\frac{\ve}{\kappa}\bigl(\sigma_2,-\theta q-\mu_1
q-rq^2 +\mu_2r^3q +2 r^3q^2 -\mu_3 r^4q\bigr)\f,
\end{equation*}
up to higher order terms.
On the other hand
\begin{equation*}
\bigl(\calb_{12}s_1,s\f\bigr) \sim
\sqrt{\frac{\ve}{\kappa}}\bigl(\calb_{12}\qf_1,s\f\bigr)
\sim\frac{\ve}{\kappa}\bigl((-\theta -\mu_1 -2rq +\mu_2 r^2 +
4r^3 q-\mu_3 r^4\bigr)q_1,\sigma\bigr)\f.
\end{equation*}
We find finally the same result as in~\eqref{5.37}.

The computation of $\partial a/\partial L$ proceeds along the 
lines of Theorem~\ref{Theorem5.7}. The principal part of $\partial
a/\partial L$ is
\begin{equation*}
4m \int_0^\infty r\sigma \phi'\bigl(\phi''+\phi'_1\bigr)\, dx
\end{equation*}
and
\begin{equation*}
m\bigl(\phi''+\phi'_1\bigr)\sim -\frac{\theta_1}{R(x) }\int_x^\infty
R\, dy.
\end{equation*}
Therefore, it suffices to find an equivalent of
\begin{equation*}
2\theta_1 \int_0^\infty \frac{\tilde R'(x-L)\phi'(x)}{\tilde
R(x-L)}\int_x^\infty \tilde R(y-L)\, dy\, dx.
\end{equation*}
We replace $\phi'$ by its expansion~\eqref{eq:48} which implies
that the expression that we wish to estimate is
\begin{equation*}
\begin{split}
&-\frac{3\theta_1}{4} \biggl[\frac{c_1}{L}\int_{-L}^\infty
2e^{2x}\bigl[\ln\bigl(1+e^{-2x}\bigr)\bigr]^2\, dx + c_2
\int_{-L}^\infty \frac{2e^{2x}\ln\bigl(1+e^{-2x}\bigr)}{1+e^{2x}}\, dx\\
&\qquad+ c_3\int_{-L}^\infty \frac{2
e^{2x}\ln\bigl(1+e^{-2x}\bigr)}{\bigl(1+e^{2x}\bigr)^2}\, dx\biggr].
\end{split}
\end{equation*}
Each of the three integrals is equal, up to terms of order
$O\bigl(\sqrt\nu\,\bigr)$ to the following values:
\begin{equation*}
\pi^2/3, \quad \pi^2/6, \quad 1;
\end{equation*}
the computation is straightforward, and uses the change of variable
\begin{equation*}
1+e^{2x}=\frac{1}{1-y}.
\end{equation*}
All these calculations lead to the following principal part:
\begin{equation*}
\frac{\partial a}{\partial L} = \frac{1}{L^2}\biggl[\frac{3\mu_2}8
-\frac{27\mu_3}{64}\biggr] \biggl[\frac{\pi^2\mu_2}{4}
-\frac{3\pi^2\mu_3}{16} +\frac{9\mu_3}{16}\biggr] + O\bigl(L^{-5/2}\bigr).
\end{equation*}

Finally, we find that
\begin{equation*}
M_{11}\bigl(\sqrt\kappa\,\bigr)=-\lambda\f+\frac{\kappa\chi(\mu)}{(L+y)^2}
 + (\kappa+\nuf) O\bigl(L^{-5/2}\bigr).
\end{equation*}
where $\chi(\mu)$ is given by
\begin{equation}
\chi(\mu)=\biggl[\mu_2 -\frac{9\mu_3}{8}
\biggr]\biggl[\frac{\pi^2\mu_2}{4} 
-\frac{3\pi^2\mu_3}{16}
+\frac{9\mu_3}{16}\biggr]\frac{1} {\bigl(2\mu_2
-15\mu_3/8\bigr)^2}.\label{eq:53}
\end{equation}
Thus, the skew stabilization takes place if
\begin{equation*}
\chi(\mu)>0
\end{equation*}
The critical value $\nu_c$ for which the stabilization takes place is 
given by
\begin{equation}
\nu_c\sim \nu \frac{2\chi(\mu)}{3L^2}.\label{eq:54}
\end{equation}
We find the value given at proposition~\ref{Proposition5.9}, when
$\mu_2=1$ and the other $\mu_j$'s vanish.

When the coefficient $m$ of $-\partial^2/\partial x^2$ is replaced by
$m+i\sqrt\ve \mu_0$, the skew stabilization also takes place under
the same condition and for the same critical value of $\ve$, up to
higher order terms.

Therefore, we have proved the final result of this article:

\begin{proposition}\label{thr:1}
In the general case i.e. for the equation
\begin{equation*}
u_t =\bigl(m+i\sqrt\ve\, \mu_0\bigr)u_{xx} -\bigl(m+i\sqrt\ve\,\mu_1\bigr)u
+\bigl(1+i\sqrt\ve,\mu_2\bigr) \abs{u}^2 u
-\bigl(1+i\sqrt\ve\mu_3\bigr) \abs{u}^4u,
\end{equation*}
there exists a pulse-like solution for all small enough $\nu$, for all
$p>0$, and for all $y\le L^p$; provided that $\chi(\mu)$ defined
by~\eqref{eq:53} is strictly positive; in this case, 
the skew stabilization takes place, at a
critical parameter $\nu_c$ whose equivalent is given by~\eqref{eq:54}.
\end{proposition}

\section{Appendix. An implicit function Theorem with
estimates}\label{sec:Append-An-impl} 

The purpose of this appendix is to prove a version of the implicit function Theorem 
which is appropriate for singular perturbations.

\begin{lemma}
\label{Lemma7.1} Let $X$ and $Z$ be Banach spaces, and let $f$ be a $C^2$ function from a
neighborhood $\calu$ of $x_0\in X$ to $Z$. Let $z_0 = f(x_0)$.
Assume that $A= Df(x_0)$ has a bounded inverse $A^{-1}$.
Assume that the ball of radius $\rho$ and of center $x_0$ is included in $\calu$. Let
\begin{equation*}
M = \sup_{|\xi|\le \rho}\|A^{-1}D^2f(x_0+\xi)\|.
\end{equation*}
There exist constants $a$ and $K$
given by
\begin{equation*}
a = \min(1,(2\rho M)^{-1}), \quad K = \frac{3 a\rho}{ 4}
\end{equation*}
such that if
$|A^{-1}z_0|\le K$, the equation 
\begin{equation*}
f(x)=0
\end{equation*}
possesses a unique solution in the ball $ \{|x-x_0|\le a\rho\}$; moreover, this solution
satisfies
\begin{equation*}
|x-x_0|\le 2|A^{-1}z_0| \text{ and } |x-x_0
+A^{-1}z_0|\le 2 M|A^{-1}z_0|^2.
\end{equation*}
\end{lemma}

\begin{proof}
Let
\begin{equation*}
F(\xi,t) = \xi - A^{-1}f(x_0+\xi) + (1-t)A^{-1}z_0.
\end{equation*}
It is equivalent to find a fixed point of $F$ for $t=1$ and to find a solution of
$f(x)=0$. Introducing $F$ enables us to start from $t=0$, where we have $F(0,0)=0$.
The implicit function Theorem in its classical formulation would give us the existence of
an interval of $t$ on which we can find a solution of $\xi -F(\xi,t) = 0$, but it would
not
ensure that we have a solution up  to  $t=1$. The purpose of our estimates is to show that
we can go as far as $t=1$.

Let us determine $a\in (0,1]$ such that $\xi\mapsto F(\xi,t)$ is a contraction
of ratio $1/2$ for $|\xi|\le a\rho$ and $0\le t\le 1$:
\begin{align*}
&F(\xi',t)& - F(\xi,t) \\
&= \xi' - \xi - A^{-1}\int_0^1 Df\bigl(x_0 + \xi + t(\xi' -
\xi)\bigr)(\xi'-\xi)\,dt \\
&=\xi' - \xi - A^{-1}\int_0^1 \biggl[Df(x_0)(\xi'-\xi) \\
&\quad + \int_0^1 D^2f\bigl(x_0 + s(\xi +
t(\xi'-\xi))\bigr)  (\xi + t(\xi'-\xi))\otimes(\xi'-\xi)\,ds\biggr]\,dt \\
&=-\int_0^1\int_0^1  D^2f\bigl(x_0 + s(\xi +
t(\xi'-\xi))\bigr)  (\xi + t(\xi'-\xi))\otimes(\xi'-\xi)\,ds\,dt .
\end{align*}
Therefore, we have the inequality
\begin{equation*}
|F(\xi',t) - F(\xi,t)|\le Ma\rho|\xi' - \xi|.
\end{equation*}
Thus, the first condition that we wish to impose is
\begin{equation*}
M a \rho \le \frac{1}{ 2}.
\end{equation*}

The second step is to impose a condition on $a$ and on $|B^{-1} z_0|$ such that $\xi
\mapsto F(\xi,t)$ maps the ball of center $x_0$ and radius $a \rho$ into itself, for
$t\le 1$. We have
\begin{equation*}
\begin{split}
|F(\xi,t)|& \le  \left|\xi - A^{-1}\left[Df(x_0) \xi + \int_0^1 D^2f(x_0+
t\xi)\xi^{\otimes 2}(1-t)\,dt\right] - t A^{-1}z_0\right|  \\
&\le {Ma^2 \rho^2}{2} + |A^{-1}z_0|.
\end{split}
\end{equation*}
Therefore, it is enough to require that
\begin{equation*}
|A^{-1}z_0|\le \frac{3a\rho}{ 4}.
\end{equation*}

Thanks to the strict contraction Theorem, there exists for $t\in
[0,1]$ and for $|x|\le 
a\rho$ a unique solution of $$F(\xi,t) = \xi.$$
Let us denote this solution by $\xi = g(t)$. We may estimate $g(t)$ as
follows: we have the 
inequality
\begin{equation*}
|g(t)| \le |F(g(t),t) - F(0,t)| + |F(0,t)|
\le \frac{1}{2} |g(t)| +  t|A^{-1}z_0|.
\end{equation*}
Therefore,
\begin{equation*}
|g(t) |\le 2t|A^{-1} z_0|.
\end{equation*}
But we can obtain a much better estimate, since
\begin{align*}
g(t) + tA^{-1}z_0 &= F(g(t),t) + t A^{-1}z_0 \\
&=g(t) - A^{-1}f(x_0 + g(t))   +
(1-t) A^{-1} z_0 + t A^{-1} z_0  \\
&=g(t) - A^{-1}z_0 - A^{-1} Df(x_0)g(t) \\
&\quad - \int_0^1 D^2f(x_0 + sg(t))g(t)^{\otimes
2}(1-s)\,ds + A^{-1} z_0  \\
&= - \int_0^1 D^2f(x_0 + sg(t))g(t)^{\otimes  2}(1-s)\,ds .
\end{align*}
Therefore,
\begin{equation*}
|g(t) + tA^{-1} z_0|\le \frac{M}{2}|g(t)|^2 \le 2 M
|A^{-1}z_0|^2,
\end{equation*}
which concludes the proof.
\end{proof}

\bibliography{pulse}
\bibliographystyle{plain}

\end{document}